\documentclass[11pt]{article}
\usepackage{amsfonts}
\usepackage{mathrsfs}
\usepackage{color}
\usepackage{amssymb,amsmath,amsthm,mathtools}
\usepackage{ulem}

\usepackage{appendix}
\usepackage{booktabs,tabularx,float}
\usepackage{cite}
\usepackage{multirow}
\usepackage{hyperref}
\hypersetup{
	hidelinks,
} 

\usepackage[titles]{tocloft}

 \allowdisplaybreaks

\catcode`!=11
\let\!int\int \def\int{\!int}
\let\!lim\lim \def\lim{\!lim}
\let\!sum\sum \def\sum{\!sum}
\let\!sup\sup \def\sup{\!sup}
\let\!inf\inf \def\inf{\!inf}
\let\!max\max \def\max{\!max}
\let\!min\min \def\min{\!min}

\catcode`!=12

\let\oldsection\section
\renewcommand\section{\setcounter{equation}{0}\oldsection}

\newcommand{\ba}{\begin{aligned}}
\newcommand{\ea}{\end{aligned}}
\newcommand{\be}{\begin{equation}}
\newcommand{\ee}{\end{equation}}
\def\p{\partial}
\def\pf{\it{Proof.}\rm\quad}

\newcommand{\ep}{\varepsilon}

\newcommand{\na}{\nabla}

\def\R{\mathbb{R}}

\newtheorem{thm}{Theorem}[section]
\newtheorem{defn}[thm]{Definition}

\newtheorem{pro}[thm]{Proposition}
\newtheorem{lem}[thm]{Lemma}
\newtheorem{cor}[thm]{\textbf Corollary}
\newtheorem{re}[thm]{Remark}

\newtheorem*{heur}{Heuristic}
\everymath{\displaystyle}

\newcommand{\Z}{\mathbb Z}
\newcommand{\B}{B_\beta}

\newcommand{\A}{\mathcal A}
\newcommand{\Boot}{\mathscr B}
\newcommand{\PP}{\mathbb P}
\newcommand{\dd}{\,\mathrm dt}

\newcommand{\NL}{\mathcal N}

\newcommand{\muEight}{-\dot M_7/M_7}
\newcommand{\sgn}{\operatorname{sgn}}
\newcommand{\Rea}{\operatorname{Re}}

\begin{document}
\title{\bf Stability threshold for 3D Boussinesq equations with rotation near the Couette flow and vertical stratified temperature
 		\thanks{hwting702@163.com (Wenting Huang), luozekailzk@163.com (Zekai Luo), sunying@bnu.edu.cn (Ying Sun), jwu29@nd.edu(Jiahong Wu), xjxu@bnu.edu.cn (Xiaojing Xu)}}
\author{
{Wenting Huang$^1$, \quad Zekai Luo$^1$, \quad Ying Sun$^2$, \quad Jiahong Wu$^3$, \quad    Xiaojing Xu$^1$ } \\
\small   $^1$ School of Mathematical Sciences, Laboratory of Mathematics and Complex Systems, MOE,   \\  \small Beijing Normal University, Beijing, 100875, China.\\  \small   $^2$ School of Mathematical Sciences, Beijing University of Posts and Telecommunications, \\ \small Beijing, 100876, China.
\\
\small $^3$Department of Mathematics, University of Notre Dame, Notre Dame, IN 46556, USA}	
\date{}
\maketitle

\noindent{\bf Abstract.}
The most important features that distinguish fluid flow in the atmosphere and ocean are the effects of rotation and stratification, which are described by the Boussinesq equations with rotation. This paper examines the stability threshold at high Reynolds number $\mathbf{Re}$ for the 3D Boussinesq equations with rotation near Couette flow and a vertical stratified temperature. For the zero modes, the combined effects of rotation and stratification sustain velocity-temperature oscillations with frequencies bounded away from zero. In the stable rotational regime, rotation suppresses the linear growth associated with the lift-up effect, while stratification sustains oscillations in the double zero modes. Although the zero modes dynamics are not fully dispersive, we exploit the structure of the oscillatory phases and cancellations in the equations to obtain a modified energy estimate for the entire zero modes system. Combining this estimate with enhanced dissipation and inviscid damping of the nonzero modes, we establish a bound for the stability threshold, with order $\mathbf{Re}^{-2/3}(1+\log\mathbf{Re})^{-1-a_0}$ for any fixed $a_0\in(0,1/2)$. 




\section{Introduction}\label{HSXSEC1}

\qquad In this paper, we consider the three-dimensional (3D) Boussinesq equations with rotation, which model the motion of incompressible ﬂuids under both the Coriolis force and buoyancy:
\begin{equation}\label{1.1}
 \begin{cases}
\partial_{t} V+V \cdot \nabla V-\nu \Delta V+\nabla P+\beta \mathrm{e}_3 \times V=\Theta \mathbf{g} \mathrm{e}_3, \\
\partial_{t} \Theta +V \cdot \nabla \Theta-\mu \Delta \Theta=0, \\
\nabla \cdot V=0.
\end{cases}
\end{equation}
where $(x_1, y, z) \in \mathbb{T} \times \mathbb{R} \times \mathbb{T}$ (the torus $\mathbb{T}$ is the periodized interval $[0,2\pi]$), $V=\left(V^1, V^2, V^3\right)^{\mathrm{T}}$ denotes the velocity field, $P$ is the pressure, $\Theta$ is the temperature, $\mathrm{e}_3=(0, 0, 1)^{\mathrm{T}}$ is the unit vector in the vertical direction. The term $\beta \mathrm{e}_3 \times V$ represents the Coriolis force, where $ \beta \in \mathbb{R}\backslash\{0\}$ is the rotation speed. $\Theta \mathbf{g} \mathrm{e}_3$ characterizes the buoyancy with the  gravity $\mathbf{g}$. $\nu$ represents the viscosity coefficient, i.e., $\nu=\mathbf{Re}^{-1}$, where $\mathbf{Re}$ is the Reynolds number, and $\mu$ is the thermal diffusivity. For simplicity, we take $\mu=\nu$ and $\mathbf{g}=1$. The derivation of  \eqref{1.1} can be referred to
 \cite{Majda2003, Pedlosky1987}. The system \eqref{1.1} is supplemented with initial data:
\begin{align}\label{invalue}
(V, \Theta)(0, x_1, y, z)=(V_{\mathrm{in}}, \Theta_{\mathrm{in}})(x_1, y, z).
\end{align}

We  study the dynamic stability behavior of solutions to \eqref{1.1}
around the nontrivial stationary solution, which is known as Couette flow with the linearly stratified temperature profile 
\begin{equation}\label{1.2}
V_s(y)=(y, 0, 0)^{\mathrm{T}},\,\,\,\,\Theta_s(z)=1+\alpha^2 z,
\end{equation}
where $\alpha > 0$ represents the Brunt-Väisälä frequency, which indicates the fluid's response strength to displacements in the direction of gravity. This scenario corresponds to a stable stratification where warmer fluid overlays colder fluid, with gravity acting as a restoring force for vertical perturbations and thereby supporting internal gravity waves.
The associated pressure  of this steady state is
\begin{equation}\label{1.4}
P_s(y, z)=z+\frac{1}{2} (\alpha^2z^2-\beta y^2).
\end{equation}

If we  introduce the perturbations
   $u=V-(y, 0, 0)^{\mathrm{T}},  \theta=-\frac{1}{\alpha}\left(\Theta-\Theta_s \right),  p=P-P_s$,
then the perturbed $(u, p, \theta)$ satisfies the following system
\begin{equation}\label{1.6}
\begin{cases}
\partial_{t} u+y \partial_{x_1} u-\nu \Delta u+\nabla p^L+\begin{pmatrix}
(1-\beta) u^2 \\
\beta u^1 \\
\alpha \theta
\end{pmatrix}=-u \cdot \nabla u-\nabla p^{NL}, \\
\partial_{t} \theta+y \partial_{x_1} \theta-\nu \Delta \theta-\alpha u^3=-u \cdot \nabla \theta, \\
\nabla \cdot u=0, \\
(u, \theta)(0, x_1, y, z)=(u_{\mathrm{in}}, \theta_{\mathrm{in}})(x_1, y, z),
\end{cases}
\end{equation}
where the initial data $u_{\mathrm{in}}=V_{\mathrm{in}}-(y, 0,0)^{\mathrm{T}}$, $\theta_{\mathrm{in}}=-\frac{1}{\alpha}\left(\Theta_{\mathrm{in}}-(1+\alpha^2 z)\right)$. 
The pressure $p^{L}$ and $p^{NL}$ are determined by 
\begin{align*}
\Delta p^{L}=-2\partial_{x_1} u^2+\beta\left(\partial_{x_1} u^2 -\partial_{y} u^1\right)-\alpha \partial_{z} \theta, \quad 
\Delta p^{NL}=-\partial_{i}u^j \partial_{j}u^i.
\end{align*} 

The Coriolis force, in its role concerning fluid stability, is influenced by the variation of the parameter $\beta$. The relevant physical quantity for quantifying the effects of rotation is known as the Bradshaw-Richardson number (see \cite{Bradshaw1969,Huang2018}) $$B_{\beta}:= \beta(\beta-1)\in [-\frac{1}{4},\infty),$$ 
which effectively captures and quantifies the stable or unstable nature of dynamics described in equation \eqref{1.6}.

\subsection{Background and Motivation}
\qquad Before elaborating on the background and motivation,   some relevant notations are introduced. The Fourier transform of a function $f= f(x_1, y, z)$  will be denoted as follows.   Given $(k,\eta,l)\in \mathbb{Z}\times \mathbb{R}\times \mathbb{Z}$, define 
\begin{equation}\label{2.1}
	\mathcal{F}(f)=\widehat{f}\left(k, \eta, l\right)=\frac{1}{(2 \pi)^{3}}\int_{x_1 \in \mathbb{T}}\int_{y \in \mathbb{R}} \int_{z \in \mathbb{T}}f\left(x_1, y, z\right)e^{- i  \left(k x_1+ \eta y+ l z\right)}  dx_1dydz.
\end{equation}
Denoting the inverse operation to the Fourier transform by $\mathcal{F}^{-1}$ or $\check{f}$
\begin{equation}\label{2.2}
	\mathcal{F}^{-1}(f)=\check{f} \left(x_1, y, z\right)=\sum_{k \in \mathbb{Z}}\int_{\eta \in \mathbb{R}} \sum_{l \in \mathbb{Z}} f \left(k, \eta, l\right)e^{ i \left(k x_1+ \eta y+ l z\right)}  d\eta.
\end{equation}
For a general Fourier multiplier with symbol  $m(k,\eta,l)$, we write $m(D)f$ to denote $\mathcal{F}^{-1}(m(k,\eta,l)\allowbreak\hat{f})$,  that is
\begin{equation*}
	\widehat{m(D)f}=m(k,\eta,l)\widehat{f}.
\end{equation*}

 We define the projections on the zero modes and the non-zero modes in $x_1$ of a function $f$  as follows:
\begin{equation}\label{HSX333}
	P_{k=0}f=f_{0}=\frac{1}{2 \pi}\int_{\mathbb{T}} f(x_1,y,z)d{x_1}, \quad	P_{\neq}f=f_{\neq}=f-f_{0}. 
\end{equation}
Furthermore, the projection onto the zero modes in $x_1$ and $z$ of a function $f(x_1,y,z)$ is denoted by
\begin{equation*}
	\overline{f}_0=\frac{1}{(2 \pi)^2}\int_{\mathbb{T}} \int_{\mathbb{T}} f(x_1,y,z)d{x_1} dz, \quad \quad \widetilde{f}_0=f_0-\overline{f}_0,
\end{equation*}
where $\overline{f}_{0}$ represents the double zero modes and $\widetilde{f}_{0}$ stands for the simple-zero modes of $f$.

\subsubsection{Known results for the Navier-Stokes equation with/without rotation}
\qquad The stability of plane Couette flow has been extensively studied since the pioneering works of Rayleigh \cite{R1879} and Kelvin \cite{K1887}. It is well established that plane Couette flow remains linearly stable for all Reynolds numbers, as noted in \cite{DR, R1973}. However, at high Reynolds numbers, the flow can exhibit nonlinear instability and transition to turbulence under small but finite perturbations—a phenomenon known as the Sommerfeld paradox. One approach to resolving this paradox involves examining the transition threshold problem, a concept rigorously formulated by Trefethen et al. \cite{T1993}. Bedrossian, Germain, and Masmoudi \cite{MR3612004, MR3974608} developed a comprehensive mathematical framework for this issue that can be articulated as follows:

Given a norm $\|\cdot\|_{X}$, determine an exponent $\gamma = \gamma(X)$ such that
\begin{align*}
	&\left\|  u_{\mathrm{in}}  \right\|_{X} \lesssim \; \mathbf{Re}^{-\gamma } \,\, \Rightarrow  \text{ stability},\\
	&
	\left\|  u_{\mathrm{in}}  \right\|_{X}  \gg  \;\mathbf{Re}^{-\gamma } \,\, \Rightarrow  \text{possible instability},
\end{align*}
where $\gamma$ is called the \textit{transition threshold}. 

Without rotation and temperature, \eqref{1.1} reduces to the Navier--Stokes equations. For 2D case on the domain
$\mathbb{T}\times \mathbb{R}$,    Bedrossian, Masmoudi and Vicol \cite{MR3448924} demonstrated that the threshold  $\gamma\leq 0$ in Gevrey space.  Bedrossian, Vicol and Wang \cite{MR3867637} showed $\gamma\leq 1/2$  in Sobolev space $H^N$ with  $N>1$. Wei and Zhang \cite{WZ2023}  obtained $\gamma\leq 1/3$ in $H^3$. On the bounded domain $\mathbb{T}\times [-1,1]$, the results  depend on the boundary conditions. Under Navier-slip boundary conditions,  Wei and Zhang \cite{WZ2026}  proved that $\gamma\leq {1}/{3}$ in Sobolev class, and Bedrossian, He, Iyer and Wang \cite{BHIW2-2024} showed $\gamma\leq 1/2 $ in an anisotropic Sobolev space, allowing the shear component to be small independently of $\nu$; Under no-slip boundary conditions, $\gamma\leq {1}/{2}$ in Sobolev space \cite{MR4121130}.  Further related results include \cite{AB2025, MR3974608, BHIW1-2024, BHIW2-2024,   DWZ2021, LLZ2025, LMZ2025, MR4176913}. The 3D case is considerably more challenging due to the lift-up effect.  Bedrossian, Germain and Masmoudi \cite{MR3612004} designed the  Fourier multipliers that precisely encode the interplay between the dissipation and possible growth  and   obtained the stability threshold index $\gamma\leq 3/2$ for Sobolev space $H^{\sigma}$ $(\sigma>9/2)$. For a related work also see  \cite{MR4126259, BGM2022} and so on. Later on, Wei and Zhang \cite{MR4373161} showed that the transition threshold $\gamma \leq1$ of 3D Navier-Stokes equations in $H^2$.  Chen, Wei and Zhang \cite{CWZ2024}  demonstrated that the transition threshold index is also $\gamma\leq 1$ in the finite channels with the non-slip boundary condition.

Without thermal stratification (i.e., $\Theta=0$), the system \eqref{1.1}  corresponds to the Navier-Stokes equations with rotation. In contrast to the standard Navier-Stokes equations, the Navier-Stokes equation with rotation  not only
represents a fundamental physical model, but also exhibits stabilizing effects induced by the Coriolis force. In \cite{HSX2024,HSX2024-2}, the  authors studied the stability of the Navier-Stokes equation with rotation near Couette flow on the domain $\mathbb{T} \times \mathbb{R} \times \mathbb{T}$. The analysis considered two distinct parameter regimes: one with significant dispersive effects for $B_{\beta}>0$, and another exhibiting a strong lift-up effect at $B_{\beta}=0$. In the dispersive regime, they established the stability threshold  $\gamma\leq 1$.  Under the additional condition
\begin{align}\label{0in}
	\overline{u_0}(t=0,y):=\int_{\mathbb{T}^2} u_{\mathrm{in}}(x_1,y,z) \mathrm{d}{x_1} \mathrm{d}z=0,
\end{align}
Coti Zelati, Del Zotto and Widmayer \cite{CZDZW2025} obtained
nonlinear dispersive estimates and the thresholds $\gamma\le8/9$
for $B_\beta>0$ and $\gamma\le5/6$ for
$B_\beta\gtrsim\nu^{-1}$. Li--Sun--Wang--Wei--Zhang
\cite{LSWZ2025} removed \eqref{0in} by proving the Strichartz estimate
\[
\left\| e^{\mathcal{A} t} h^l \right\|_{L_t^2 (L_y^{\infty}(\mathbb{R}))} 
+ \nu^{1/2} |l| \left\| e^{\mathcal{A} t} h^l \right\|_{L_t^1 (L_y^{\infty}(\mathbb{R}))} 
\lesssim \nu^{-1/3} |l|^{-\frac{1}{6}} \| h^l \|_{L^2(\mathbb{R})},
\]
where $\mathcal A=\sqrt{B_\beta}\mathcal R_3^l+\nu(\partial_y^2-l^2)$
and $\mathcal R_3^l=-il(l^2-\partial_y^2)^{-1/2}$,
and obtained $\gamma\le5/6$ in anisotropic Sobolev spaces on
$\mathbb T\times\mathbb R\times\mathbb T$. They also obtained a
threshold $2/3+\delta$, $\delta>0$, on
$\mathbb T\times\mathbb R^2$.

\subsubsection{Vertical stratification}\label{s1.1.2}
\qquad Coti Zelati and Del Zotto \cite{Zelati2023,ZZW2024}
studied stability near the horizontally stratified Couette state
$V_s=(y,0,0)$, $\Theta_s=1+\alpha^2y$. Under the vanish condition for the initial data of double zero modes, Coti Zelati, Del Zotto and Widmayer \cite{ZZW2024} obtained
$\gamma\le11/12$ for $\alpha>1/2$ and $\gamma\le8/9$ for
$\alpha\gtrsim\nu^{-1/3}$. For the more relevant results of the Boussinesq equation, see \cite{Arbon2025, BBCD2023, CWW2025, DWZ2021, MSZ2022, NZ2024, ZhaiZ2023, ZZ2023} for instance.

In geophysical flows, rotation almost always coexists with vertical thermal stratification. The interplay between rotation and vertical thermal stratification establishes a fundamental geophysical framework that underpins numerous atmospheric, oceanic, and astrophysical flows. Oxley and Kerswell~\cite{OK2024} studied the linear stability of stratified, rotating, viscous plane Couette--Poiseuille flow. When the rotational effect is incorporated, the dispersion relation is reshaped and the lift-up effect in zero modes is suppressed. Notably, the wave branches maintain a strictly positive lower bound on frequency, which holds all the way up to the  non-dispersive double zero modes.

The dynamics of fluids are significantly influenced by rotation and temperature stratification. Rapid rotation excites inertial waves, while perturbations in a stably stratified state generate internal gravity waves. It is also important to note that all these dispersive waves exhibit the oscillatory effects. This phenomenon has already been applied to investigate the global well-posedness of equation \eqref{1.1}, as referenced in works such as \cite{BMN1999, HZ2023, IMT2017, KMY2012, MS2022, SLF2025}.
\subsection{Summary of main result}
\qquad We first introduce a coordinate transform described in detail in \cite{MR4126259,MR3612004,Zelati2023}, which is a central
tool for our analysis. We define a linear change of variable
\begin{equation}\label{2.5}
x=x_{1}- t y, \quad \quad y=y, \quad \quad z=z.
\end{equation}
In the new coordinates, the corresponding differential operators become
\begin{equation}\label{2.6}
\nabla =(\partial_{x_1},\partial_{y},\partial_{z}) \mapsto \nabla_{L}:=(\partial_{x},\partial_{y}^{L},\partial_{z}), \quad\quad \Delta\mapsto \Delta_{L}
\end{equation}
where 
\begin{equation}\label{2.77}
\partial_{y}^{L}:=\partial_{y}-t\partial_{x},\quad  \Delta_{L}:=\partial_{x}^{2}+(\partial_{y}^{L})^{2}+\partial_{z}^{2}.
\end{equation}
Furthermore, we define 
\[  
\nabla_{L,h}:=(\partial_{x},\partial_{y}-t\partial_{x}),\, \quad
\Delta_{L,h}:=(\p_x)^2+(\p_y-t\p_x)^2.
\]
The corresponding Fourier lengths are
\begin{equation}\label{HSX444}
	|\widehat{\nabla_{L}}|=|k,\eta-kt,l|,\quad		|\widehat{\nabla_{L,h}}|=|k,\eta-kt|,
\end{equation}
and 
\begin{equation}\label{2.8}
	p(t, k, \eta, l):=\widehat{-\Delta_{L}}=k^2+(\eta- kt)^2+l^2,\quad
	p_{h}(t, k, \eta):= \widehat{-\Delta_{L,h}}=k^2+(\eta- kt)^2.
\end{equation}
Both symbols have derivative
\begin{equation}\label{2.9}
	\partial_{t}	{p}= \partial_{t}{p_{h}} =-2 k (\eta- kt).
\end{equation}

Combining the definition of the new coordinates as in (\ref{2.5}) and denoting $U^i(t, x, y, z)=u^i(t, x_1, y, z)$ $(i=1, 2, 3)$, $  \Theta(t, x, y, z)=\theta(t, x_1, y, z)$, \eqref{1.6}  can be  rewritten as

\begin{align}\label{2.11}
\begin{cases}
\partial_{t} U-\nu \Delta_{L} U+\begin{pmatrix}
\left(1-\beta \right) U^2 \\
\beta U^1 \\
\alpha\Theta\end{pmatrix}+\left( \beta-2 \right) \nabla_{L} \Delta_{L}^{-1} \partial_{x} U^2-\beta \nabla_{L} \Delta_{L}^{-1} \partial_{y}^{L} U^1\\
\quad-\alpha\nabla_{L}\Delta_{L}^{-1}\partial_{z}\Theta=-U\cdot \nabla_{L} U+\nabla_{L} \Delta_{L}^{-1}(\partial_{i}^{L} U^j \partial_{j}^{L} U^i),\\
\partial_{t}\Theta-\nu\Delta_{L}\Theta-\alpha U^{3}=-U\cdot \nabla_{L}\Theta,\\
U(t=0)=U_{\text{in}},\quad \Theta(t=0)= \Theta_{\text{in}}.
\end{cases}
\end{align}

It is important to note that under the coordinate transformation \eqref{2.5}, the initial value conditions presented in \eqref{2.11} are fundamentally equivalent to the initial values stated in \eqref{1.6}. Specifically, we have $U_{\text{in}}=u_{\text{in}}$ and $\Theta_{\text{in}}=\theta_{\text{in}}$. In the following section, we will not make a special distinction regarding this matter. We will examine the behavior of solutions derived from the linearization of \eqref{1.6}. Our first main result gives the linear stability of \eqref{1.1} near the Couette flow and vertical stratification as follows. 

\begin{thm}[\textbf{Linear stability}]\label{thm:linear}
Assume   that $s\geq0, 1\geq\nu>0$, $\alpha>0$, $\beta \in \mathbb{R}$ with $B_{\beta}>0$. Let $u_{\mathrm{in}}$ and $\theta_{\mathrm{in}}$ be smooth initial data satisfying $\na\cdot u_{\mathrm{in}}=0$. Then the solutions $u=u_{\neq}+\widetilde{u}_0+\overline{u}_0$ and $\theta=\theta_{\neq}+\widetilde{\theta}_0+\overline{\theta}_0$ to the linearized  equation for \eqref{1.6}  satisfy the following  linear stability estimates:

	$(1)$  The enhanced dissipation  and inviscid damping of the non-zero  modes $(U_{\neq}, \Theta_{\neq})$:
\begin{align}
 \left\| U^1_{\neq} (t) \right\|_{H^{s}}  &\lesssim\max \left\{c_{\beta},c_{\beta}^{-1} \right\} e^{-\frac{1}{16}  \nu^{\frac{1}{3}} t} \left(\left\|  u_{\mathrm{in}}\right\|_{H^{s+2}}+\left\|  \theta_{\mathrm{in}}\right\|_{H^{s+2}} \right) , \label{eq:linear-decay}  \\
  {\left\langle {t} \right\rangle} \left\|  U^2_{\neq} (t) \right\|_{H^{s}} &\lesssim\max \left\{c_{\beta},c_{\beta}^{-1} \right\} e^{-\frac{1}{16} \nu^{\frac{1}{3}} t} \left(\left\|  u_{\mathrm{in}}\right\|_{H^{s+3}}+\left\|  \theta_{\mathrm{in}}\right\|_{H^{s+3}} \right) , \label{eq:linear-damping}  \\
  \left\| (U^3_{\neq}, \Theta_{\neq}) (t) \right\|_{H^s} &\lesssim   e^{-\frac{1}{16} \nu^{\frac{1}{3}} t} \left(\left\|  u^3_{\mathrm{in}}\right\|_{H^{s+2}} + c_{\beta} \left\|  \p_xu^2_{\mathrm{in}}-\p_yu^1_{\mathrm{in}}\right\|_{H^{s+1}}+\left\|  \theta_{\mathrm{in}}\right\|_{H^{s+2}} \right),\label{A.33x}  
\end{align}
where $c_\beta=\sqrt{\beta/(\beta-1)}$.

$(2)$ Cancellation of lift-up effect in the Bradshaw-Richardson stable regime: 

	For $k=0$,  if $l \neq 0$, then the simple-zero modes $(\widetilde{u}_{0}, \widetilde{\theta}_{0})$ on the Fourier space are written as
    \begin{align}
    	\widehat{\widetilde{({u}^{1})_0}}(t, \eta, l)=&\left[ \frac{\alpha^2 \eta^2}{(\eta^2+l^2)\omega_0^2} e^{\lambda_1 t}+\frac{\beta (\beta-1) l^2}{2(\eta^2+l^2)\omega_0^2} (e^{\lambda_2 t}+e^{\lambda_3 t})\right] \widehat{\widetilde{({u^{1}_{\mathrm{in}}})_0}} \nonumber \\
    	&+\frac{i (\beta-1) }{2\omega_0} (e^{\lambda_3 t}-e^{\lambda_2 t})\widehat{\widetilde{({u^{2}_{\mathrm{in}}})_0}} \nonumber \\
    	&+\left[\frac{\alpha (\beta-1) \eta l}{(\eta^2+l^2)\omega_0^2} e^{\lambda_1 t}-\frac{\alpha (\beta-1) \eta l}{2(\eta^2+l^2)\omega_0^2}(e^{\lambda_2 t}+e^{\lambda_3 t})\right]\widehat{\widetilde{({\theta_{\mathrm{in}}})_0}}, \label{3.2.4} \\
    	\widehat{\widetilde{({u}^{2})_0}}(t, \eta, l)=&\left[ \frac{i \beta l^2 }{2(\eta^2+l^2)\omega_0} (e^{\lambda_2 t}-e^{\lambda_3 t})\right] \widehat{\widetilde{({u^{1}_{\mathrm{in}}})_0}}+\frac{1}{2} (e^{\lambda_2 t}+e^{\lambda_3 t})\widehat{\widetilde{({u^{2}_{\mathrm{in}}})_0}}  \nonumber \\
    	&+\frac{i\alpha \eta l}{2(\eta^2+l^2)\omega_0}(e^{\lambda_3 t}-e^{\lambda_2 t})\widehat{\widetilde{({\theta_{\mathrm{in}}})_0}}, \label{3.2.5}  \\
    	\widehat{\widetilde{({u}^{3})_0}}(t, \eta, l)=&-\frac{\eta}{l}\widehat{\widetilde{({u}^{2})_0}},
    	\label{3.2.6}\\
    	\widehat{\widetilde{(\theta)_0}}(t, \eta, l)=&\left[ \frac{\alpha \beta \eta l}{(\eta^2+l^2)\omega_0^2} e^{\lambda_1 t}-\frac{\alpha \beta \eta l}{2(\eta^2+l^2)\omega_0^2} (e^{\lambda_2 t}+e^{\lambda_3 t})\right] \widehat{\widetilde{({u^{1}_{\mathrm{in}}})_0}} \nonumber\\
    	&+\frac{i \alpha \eta }{2\omega_0 l} (e^{\lambda_2 t}-e^{\lambda_3 t})\widehat{\widetilde{({u^{2}_{\mathrm{in}}})_0}} \nonumber\\
    	&+\left[\frac{\beta (\beta-1) l^2}{(\eta^2+l^2)\omega_0^2} e^{\lambda_1 t}+\frac{\alpha^2 \eta^2}{2(\eta^2+l^2)\omega_0^2}(e^{\lambda_2 t}+e^{\lambda_3 t})\right]\widehat{\widetilde{({\theta_{\mathrm{in}}})_0}},\label{3.2.7} 
    \end{align}
	where $\omega_0=\omega_0(\alpha, \beta, \eta, l)=\sqrt{\frac{\alpha^2 \eta^2+ B_{\beta} l^2}{\eta^2+l^2}}$, $\lambda_1=-\nu (\eta^2+l^2)$, $\lambda_2=-\nu (\eta^2+l^2)+i\, \omega_0$ and $\lambda_3=-\nu (\eta^2+l^2)-i\, \omega_0$. 

	If  $l=0$, then the double zero modes $(\overline{u}_0, \overline{\theta}_0)$ satisfy
	\begin{align}
	\label{thm1.1.1.1}&\overline{u}_0^1(t, y)=e^{\nu \partial_{yy} t} \overline{(u_{\mathrm{in}}^1)_0}, \qquad  
	\overline{u}_0^2(t, y)=0, \\ \label{thm1.1.2}
	&\overline{u}_0^3(t, y)=e^{\nu \partial_{yy} t} \left[ \cos(\alpha t)\overline{(u_{\mathrm{in}}^3)_0}-\sin(\alpha t)\overline{(\theta_{\mathrm{in}})_0}\right],  \\
	&\overline{\theta}_0(t, y)=e^{\nu \partial_{yy} t} \left[ \sin(\alpha t)\overline{(u_{\mathrm{in}}^3)_0}+\cos(\alpha t)\overline{(\theta_{\mathrm{in}})_0}\right].  \label{thm1.1.3}
	\end{align}

\end{thm}

\begin{re}
Linear stability requires only $B_\beta>0$, its nonzero modes proof
is given in appendix~\ref{secA}. In section~\ref{sec3.1}, the stronger condition
$B_\beta>1/4$ makes the cross energy in
\eqref{00} coercive and gives the weighted estimate
used in the nonlinear proof. The linear stability bounds also hold
when $B_\beta=\alpha^2$.
\end{re}

Our second result establishes nonlinear stability in Sobolev spaces.
\begin{thm}[\textbf{Nonlinear stability}]\label{1..1}
Let $0<\nu\le1$, $\alpha>0$, and $\beta\in\R$ satisfy
$B_\beta>1/4$ and $B_\beta\ne\alpha^2$. For fixed $0<a_0<1/2$ and any integer $N\ge4$,
there exist constants $\delta$ and $C$, depending only on
$N,\alpha,\beta,a_0$, such that the following holds. If the initial data satisfy $\na\cdot u_{\mathrm{in}}=0$ and
\begin{equation}\label{1.23}
\begin{aligned}
\|u_{\mathrm{in}}\|_{H^{N+1}}+
\|\theta_{\mathrm{in}}\|_{H^{N+1}}=\ep\le\delta\nu^{2/3}/(1+|\log\nu|)^{1+a_0},
\end{aligned}
\end{equation}
then \eqref{1.6} has a unique global strong solution satisfying:

$(1)$ Global stability:
\begin{equation}\label{eq:main-bounds}
\|(U,\Theta)\|_{L^\infty_t H^N}
+\nu^{1/6}\|(U_{\neq},\Theta_{\neq})\|_{L^2_tH^N}
+\nu^{1/2}\|\na_L(U,\Theta)\|_{L^2_tH^N}\le C\ep.
\end{equation}

$(2)$ Enhanced dissipation and inviscid damping:
\begin{align}\label{eq:physical-decay}
\|(U_{\neq},\Theta_{\neq})(t)\|_{H^N}
+\langle t\rangle\|U^2_{\neq}(t)\|_{H^{N-2}}
\le C\ep e^{-\lambda\nu^{1/3}t},
\end{align}
where $\lambda=\frac{1}{16}\frac{2\sqrt{B_{\beta}}-1}{2\sqrt{B_{\beta}}+1}$.

$(3)$ Higher-order zero modes bounds: for $0\le\sigma\le1$,
\begin{equation}\label{eq:main-high}
\|(U_0,\Theta_0)\|_{L^\infty_t H^{N+\sigma}}
+\nu^{1/2}\|\na(U_0,\Theta_0)\|_{L^2_t H^{N+\sigma}}
\le C\ep.
\end{equation}

\end{thm}

\begin{re}
(1) Theorem~\ref{1..1} requires neither an additional
$L^1$ integrability assumption nor vanishing double zero modes
initial data since our work is based on an energy method.

(2) Theorem \ref{1..1} states the result when the thermal diffusivity $\mu$ and the viscosity $\nu$ are equal. Indeed, with suitable changes to the diffusion weights and estimates, one can expect the stability threshold for $0<\nu,\mu\le1$ with order 
\[
	O\left(\min\left\{\min\{\nu,\mu\}^{\frac23}(1+|\log \min\{\nu,\mu\}|)^{-1-a_0},\, \max\{\nu,\mu\}^{-\frac13}\min\{\nu,\mu\}\right\}\right).
\]

(3) When $B_\beta=\alpha^2$, we have $\omega_0\equiv\alpha$, where the group velocity and the curvature of the dispersion curves vanish. Hence, our nonlinear zero modes argument does not apply directly.
\end{re}

\subsection{The improvement in the nonlinear estimates}
\qquad We now present the strategy for proving theorem~\ref{1..1}.
For the zero modes, we use the two oscillatory branches and the
stationary branch identified in section~\ref{sec3.2}.
For the nonzero modes, we use three good unknowns and the weighted
energy estimate \eqref{3x.0} in lemma~\ref{lem3.4}.

$(1)$ \textit{Zero modes.}
The linear zero modes dynamics consist of two oscillatory branches and one stationary branch. We preserve the transport cancellation in the total energy. Away from resonance, oscillations allow us to absorb the cubic energy flux into a small energy correction. Near resonance, we combine structural cancellations in the energy coefficients with geometric bounds on the relevant frequency sets. These mechanisms enter directly through energy estimates, without relying on pointwise dispersive or Strichartz estimates.

We construct a small cubic correction to the total energy. It absorbs the nonresonant cubic flux
and leaves controllable cubic viscous and quartic nonlinear
remainders. Near resonance, two complementary structures are used.
For interactions involving two waves and a stationary component,
incompressibility condition makes the complete energy coefficient vanish
when the wave vectors are parallel. In the configuration of
opposite $z$-wavenumbers, the phase and the even wave-weight
difference share a factor, whose cancellation yields the remaining
near-resonant bound.
For three-wave interactions, we design a customized energy weight $L_s(\xi)\asymp\langle\xi\rangle^{2s}$ built from the group velocity $|\p_{\eta}\omega_0(\xi)|$. Near space--time resonance,  this weight is invariant across the three interacting waves, so that the weight differences vanish at resonance and are bounded near resonance by differences of the signed derivatives $\sigma_j\p_{\eta}\omega_0(\xi_j)$. The change of variables $u=\Psi$ in the integrals over $|\Psi|\le\rho$ then yields the gain $\rho^{1/2}$ in lemma~\ref{lem:phase-blocks}, while the two critical ratios $\alpha=2\sqrt{B_\beta}$ and $\sqrt{B_\beta}=2\alpha$ are treated directly via frequency-set bounds. Employing the scale-dependent cutoff $\chi(\Psi/(\vartheta R))$ in \eqref{eq:block-cutoffs} leads to the smallness condition
\[
 \frac{\ep}{\vartheta}
 +\frac{\ep\sqrt\vartheta}{\nu}
 +\frac{\ep^2}{\nu\sqrt\vartheta}\ll1.
\]

Keeping $\vartheta$ free throughout the estimates and choosing $\vartheta=K\ep$ with $K$ fixed and sufficiently large reduces the dissipative smallness condition to $\ep^{3/2}\ll\nu$.

The single zero and double zero modes are treated in a
modified-energy estimate. The stratification and rotation effects
simplify the geometric analysis of degenerate frequency
configurations. For $\alpha^2\ne B_\beta$, the resulting coefficient estimates are uniform over the spatial frequency blocks (lemma~\ref{lem:phase-blocks}), including blocks with arbitrarily small lowest frequency. The uniform $H^{N+1}$ zero modes bound controls the coupling to nonzero modes.

$(2)$ \textit{Nonzero modes and their coupling to the zero modes.}
Multiplication couples nonzero modes with different critical times. We use the Fourier multiplier $m_3$ introduced in \cite{LSWZ2025} and renormalize it as $M_7$ to bound the reaction term appearing in the nonzero modes energy estimates. Using the initial normalization $m(0,\xi)=|\xi|^{1/2}$ with the truncated multiplier $m$ and an integrable logarithmic kernel in $M_7$ gives the coefficient $(1+|\log\nu|)^{1+a_0}\nu^{-2/3}$ in the nonzero modes estimate. Physical variables are recovered using $mp^{1/4}\ge1$. For mixed interactions, incompressibility condition places the derivative on the nonzero modes test function, and the multiplier comparison at the same $x$-wavenumber $k$ requires an $H^{N+1}$ bound for the zero modes.
To propagate this highest-order norm without derivative loss, we estimate the zero-mode energy weighted by $M_7(t,0,\p_y)$ up to the switching time $t_*=A_*\nu^{-1/3}(1+|\log\nu|)$. Beyond $t_{*}$ enhanced dissipation renders the Reynolds stress small enough to close the modified energy estimate in $H^{N+1}$.

The rest of the paper is organized as follows. Section~\ref{Preliminaries}
introduces the notation and functional estimates.
Section~\ref{section3} develops the linear framework for $B_\beta>1/4$
and analyzes the zero modes; appendix~\ref{secA} proves the
nonzero modes estimates for the full range $B_\beta>0$. Section~\ref{section4} gives the bootstrap argument for theorem~\ref{1..1}. The nonlinear estimates for zero and nonzero modes are proved in sections~\ref{sec:zero} and \ref{sec:nonzero-nonlinear}, respectively.

\section{Preliminaries}\label{Preliminaries} \label{sec2}

\subsection{Notation}
\qquad Given two quantities $A$ and $D$, we denote $A \lesssim D$, if there exists a constant $C$ such that $A\leqslant CD$ where  $C>0$ depends only on the fixed parameters $\alpha,\beta,N,a_0$, but not on $\delta$, $\nu$ or the time interval.  { Let $\Boot$ be a fixed bootstrap constant and the notation $\lesssim_{\Boot}$ indicates that the constant $C$ depends on $\Boot$, which can be denoted as $C_{\Boot}$.} We similarly denote $A \ll D$ if $A \leqslant cD$ for a small constant $c \in(0, 1)$ to emphasize the small size of the implicit constant. We also denote $A\asymp D$, if there exists a constant $C$ such that $C^{-1}D\leq A\leq CD$. We define $\xi=(k,\eta,l)$ and $\xi_L=(k,\eta-kt,l)$, then $(1+|\xi|^2)^{1/2}$ is represented by $\langle\xi\rangle$ 
and  $\sqrt{k^2+\eta^2+l^2}$ is denoted by $|k,\eta,l|$. For two functions $f$ and $g$ and a norm $\|\cdot\|_{X}$, we write 
$$\|(f,g)\|_{X}=\sqrt{\|f\|_{X}^{2}+\|g\|_{X}^{2}}.$$
Unless specified otherwise, in the rest of this section $f$ and $g$ denote functions from $\mathbb{T} \times \mathbb{R} \times \mathbb{T}$ to $\mathbb{R}^{n}$ for some $n\in \mathbb{N}$.

\subsection{Functional spaces}
\qquad The Schwartz space $\mathscr{S}(\mathbb{R}^{n})$ is the set of smooth functions $f$ on $\R^{n}$ such that 
$$p_{m,\gamma}(f)\triangleq \sup_{x\in\mathbb{R}^{n}}(1+|x|)^m|\partial^\gamma f(x)|<\infty,\quad \forall\ m\in\mathbb{N}_{0},\ \gamma\in\mathbb{N}_{0}^n.$$

The Sobolev space $H^{s} (s\geqslant0) $ is given  by  the norm 
\begin{equation}\label{HSX-A}
	\left\| f \right\|_{H^{s}}=\left\| {\left\langle {D} \right\rangle}^{s} f \right\|_{L^2}=\left\| {\left\langle {k, \eta, l} \right\rangle}^{s} \widehat{f} \right\|_{L^2}.
\end{equation}

Recall that, for $s>\frac{3}{2}$, $H^{s}$ is an algebra. Hence, for any $f,g\in H^{s}$, one has
\begin{align}\label{2.3}
	\left\|f g\right\|_{H^{s}} \lesssim \left\|f \right\|_{H^{s}}\left\|g\right\|_{H^{s}}.
\end{align}
We write the associated inner product, also for complex-valued functions, as
$$ \langle f,g \rangle_{H^{s}}=\int \langle\nabla\rangle^{s}f\cdot \overline{\langle \nabla \rangle^{s} g} dV. $$
For function $f(t,x,y,z)$ of space and time defined on the time interval $(a,b)$, we define the Banach space $L_{t}^{p}(a,b;H^{s})$ for $1\leqslant p\leqslant \infty$ by the norm
$$\|f\|_{L_{t}^{p}(a,b; H^{s})}=\|\|f\|_{H^{s}}\|_{L_{t}^{p}(a,b)}.$$
For simplicity of notation, we usually simply write $\|f\|_{L_{t}^{p}H^{s}}$ since the time-interval of integrating in this work will be the same basically everywhere.

\section{The proof of theorem \ref{thm:linear}}\label{section3}
\qquad  The nonzero-mode estimates in theorem~\ref{thm:linear},
which hold for $B_\beta>0$,
are proved in appendix~\ref{secA}.
In this section, we  first establish the  weighted linear
estimates under the stronger assumption $B_\beta>1/4$
that are used in the nonlinear argument. The linearized equation is
\begin{equation}\label{1.66}
	\begin{cases}
		\partial_{t} u+y \partial_{x_1} u-\nu \Delta u+\begin{pmatrix}
			(1-\beta) u^2 \\
			\beta u^1 \\
			\alpha \theta
	    	\end{pmatrix}+ (\beta-2)\nabla \Delta^{-1}\partial_{x_1}u^{2}\\ \qquad -\beta\nabla\Delta^{-1}\partial_{y}u^{1}-\alpha\nabla\Delta^{-1}\partial_{z}\theta  =0, \\
		\partial_{t} \theta+y \partial_{x_1} \theta-\nu \Delta \theta-\alpha u^3=0, \\
		u(0, x_1, y, z)=u_{\mathrm{in}}, \quad \theta(0, x_1, y, z)=\theta_{\mathrm{in}}.
	\end{cases}
\end{equation}

By the moving frame \eqref{2.5}, we recall the definition of $ U(t,x,y,z)$ and $ \Theta(t,x,y,z) $, and denote $W^3\triangleq\p_xU^2-\p_y^LU^1$. Here, we introduce the new unknowns
\begin{equation}\label{yyy}
Q_{\neq}=-|\na_L|^{3/2}|\na_{L,h}|^{-1}U^3_{\neq},
H_{\neq}=-|\na_L|^{1/2}\Theta_{\neq},\ K_{\neq}=-ic_\beta|\na_L|^{1/2}|\na_{L,h}|^{-1}W^3_{\neq}.
\end{equation}
where $c_\beta=\sqrt{\beta/(\beta-1)}$. After taking the Fourier transform with respect to the spatial variables, for $k\neq 0$, the linearized system of \eqref{1.66} can be written as
\begin{align}\label{Non0L}
	\begin{cases}
		\p_t \widehat{Q_{\neq}} + \nu p \widehat{Q_{\neq}} - \frac{1}{4} \frac{\p_t p}{p}\widehat{Q_{\neq}} - \frac{1}{2} \frac{l^2}{p}\frac{\p_t p_h}{p_h} \widehat{Q_{\neq}} \\
		\qquad\ \ - { \operatorname{sgn}(\beta)}\sqrt{B_{\beta}}\frac{l}{p^{1/2}}\widehat{K_{\neq}} + \alpha \frac{p_h^{1/2}}{p^{1/2}}\widehat{H_{\neq}} + 2\sqrt{\frac{\beta-1}{\beta}}\frac{k^2}{p_h}\frac{l}{p^{1/2}}\widehat{K_{\neq}}  = 0, \\
		\p_t\widehat{K_{\neq}} + \nu p\widehat{K_{\neq}} + \frac{1}{4}\frac{\p_t p}{p}\widehat{K_{\neq}} +\frac{1}{2} \frac{l^2}{p}\frac{\p_t p_h}{p_h} \widehat{K_{\neq}}+ { \operatorname{sgn}(\beta)}\sqrt{B_{\beta}}\frac{l}{p^{1/2}}\widehat{Q_{\neq}} = 0,\\
		\p_t\widehat{H_{\neq}} + \nu p\widehat{H_{\neq}} - \frac{1}{4}\frac{\p_t p}{p} \widehat{H_{\neq}} - \alpha \frac{p_h^{1/2}}{p^{1/2}}\widehat{Q_{\neq}} = 0.
	\end{cases}
\end{align}

\subsection{Dynamics of nonzero modes}\label{sec3.1}
\qquad \textit{The multiplier $m$.}
As observed in \cite{MR3612004}, to capture the associated effects of inviscid damping and enhanced dissipation, it is convenient to use a suitably constructed Fourier multiplier $m$ with symbol $m=m(t,k,\eta,l)$.
Set $\tau_\nu=1000\nu^{-1/3}$ and $m(t,0,\eta,l)=1$.
For $k\ne0$, define $m$ as the positive solution of
\begin{equation}\label{eq:m-formula}
-\frac{\dot m}{m}
=\begin{cases}
\displaystyle\frac{|k(\eta-kt)|}{2p},&|t-\eta/k|\le\tau_\nu,\\
0,&|t-\eta/k|>\tau_\nu,
\end{cases}
\end{equation}
with initial value
\begin{equation}\label{eq:m-initial}
 m(0,k,\eta,l)=(k^2+\eta^2+l^2)^{1/4}.
\end{equation}

The central role of $m$ is to deal with the slowly decaying term
$\partial_x\partial_y^L|\nabla_L|^{-2}\{Q,K,H\}$ in \eqref{Non0L},
balancing the growth that solutions experience in the transition between
the inviscid damping regime $t<\eta/k$ and the dissipative regime
$t\gtrsim\eta/k+\nu^{-1/3}$.

\begin{lem}\label{lem3.1}
$(1)$ The multiplier $m(t,k,\eta,l)$ can be given by the following exact formula:
\begin{itemize}
\item $m(t,k=0,\eta,l)=1$
\item If $k\ne0$, $\eta/k<-\tau_\nu$: $m(t,k,\eta,l)\equiv(k^2+\eta^2+l^2)^{1/4}$
\item If $k\ne0$, $-\tau_\nu\le\eta/k\le0$:
\[
 m=\begin{cases}
 \displaystyle\left(\frac{(k^2+\eta^2+l^2)^2}{p}\right)^{1/4},
       &0\le t\le\eta/k+\tau_\nu,\\
 \displaystyle\left(\frac{(k^2+\eta^2+l^2)^2}{k^2+(k\tau_\nu)^2+l^2}\right)^{1/4},
       &\eta/k+\tau_\nu\le t.
 \end{cases}
\]

\item If $k\ne0$, $0<\eta/k<\tau_\nu$:
\[
 m=\begin{cases}
 p^{1/4},&0\le t\le\eta/k,\\
 \displaystyle\left(\frac{(k^2+l^2)^2}{p}\right)^{1/4},
       &\eta/k\le t\le\eta/k+\tau_\nu,\\
 \displaystyle\left(\frac{(k^2+l^2)^2}{k^2+(k\tau_\nu)^2+l^2}\right)^{1/4},
       &\eta/k+\tau_\nu\le t.
 \end{cases}
\]

\item If $k\ne0$, $\tau_\nu\le\eta/k$:
\[
 m=\begin{cases}
 (k^2+\eta^2+l^2)^{1/4},&0\le t\le\eta/k-\tau_\nu,\\
 \displaystyle\left(\frac{(k^2+\eta^2+l^2)p}{k^2+(k\tau_\nu)^2+l^2}\right)^{1/4},
       &\eta/k-\tau_\nu\le t\le\eta/k,\\
 \displaystyle\left(\frac{(k^2+\eta^2+l^2)(k^2+l^2)^2}{[k^2+(k\tau_\nu)^2+l^2]p}\right)^{1/4},
       &\eta/k\le t\le\eta/k+\tau_\nu,\\
 \displaystyle\left(\frac{(k^2+\eta^2+l^2)(k^2+l^2)^2}{[k^2+(k\tau_\nu)^2+l^2]^2}\right)^{1/4},
       &\eta/k+\tau_\nu\le t.
 \end{cases}
\]

\end{itemize}

$(2)$ In particular, $m(t,k,\eta,l)$ is continuous in time and, for $k\ne0$,
\begin{equation}\label{4.11-1}
 \nu^{\frac{1}{6}}\lesssim m(t,k,\eta,l)\le(k^2+\eta^2+l^2)^{\frac{1}{4}}.
\end{equation}

$(3)$ $m(t,k,\eta,l)$ and the frequency have the following relationship for $k\ne0$:
\begin{equation}\label{4.11-2}
 \max\left\{\frac{|k,l|}{p^{1/4}},\ (1+t)^{-1/2},\
 \frac{p^{1/4}}{1+t},\ \frac{|k,\eta,l|}{(1+t)p^{1/4}}\right\}
 \lesssim m(t,k,\eta,l)\lesssim\left(1+\frac{t}{\tau_\nu}\right)p^{1/4}.
\end{equation}

$(4)$ Moreover, for all $t$, $k\ne0$, $\eta,\eta',l,l'$, the following product estimate holds:
\begin{equation}\label{4.11-3}
 m(t,k,\eta,l)p(t,k,\eta,l)^{1/4}
 \lesssim\langle\eta-\eta',l-l'\rangle
 m(t,k,\eta',l')p(t,k,\eta',l')^{1/4}.
\end{equation}

\end{lem}
\begin{proof}
$(1)$ Integrating \eqref{eq:m-formula} gives
\[
 m(t,k,\eta,l)=(k^2+\eta^2+l^2)^{1/4}
 \exp\left(-\int_0^t\frac{|k(\eta-ks)|}{2p(s,k,\eta,l)}
 \mathbf1_{\{|s-\eta/k|\le\tau_\nu\}}\,\mathrm ds\right),
\]
which gives the above formulas and their agreement at the endpoints.

$(2)$ By $(1)$ and $\dot m\le0$, one has
\[
 \nu^{1/6}\lesssim(1+\tau_\nu^2)^{-1/4}
 \le\frac{|k,l|}{[k^2+(k\tau_\nu)^2+l^2]^{1/4}}
 \le m(t,k,\eta,l)\le m(0,k,\eta,l).
\]

$(3)$ The formulas in $(1)$ give
\begin{align*}
 m^2\sqrt p
 &=\min\left\{\min_{0\le s\le t}p(s,k,\eta,l),\ k^2+(k\tau_\nu)^2+l^2\right\}\\
 &\quad\times\max\left\{1,\frac{\sqrt{p(t,k,\eta,l)}}{\sqrt{k^2+(k\tau_\nu)^2+l^2}}\right\}
 \max\left\{1,\frac{\sqrt{k^2+\eta^2+l^2}}{\sqrt{k^2+(k\tau_\nu)^2+l^2}}\right\}.
\end{align*}
Hence, 
\[
 |k,l|\le\min_{0\le s\le t}\sqrt{p(s,k,\eta,l)}
 \le mp^{1/4}\le\left(1+\frac{t}{\tau_\nu}\right)\min_{0\le s\le t}\sqrt{p(s,k,\eta,l)}.
\]
Moreover, since
\[
 \max\{\sqrt p,|k,\eta,l|\}
 \le\min_{0\le s\le t}\sqrt{p(s,k,\eta,l)}+|k|t
 \le(1+t)\min_{0\le s\le t}\sqrt{p(s,k,\eta,l)}.
\]
Consequently,
\[
 m\ge\frac{\min_{0\le s\le t}\sqrt{p(s,k,\eta,l)}}{p^{1/4}}
 \ge(1+t)^{-1/2}\min_{0\le s\le t}p(s,k,\eta,l)^{1/4}
 \ge(1+t)^{-1/2},
\]
which proves \eqref{4.11-2}.

$(4)$ For every $0\le s\le t$,
\begin{align*}
 |\sqrt{p(s,k,\eta,l)}-\sqrt{p(s,k,\eta',l')}|
 &\le|\eta-\eta',l-l'|,\\
 \left|\sqrt{k^2+(k\tau_\nu)^2+l^2}
       -\sqrt{k^2+(k\tau_\nu)^2+(l')^2}\right|
 &\le|l-l'|.
\end{align*}

Using these inequalities and $\min_{0\le s\le t}\sqrt{p(s,k,\eta,l)}\ge|k,l|\ge1$, we obtain
\begin{align*}
 &\frac{\displaystyle\min\left\{\sqrt{p(t,k,\eta,l)},\
 \min_{0\le s\le t}\sqrt{p(s,k,\eta,l)}
 \max\left\{1,\frac{\sqrt{p(t,k,\eta,l)}}{\sqrt{k^2+(k\tau_\nu)^2+l^2}}\right\}\right\}}{\displaystyle\min\left\{\sqrt{p(t,k,\eta',l')},\
 \min_{0\le s\le t}\sqrt{p(s,k,\eta',l')}
 \max\left\{1,\frac{\sqrt{p(t,k,\eta',l')}}{\sqrt{k^2+(k\tau_\nu)^2+(l')^2}}\right\}\right\}}\\
 &\quad\le\begin{cases}
 2|\eta-\eta',l-l'|,\quad &\text{if }\sqrt{p(t,k,\eta,l)}>2\sqrt{p(t,k,\eta',l')},\\
 4(1+|\eta-\eta',l-l'|),\quad &\text{if }\begin{gathered}\sqrt{p(t,k,\eta,l)}\le2\sqrt{p(t,k,\eta',l')}\text{ and } |\eta-\eta',l-l'|\le\tau_\nu,\end{gathered}\\
 2\langle\eta-\eta',l-l'\rangle,\quad &\text{if }\begin{gathered}\sqrt{p(t,k,\eta,l)}\le2\sqrt{p(t,k,\eta',l')}\text{ and }|\eta-\eta',l-l'|>\tau_\nu.\end{gathered}
 \end{cases}
\end{align*}
Here the three cases use, respectively,
\[
 \sqrt{p(t,k,\eta,l)}<2|\eta-\eta',l-l'|,
 \frac{\sqrt{k^2+(k\tau_\nu)^2+(l')^2}}{\sqrt{k^2+(k\tau_\nu)^2+l^2}}\le2, \frac{\sqrt{k^2+(k\tau_\nu)^2+(l')^2}}{|k,l'|}
 \le\langle\eta-\eta',l-l'\rangle.
\]
By $(1)$, we have
\[
 m^2\sqrt p=\prod_{\zeta=\eta,\,kt-\eta}
 \min\left\{\sqrt{p(t,k,\zeta,l)},\
 \min_{0\le s\le t}\sqrt{p(s,k,\zeta,l)}
 \max\left\{1,\frac{\sqrt{p(t,k,\zeta,l)}}{\sqrt{k^2+(k\tau_\nu)^2+l^2}}\right\}\right\}.
\]
Applying the preceding estimate to $(\eta,\eta')$ and $(kt-\eta,kt-\eta')$ gives \eqref{4.11-3}.
\end{proof}

\noindent\textit{The ghost multipliers $M_1,\ldots,M_6$.}
For $k\ne0$, define the multipliers $M_{i}(t, k, \eta, l), i=1,2,...,6$,  by $M_{i}(t=0, k, \eta, l)=1$ and the following equations. For $k=0$, set $M_j(t,0,\eta,l)=1$.
{\small
\begin{align}
-\frac{\dot M_1}{M_1}
 &=\frac{\nu^{1/3}}{1+\nu^{2/3}(t-\eta/k)^2}, &-\frac{\dot M_2}{M_2}&=\frac{1}{4\sqrt{B_\beta}}
    \left|\frac{\p_t p}{p}\frac{l}{\sqrt p}\frac{\p_t p_h}{p_h}\right|, &-\frac{\dot M_3}{M_3}
 &=\left|\p_t\left(-\frac{1}{2\sqrt{B_\beta}}\frac{l}{\sqrt p}\frac{\p_t p_h}{p_h}\right)\right|,\nonumber\\
-\frac{\dot M_4}{M_4}
 &=\left|\frac{2}{\beta}\frac{l^2}{p}\frac{k^2}{p_h}\frac{\p_t p_h}{p_h}\right|, &-\frac{\dot M_5}{M_5}&=2\sqrt{\frac{\beta-1}{\beta}}\frac{k^2}{p_h}\left|\frac{l}{\sqrt p}\right|, &-\frac{\dot M_6}{M_6}
 &=\frac{\alpha}{2\sqrt{B_\beta}}\left|\frac{l}{p}\frac{\p_t p_h}{\sqrt{p_h}}\right|.\label{4..5}
\end{align}
}

\begin{lem}\label{lem4.2}
For $1\le j\le6$, the above multipliers satisfy
\begin{equation}\label{eq:ghost-bounds}
C_{\alpha,\beta}^{-1}\le M_j\le1,\qquad
\int_0^\infty-\frac{\dot M_j}{M_j}\,\mathrm dt
\le C_{\alpha,\beta},
\end{equation}
where $C_{\alpha,\beta}$ is independent of $\nu,t,k,\eta,l$.
Moreover, for $k\ne0$ one has
\begin{equation}\label{eq:enhanced-rate}
\nu p-\frac{\dot M_1}{M_1}\ge\frac{1}{2}\nu^{1/3}.
\end{equation}
\end{lem}
\begin{proof}
For $k\ne0$, recall $p=k^2+(\eta-kt)^2+l^2=p_h+l^2$. By \eqref{4..5},
\[
 -\frac{\dot{M_2}}{M_2}\le\frac{k^2(\eta-kt)^2}{\sqrt{B_\beta}\,p_h^2},
 \quad
 -\frac{\dot{M_4}}{M_4}\le\frac{4}{|\beta|}\frac{|k|^3|\eta-kt|}{p_h^2},
 \quad
 -\frac{\dot{M_5}}{M_5}\le\frac{2k^2}{c_\beta p_h},\quad -\frac{\dot{M_6}}{M_6}\le\frac{\alpha}{\sqrt{B_\beta}}\frac{|k||l|}{p}.
\]
For $M_3$, we define $f(w)=\frac{|l/k|w}{(1+w^2)\sqrt{1+w^2+|l/k|^2}}$, then $-\dot{M_3}/M_3=|f'(t-\eta/k)|/\sqrt{B_\beta}$ and
\[
 f'(w)=
 \frac{|l/k|\big[(1+|l/k|^2)(1-w^2)-2w^4\big]}
 {(1+w^2)^2(1+w^2+|l/k|^2)^{3/2}}.
\]
For $l\neq 0$, the function $f(w)$ has exactly one critical
point on $(0,\infty)$, which is odd and tends to zero at infinity. Moreover, we have $0\le f(w)\le \frac{w}{1+w^2}\le\frac12$ for $w\ge0$. Thus $\int_{\mathbb R}|f'(w)|\,\mathrm dw
 =4\max_{w\ge0}f(w)\le2$. The same bound can be  immediately obtained for $l=0$.

Integrating the preceding bounds with respect to $t$ yields
\[
\begin{aligned}
 \int_0^\infty -\frac{\dot{M_1}}{M_1}\,\mathrm dt&\le\pi,
 &\int_0^\infty -\frac{\dot{M_2}}{M_2}\,\mathrm dt&\le\frac{\pi}{2\sqrt{B_\beta}},
 &\int_0^\infty -\frac{\dot{M_3}}{M_3}\,\mathrm dt&\le\frac{2}{\sqrt{B_\beta}},\\
 \int_0^\infty -\frac{\dot{M_4}}{M_4}\,\mathrm dt&\le\frac{4}{|\beta|},
 &\int_0^\infty -\frac{\dot{M_5}}{M_5}\,\mathrm dt&\le\frac{2\pi}{c_\beta},
 &\int_0^\infty -\frac{\dot{M_6}}{M_6}\,\mathrm dt&\le\frac{\alpha\pi}{\sqrt{B_\beta}}.
\end{aligned}
\]
Since $-\frac{\dot{M_j}}{M_j}\ge0$ and $M_j(t)=\exp\bigg(-\int_0^t \big(-\frac{\dot{M_j}}{M_j}\big)(s)\,\mathrm ds\bigg)$,
these bounds give \eqref{eq:ghost-bounds}.

Finally, if $|\nu^{1/3}(t-\eta/k)|\le1$, then
$-\dot{M_1}/M_1\ge\nu^{1/3}/2$. If $|\nu^{1/3}(t-\eta/k)|>1$, then
$\nu p\ge\nu (t-\eta/k)^2\ge\nu^{1/3}$.
This proves \eqref{eq:enhanced-rate}.
\end{proof}

The following corollary follows directly from lemma~\ref{lem4.2}.
\begin{cor}[\cite{MR3612004}]\label{cor4.1}
 	For any $g$ and $s \geqslant 0$, the following inequality holds
\begin{equation}\label{4.17}
		\left\|  g_{\neq}  \right\|_{L_t^2 H^s} \leq 2 \nu^{-\frac{1}{6}} \left( \left\|  \sqrt{-\frac{\dot{M_{1}}}{M_{1}}} g_{\neq} \right\|_{L_{t}^2 H^s}+ \nu^{\frac{1}{2}} \left\|  \nabla_{L} g_{\neq} \right\|_{L_t^2 H^s}\right),
	\end{equation}
	which will be used frequently below.
\end{cor}

\noindent\textit{The multiplier $M_7$.}
Inspired by Knobel \cite{Kno2505} and Li--Sun--Wang--Wei--Zhang \cite{LSWZ2025}, we use the multiplier $m_3$ introduced by \cite[(2.21)]{LSWZ2025}, with the reciprocal convention and normalization at $t=0$. The precise relation is
\[
 M_7(t,k,\eta)=\frac{m_3(0,k,\eta)}{m_3(t,k,\eta)},
 \qquad -\frac{\dot M_7}{M_7}=\frac{ \dot{m_3}}{m_3}.
\]
Equivalently, $M_7$ is defined by
\begin{equation}\label{eq:M8}
 -\frac{\dot M_7}{M_7}(t,k,\eta)
 =\sum_{n\in\Z\setminus\{0\}}
 \frac{\langle k-n\rangle^{-1-a_0}}{\langle t-\eta/n\rangle[\log(e+\langle t-\eta/n\rangle)]^{1+a_0}},
 \qquad M_7(0,k,\eta)=1,
\end{equation}
for all $t\ge0$, $(k,\eta,l)\in\mathbb Z\times\R\times\mathbb Z$, where $0<a_0<1/2$ is fixed as in theorem~\ref{1..1}. 

In particular, $M_7$ satisfies the following lemma.
\begin{lem}\label{lem:transfers}
Let $M_7$ be the Fourier multiplier defined in \eqref{eq:M8}. Then $M_7$ is bounded from below and above uniformly with respect to $0<\nu<1$ and $(t,k,\eta)\in\R^+\times\mathbb Z\times\R$ as
\begin{equation}\label{eq:M7-bounds}
 e^{-Ca_0^{-2}}\le M_7\le1,\qquad
 \int_0^\infty-\frac{\dot M_7}{M_7}(t,k,\eta)\,\mathrm dt\le Ca_0^{-2}.
\end{equation}

For $\varsigma=(k,\eta)$ and $\varsigma'=(k',\eta')\in\mathbb Z\times\R$, we have
\begin{equation}\label{eq:single-transfer}
\begin{aligned}
 \frac{1}{\langle t-\eta/k\rangle[\log(e+\langle t-\eta/k\rangle)]^{1+a_0}}
 &\le-\frac{\dot M_7}{M_7}(t,\varsigma),\qquad k\ne0,\\
 p_h(t,\varsigma)^{-1/2}
 &\lesssim_{a_0}\langle\varsigma-\varsigma'\rangle^{1+a_0}
 \sqrt{-\dot M_7/M_7(t,\varsigma')},\qquad k\ne0,\\
 \frac{1}{\langle t-\eta'/k'\rangle[\log(e+\langle t-\eta'/k'\rangle)]^{1+a_0}}
 &\lesssim_{a_0}\langle\varsigma-\varsigma'\rangle^{1+a_0}
 \sqrt{\left(-\frac{\dot M_7}{M_7}\right)(t,\varsigma')
       \left(-\frac{\dot M_7}{M_7}\right)(t,\varsigma)},\quad k'\ne0.
\end{aligned}
\end{equation}

\end{lem}
\begin{proof}
The bounds \eqref{eq:M7-bounds} and $\eqref{eq:single-transfer}_1$
follow from \cite[lemma~2.3]{LSWZ2025} and the normalization above.
For the last two inequalities, we retain the slightly sharper frequency
loss $1+a_0$, instead of $1+a_0+a_0^2/2$ obtained directly from the
comparison estimate in that lemma. We give only this refinement.
Set $w(s)=\langle s\rangle[\log(e+\langle s\rangle)]^{1+a_0}$.
For $k\ne0$, Peetre's inequality and logarithmic translation give
\begin{equation}\label{eq:psi-translation}
\begin{aligned}
 \frac{w(t-\eta'/k)}{w(t-\eta/k)}
 &\lesssim\langle\eta-\eta'\rangle
 [\log(e+\langle\eta-\eta'\rangle)]^{1+a_0}
 \lesssim_{a_0}\langle\eta-\eta'\rangle^{1+a_0}.
\end{aligned}
\end{equation}
Choosing the summand $n=k$ in \eqref{eq:M8} therefore yields
\begin{align*}
 -\frac{\dot{M_7}}{M_7}(t,\varsigma')
 &\ge\frac{\langle k'-k\rangle^{-1-a_0}}{w(t-\eta'/k)}
 \gtrsim_{a_0}\frac{\langle\varsigma-\varsigma'\rangle^{-2(1+a_0)}}{w(t-\eta/k)}.
\end{align*}
Since $p_h(t,\varsigma)^{-1}\le\langle t-\eta/k\rangle^{-2}
\lesssim_{a_0}w(t-\eta/k)^{-1}$, taking square roots proves
$\eqref{eq:single-transfer}_2$.
For $k'\ne0$, interchange $\varsigma$ and $\varsigma'$ in the preceding
comparison and combine it with
$w(t-\eta'/k')^{-1}\le -\big(\dot{M_7}/M_7\big)(t,\varsigma')$.
Taking the geometric mean gives $\eqref{eq:single-transfer}_3$.
\end{proof}

The following lemma  provides the weighted nonzero-mode
estimates used in the nonlinear analysis.
Its decay rate and  regularity requirements differ from those
of theorem~\ref{thm:linear}, whose nonzero-mode estimates
are proved in appendix~\ref{secA}. 
\begin{lem}\label{lem3.4}
	Assume that $0<\nu\le1$, $\alpha>0$, and  $\beta\in\mathbb R$ satisfy $B_\beta>1/4$. Fix $0<a_0<1/2$.
Then the following estimates hold:
	\begin{align}
		&\left\| U^1_{\neq} (t) \right\|_{L^{2}} \lesssim_{\alpha, \beta, a_0} \, e^{-\lambda  \nu^{\frac{1}{3}} t} \left(\left\|  u_{\mathrm{in}}\right\|_{H^{1}}+\left\|  \theta_{\mathrm{in}}\right\|_{H^{1}} \right) , \label{3.1x}  \\
		&{\left\langle {t} \right\rangle} \left\|  U^2_{\neq} (t) \right\|_{L^{2}} \lesssim_{\alpha, \beta, a_0} \, e^{-\lambda \nu^{\frac{1}{3}} t} \left(\left\|  u_{\mathrm{in}}\right\|_{H^{2}}+\left\|  \theta_{\mathrm{in}}\right\|_{H^{2}} \right) , \label{3.2x}  \\
		&\left\| (U^3_{\neq}, \Theta_{\neq}) (t) \right\|_{L^2} \lesssim_{\alpha, \beta, a_0} \, e^{-\lambda \nu^{\frac{1}{3}} t} \left(\left\|  u^3_{\mathrm{in}}\right\|_{H^{1}} + c_{\beta} \left\|  \p_xu^2_{\mathrm{in}}-\p_yu^1_{\mathrm{in}}\right\|_{L^2}+\left\|  \theta_{\mathrm{in}}\right\|_{H^{1}} \right),\label{3.3x}  
	\end{align}
where $\lambda=\frac{1}{16}\frac{2\sqrt{B_\beta}-1}{2\sqrt{B_\beta}+1}$.
\end{lem}
\begin{proof}
To control the evolution of the nonzero modes, we introduce the Fourier multiplier
\begin{equation}\label{888}
	M=M_1\left(\prod_{j=2}^6M_j\right)^{\kappa_{\beta}}M_7,\qquad \mathcal{A}:=mMe^{\lambda \nu^{\frac{1}{3}}t},
\end{equation} 
for $\kappa_{\beta}=(1-1/(2\sqrt{B_{\beta}}))^{-1}$, which incorporate the multiplier $m$ already used in the linear analysis and additional ghost weights ${M}$ as in \eqref{888} as well as the time weight $e^{\lambda\nu^{\frac{1}{3}}t}$, which $\lambda$ can be confirmed later. By lemmas~\ref{lem4.2} and \ref{lem:transfers}, there exists
$C_M=C_M(\alpha,\beta,a_0)\ge1$, independent of
$\nu,t,k,\eta,l$, such that
\[
 C_M^{-1}\le M\le1.
\]
Recall the system  \eqref{Non0L}, energy estimates yield that 
\begin{equation}\label{eq:linear-basic-energy}
\begin{aligned}
    \frac{1}{2}&\frac{d}{dt}|\mathcal{A}(\widehat{Q_{\neq}},\widehat{K_{\neq}},\widehat{H_{\neq}})|^2  + \nu p |\mathcal{A}(\widehat{Q_{\neq}},\widehat{K_{\neq}},\widehat{H_{\neq}})|^2 +\left[(-\frac{\dot{M}}{M})+(-\frac{\dot{m}}{m})-\lambda\nu^{\frac{1}{3}}\right]|\mathcal{A}(\widehat{Q_{\neq}},\widehat{K_{\neq}},\widehat{H_{\neq}})|^2\\
    &-\left(\frac{1}{4}\frac{\p_t p}{p}+\frac{1}{2}\frac{l^2}{p}\frac{\p_t p_h}{p_h}\right)\left(|\widehat{\mathcal{A}Q_{\neq}}|^2-|\widehat{\mathcal{A}K_{\neq}}|^2\right)-\frac{1}{4}\frac{\p_t p}{p}|\widehat{\mathcal{A}H_{\neq}}|^2+\frac{2}{c_{\beta}}\frac{k^2}{p_h}\frac{l}{\sqrt{p}}\Rea\big(\widehat{\mathcal{A}K_{\neq}}\overline{\widehat{\mathcal{A}Q_{\neq}}}\big)=0.
\end{aligned}
\end{equation}
Consider the cross terms, we introduce the Fourier multiplier $G$ as 
\begin{equation}\label{eq:G}
G=\frac{\sgn(\beta)}{2\sqrt\B}\frac{l}{\sqrt p}\frac{\dot p_h}{p_h}.
\end{equation}
thus one has 
\begin{align}\label{sss}
	&\p_t\Rea\left(G\widehat{\mathcal{A}Q_{\neq}}\overline{\widehat{\mathcal{A}K_{\neq}}}\right) \nonumber\\
	&\quad=\Rea\Big\{(\p_t G)\widehat{\mathcal{A}Q_{\neq}}\overline{\widehat{\mathcal{A}K_{\neq}}} + \left[(\frac{\dot{M}}{M})+(\frac{\dot{m}}{m})+\lambda\nu^{\frac{1}{3}}\right]2G\widehat{\mathcal{A}Q_{\neq}}\overline{\widehat{\mathcal{A}K_{\neq}}}-G\left[\nu p\widehat{\mathcal{A}Q_{\neq}}\right. \nonumber\\
	&\qquad \left.-\left(\frac{1}{4}\frac{\p_t p}{p}+\frac{1}{2}\frac{l^2}{p}\frac{\p_t p_h}{p_h}\right)\widehat{\mathcal{A}Q_{\neq}}-\sgn(\beta)\sqrt{B_{\beta}}\frac{l}{\sqrt{p}}\widehat{\mathcal{A}K_{\neq}}+\alpha\sqrt{\frac{p_h}{p}}\widehat{\mathcal{A}H_{\neq}}+\frac{2}{c_{\beta}}\frac{k^2}{p_h}\frac{l}{\sqrt{p}}\widehat{\mathcal{A}K_{\neq}}\right]\overline{\widehat{\mathcal{A}K_{\neq}}} \nonumber\\
	&\qquad-G\widehat{\mathcal{A}Q_{\neq}}\left[\nu p\overline{\widehat{\mathcal{A}K_{\neq}}}+\left(\frac{1}{4}\frac{\p_t p}{p}+\frac{1}{2}\frac{l^2}{p}\frac{\p_t p_h}{p_h}\right)\overline{\widehat{\mathcal{A}K_{\neq}}}+\sgn(\beta)\sqrt{B_{\beta}}\frac{l}{\sqrt{p}}\overline{\widehat{\mathcal{A}Q_{\neq}}}\right]\Big\}.
\end{align}
This combined with \eqref{eq:linear-basic-energy} implies that
\begin{align}\label{999}
    &\frac{1}{2}\frac{d}{dt}\left(|\mathcal{A}(\widehat{Q_{\neq}},\widehat{K_{\neq}},\widehat{H_{\neq}})|^2+2G\Rea\big(\widehat{\mathcal{A}Q_{\neq}}\overline{\widehat{\mathcal{A}K_{\neq}}}\big)\right)\nonumber\\ &+\left[ \nu p+\left(-\frac{\dot{M}}{M}\right)\right]\left(|\mathcal{A}(\widehat{Q_{\neq}},\widehat{K_{\neq}},\widehat{H_{\neq}})|^2+2G\Rea\big(\widehat{\mathcal{A}Q_{\neq}}\overline{\widehat{\mathcal{A}K_{\neq}}}\big)\right)\nonumber\\
    =&\lambda\nu^{\frac{1}{3}}\left(|\mathcal{A}(\widehat{Q_{\neq}},\widehat{K_{\neq}},\widehat{H_{\neq}})|^2+2G\Rea\big(\widehat{\mathcal{A}Q_{\neq}}\overline{\widehat{\mathcal{A}K_{\neq}}}\big)\right)\nonumber\\
    &-\left[(-\frac{\dot{m}}{m})-\frac{1}{4}\frac{\p_t p}{p}\right]|\widehat{\mathcal{A}Q_{\neq}}|^2-\left[(-\frac{\dot{m}}{m})+\frac{1}{4}\frac{\p_t p}{p}\right]|\widehat{\mathcal{A}K_{\neq}}|^2-\left[(-\frac{\dot{m}}{m})-\frac{1}{4}\frac{\p_t p}{p}\right]|\widehat{\mathcal{A}H_{\neq}}|^2\nonumber\\
    &-(-\frac{\dot{m}}{m})2G\Rea\big(\widehat{\mathcal{A}Q_{\neq}}\overline{\widehat{\mathcal{A}K_{\neq}}}\big)+(\p_t G)\Rea\big(\widehat{\mathcal{A}Q_{\neq}}\overline{\widehat{\mathcal{A}K_{\neq}}}\big)-\frac{1}{\beta}\frac{l^2}{p}\frac{k^2}{p_h}\frac{\p_t p_h}{p_h}|\widehat{\mathcal{A}K_{\neq}}|^2\nonumber\\
    &-\frac{2}{c_{\beta}}\frac{k^2}{p_h}\frac{l}{\sqrt{p}}\Rea\big(\widehat{\mathcal{A}K_{\neq}}\overline{\widehat{\mathcal{A}Q_{\neq}}}\big)-\frac{\sgn(\beta)\alpha}{2\sqrt{B_{\beta}}}\frac{l}{p}\frac{\p_tp_h}{\sqrt{p_h}}\Rea\big(\widehat{\mathcal{A}H_{\neq}}\overline{\widehat{\mathcal{A}K_{\neq}}}\big),
\end{align}
where the choice of $G$ as in \eqref{eq:G} cancels the stretching terms containing
$(l^2/p)(\dot p_h/p_h)$ in \eqref{eq:linear-basic-energy}.

Denote $e_{\neq}(t):=|\mathcal{A}(\widehat{Q_{\neq}},\widehat{K_{\neq}},\widehat{H_{\neq}})|^2+2G\Rea\big(\widehat{\mathcal{A}Q_{\neq}}\overline{\widehat{\mathcal{A}K_{\neq}}}\big)$, the expression of $G$ yields that 
\begin{equation}\label{00}
\left(1-\frac{1}{2\sqrt{B_\beta}}\right)
 |\mathcal A(\widehat Q_{\neq},\widehat K_{\neq},\widehat H_{\neq})|^2
\le e_{\neq}(t)\le\left(1+\frac{1}{2\sqrt{B_\beta}}\right)
 |\mathcal A(\widehat Q_{\neq},\widehat K_{\neq},\widehat H_{\neq})|^2.
\end{equation}
Notice that the definition of $m$ as in \eqref{eq:m-formula}, we have $-\frac{\dot{m}}{m}+\frac{\nu p}{2}-\frac{1}{4}\big|\frac{\p_t p}{p}\big|\gtrsim 0$. A direct calculation yields
\begin{align}
 &\frac{1}{2}\frac{d}{dt}e_{\neq}(t)
 +\left[\frac{3}{4}\nu p+\left(-\frac{\dot M}{M}\right)\right]e_{\neq}(t)\nonumber\\
 &\quad\le\lambda\nu^{1/3}e_{\neq}(t)
 +\left|2\frac{\dot m}{m}G\right|
   \left|\Rea\big(\widehat{\mathcal AQ_{\neq}}\overline{\widehat{\mathcal AK_{\neq}}}\big)\right|
 +|\p_tG|\left|\Rea\big(\widehat{\mathcal AQ_{\neq}}\overline{\widehat{\mathcal AK_{\neq}}}\big)\right|\nonumber\\
 &\qquad+\left|\frac{1}{\beta}\frac{l^2}{p}\frac{k^2}{p_h}\frac{\p_tp_h}{p_h}\right|
       |\widehat{\mathcal AK_{\neq}}|^2
 +\frac{2}{c_\beta}\frac{k^2}{p_h}\frac{|l|}{\sqrt p}
   \left|\Rea\big(\widehat{\mathcal AK_{\neq}}\overline{\widehat{\mathcal AQ_{\neq}}}\big)\right|+\frac{\alpha}{2\sqrt{B_\beta}}
       \left|\frac{l}{p}\frac{\p_tp_h}{\sqrt{p_h}}\right|
   \left|\Rea\big(\widehat{\mathcal AH_{\neq}}\overline{\widehat{\mathcal AK_{\neq}}}\big)\right|\nonumber\\
 &\quad\le\lambda\nu^{1/3}e_{\neq}(t)
 +\frac{1}{2}\sum_{j=2}^6\left(-\frac{\dot M_j}{M_j}\right)
   |\mathcal A(\widehat Q_{\neq},\widehat K_{\neq},\widehat H_{\neq})|^2,\label{1010}
\end{align}
where we have used $\left|2G\dot m/m\right|\le-\dot M_2/M_2$ by \eqref{eq:m-formula} and \eqref{4..5} on $|t-\eta/k|\le\tau_\nu=1000\nu^{-1/3}$. Outside this interval, $\dot m=0$, so the term containing $G\dot m/m$ vanishes. 

In the following, we introduce the energy and dissipation
\begin{equation}\label{eq:EF}
\begin{aligned}
    E_{\neq,s}(t):=&\|\mathcal{A}(Q_{\neq},K_{\neq},H_{\neq})\|_{H^s}^2 +2\Rea\langle G\mathcal{A}Q_{\neq},\mathcal{A}K_{\neq}\rangle_{H^s},\\
    F_{\neq,s}(t):=&\nu\|\mathcal{A}\na_L (Q_{\neq},K_{\neq},H_{\neq})\|_{H^s}^2 +2\nu\Rea\langle G\mathcal{A}\na_L Q_{\neq},\mathcal{A}\na_L K_{\neq}\rangle_{H^s}\\
    &+\left\|\mathcal{A}\sqrt{-\frac{\dot{M}}{M}} \left(Q_{\neq},K_{\neq},H_{\neq}\right)\right\|_{H^s}^2 +2\Rea\left\langle G\mathcal{A}\sqrt{-\frac{\dot{M}}{M}} Q_{\neq},\mathcal{A}\sqrt{-\frac{\dot{M}}{M}} K_{\neq}\right\rangle_{H^s},
\end{aligned}
\end{equation}
for any $s \geq 0$. Using \eqref{00}, $  E_{\neq,s}(t) $ and $  F_{\neq,s}(t) $ satisfy
\begin{align}
    & \left(1-\frac{1}{2\sqrt{B_{\beta}}}\right)\|\mathcal{A}(Q_{\neq},K_{\neq},H_{\neq})\|_{H^s}^2\leq E_{\neq,s}(t)\leq  \left(1+\frac{1}{2\sqrt{B_{\beta}}}\right)\|\mathcal{A}(Q_{\neq},K_{\neq},H_{\neq})\|_{H^s}^2,\label{eq:Equiv1}\\
    & \left(1-\frac{1}{2\sqrt{B_{\beta}}}\right)\left(\nu\|\mathcal{A}\na_L(Q_{\neq},K_{\neq},H_{\neq})\|_{H^s}^2 + \left\|\mathcal{A}\sqrt{-\frac{\dot{M}}{M}}(Q_{\neq},K_{\neq},H_{\neq})\right\|_{H^s}^2\right)\leq F_{\neq,s}(t)\nonumber\\
	&\quad \leq \left(1+\frac{1}{2\sqrt{B_{\beta}}}\right)\left(\nu\|\mathcal{A}\na_L(Q_{\neq},K_{\neq},H_{\neq})\|_{H^s}^2 + \left\|\mathcal{A}\sqrt{-\frac{\dot{M}}{M}}(Q_{\neq},K_{\neq},H_{\neq})\right\|_{H^s}^2\right)\label{eq:Equiv2},
\end{align}
thus, recalling \eqref{1010} and going back to the physical space, we obtain that
{\small
\begin{align*}
    &\frac{\rm d}{{\rm d}t}E_{\neq,s}(t) + \left(1-\frac{1}{2\sqrt{B_{\beta}}}\right)\nu\|\mathcal{A}\na_L(Q_{\neq},K_{\neq},H_{\neq})\|_{H^s}^2 + \left(1-\frac{1}{2\sqrt{B_{\beta}}}\right)\|\mathcal{A}\sqrt{-\frac{\dot{M}}{M}}(Q_{\neq},K_{\neq},H_{\neq})\|_{H^s}^2\\
    &\quad\leq2\lambda\nu^{\frac{1}{3}}\left(1+\frac{1}{2\sqrt{B_{\beta}}}\right)\|\mathcal{A}(Q_{\neq},K_{\neq},H_{\neq})\|_{H^s}^2.
\end{align*}
}
Taking $\lambda=\frac{1}{16}\frac{2\sqrt{B_{\beta}}-1}{2\sqrt{B_{\beta}}+1}$  and integrating the above inequality in time yields that
\begin{align*}
    &\|\mathcal{A}(Q_{\neq},K_{\neq},H_{\neq})\|_{H^s}^2 + \nu\|\mathcal{A}\na_L(Q_{\neq},K_{\neq},H_{\neq})\|_{L^2_tH^s}^2 + \left\|\mathcal{A}\sqrt{-\dot{M}/M}(Q_{\neq},K_{\neq},H_{\neq})\right\|_{L^2_tH^s}^2\\
    & \quad\leq \frac{2\sqrt{B_{\beta}}+1}{2\sqrt{B_{\beta}}-1}\|\mathcal{A}(Q_{\neq},K_{\neq},H_{\neq})(t=0)\|_{H^s}^2 +\frac{\nu^{\frac{1}3}}{8}\|\mathcal{A}(Q_{\neq},K_{\neq},H_{\neq})\|_{L^2_tH^s}^2.
\end{align*}

Due to $\mathcal{A}=e^{\lambda \nu^{\frac{1}{3}} t }mM$, \eqref{eq:m-initial} and \eqref{4.17},  we have
\begin{align}
	\|mM(Q_{\neq},K_{\neq},H_{\neq})(t)\|_{H^s}^2&\lesssim \frac{2\sqrt{B_\beta}+1}{2\sqrt{B_\beta}-1}e^{-2\lambda\nu^{1/3}t}\|(Q_{\mathrm{in}},K_{\mathrm{in}},H_{\mathrm{in}})\|_{H^{s+1/2}}^2,\label{3x.0}\\
	\||\na_{x,z}|^{-1}mM(Q_{\neq},K_{\neq},H_{\neq})(t)\|_{H^s}&\lesssim e^{-\lambda\nu^{1/3}t}(\|u_{\rm in}^3\|_{H^{s+1}}+c_{\beta}\|\p_xu_{\rm in}^2-\p_yu_{\rm in}^1\|_{H^s}+\|\theta_{\rm in}\|_{H^{s+1}}).\nonumber
\end{align}
For the second inequality, apply the same pointwise energy estimate
with the time-independent weight $|k,l|^{-2}$. Indeed, we have $|k,l|^{-1}m(0)|\widehat{Q_{\mathrm{in}}}|\le |k|^{-1}|k,\eta,l||\widehat{u_{\mathrm{in}}^3}|$, $|k,l|^{-1}m(0)|\widehat{K_{\mathrm{in}}}|\le c_\beta |k|^{-1}
 |\widehat{\partial_xu_{\mathrm{in}}^2-\partial_yu_{\mathrm{in}}^1}|$ and $|k,l|^{-1}m(0)|\widehat{H_{\mathrm{in}}}|\le |k|^{-1}|k,\eta,l||\widehat{\theta_{\mathrm{in}}}|$. Here $M(0)=1$ and $|k|\ge1$ on the nonzero modes.

Next, recalling the definition of $Q$, $K$ and  $H$ and combining with \eqref{4.11-2}, we have 
\begin{align}\label{3x.1-linear}
	|\widehat{U^3_{\neq}}|&= \frac{p_h^{1/2}}{p^{3/4}}|\widehat{Q_{\neq}}| \lesssim C_M\frac{\sqrt{p_h}}{\sqrt{p}|k,l|}|mM\widehat{Q_{\neq}}|\lesssim \frac{C_M}{|k,l|}|mM\widehat{Q_{\neq}}|,\\
	|\widehat{\Theta_{\neq}}|&=p^{-\frac{1}{4}}|\widehat{H_{\neq}}|\lesssim \frac{C_M}{|k,l|}|mM\widehat{H_{\neq}}|.\label{3x.2}
\end{align}
Hence, together  \eqref{3x.0} with \eqref{3x.1-linear}-\eqref{3x.2} gives 
\begin{equation}
\begin{aligned}\label{3x.3}
	\left\|{U^3_{\neq}},{\Theta_{\neq}}\right\|_{L^2}^2\lesssim& C_M^2\left\||\na_{x,z}|^{-1} mM({Q_{\neq}},{H_{\neq}})   \right\|_{L^2}^2 \\
	\lesssim& C_M^2\frac{2\sqrt{B_{\beta}}+1}{2\sqrt{B_{\beta}}-1} e^{-2\lambda \nu^\frac{1}{3} t} \left(\| u^3_{\mathrm{in}}\|_{H^1}^2+c_{\beta}^2 \left\|  \p_xu^2_{\mathrm{in}}-\p_yu^1_{\mathrm{in}}\right\|_{L^2}^2+\| \theta_{\mathrm{in}}\|_{H^1}^2 \right).
\end{aligned}
\end{equation}

The incompressibility condition $\partial_{x} U^1+\partial_{y}^L U^2=-\partial_{z} U^3$ and the definition of $W^3$ yield that
\begin{equation}\label{3.20}
	\begin{cases}
		\left(\partial_{x x}+\partial_{y y}^L\right) U^{1}=-\partial_{y}^{L} W^{3}-\partial_{x z} U^{3},  \\
		\left(\partial_{x x}+\partial_{y y}^L\right) U^{2}=\partial_{x} W^{3}-\partial_{y z}^{L} U^{3}.
	\end{cases}
\end{equation}
By \eqref{4.11-2} and \eqref{3.20}, we get $\widehat{U_{\neq}^1}$ and $\widehat{U_{\neq}^2}$ which satisfy
\begin{align}
	|\widehat{U_{\neq}^1}|\lesssim &\frac{1}{p_h}\bigg(  \frac{|\eta-kt| p_h^{\frac{1}{2}}}{c_{\beta} p^{\frac{1}{4}}} |\widehat{K_{\neq}}|+ \frac{|k||l| p_h^{\frac{1}{2}}}{p^{\frac{3}{4}}} |\widehat{Q_{\neq}}|\bigg)\lesssim  \frac{C_M}{|k,l|} \left(c_{\beta}^{-1}|mM\widehat{K_{\neq}}| + |mM\widehat{Q_{\neq}}|\right),\label{ll6}\\
	|\widehat{U_{\neq}^2}|\lesssim &\frac{1}{p_h}\bigg(  \frac{|k| p_h^{\frac{1}{2}}}{c_{\beta} p^{\frac{1}{4}}} |\widehat{K_{\neq}}|+ \frac{|\eta-kt||l| p_h^{\frac{1}{2}}}{p^{\frac{3}{4}}} |\widehat{Q_{\neq}}|\bigg)\lesssim   \frac{C_M}{|k,l|}\bigg( \frac{|k|}{c_{\beta} \sqrt{p_h}}|mM\widehat{K_{\neq}}| +\frac{|l|}{\sqrt{p}} |mM\widehat{Q_{\neq}}|\bigg).\label{ll7-linear}
\end{align}
Similar to the estimate \eqref{3x.3}, this together with \eqref{3x.0}  leads to 
\begin{align}
	\left\|{U^1_{\neq}}(t)\right\|_{L^2}^2\lesssim& C_M^2\left(c_{\beta}^{-2}\left\||\na_{x,z}|^{-1}mMK_{\neq}\right\|_{L^2}^2 + \left\||\na_{x,z}|^{-1}mMQ_{\neq}\right\|_{L^2}^2\right) \nonumber\\
	\lesssim& C_M^2\frac{2\sqrt{B_{\beta}}+1}{2\sqrt{B_{\beta}}-1} \max \left\{ c_{\beta}^2,c_{\beta}^{-2} \right\} e^{-2\lambda \nu^{\frac{1}{3}} t} \left(\| u_{\mathrm{in}}\|_{H^1}^2+\| \theta_{\mathrm{in}}\|_{H^1}^2 \right),\label{4x.5}\\
	\left\|{U^2_{\neq}}(t)\right\|_{L^2}^2\lesssim& {\left\langle {t} \right\rangle}^{-2} C_M^2\left(c_{\beta}^{-2}\left\||\na_{x,z}|^{-1}mMK_{\neq}\right\|_{H^{1}}^2 + \left\||\na_{x,z}|^{-1}mMQ_{\neq}\right\|_{H^{1}}^2\right) \nonumber\\
	\lesssim& {\left\langle {t} \right\rangle}^{-2} \,C_M^2\frac{2\sqrt{B_{\beta}}+1}{2\sqrt{B_{\beta}}-1} \max \left\{ c_{\beta}^2,  c_{\beta}^{-2} \right\} e^{-2\lambda \nu^{\frac{1}{3}} t} \left(\| u_{\mathrm{in}}\|_{H^2}^2+\| \theta_{\mathrm{in}}\|_{H^{2}}^2 \right),\label{4x.6}
\end{align}
here we used the weighted initial-data estimate above and $\langle t\rangle(p_h^{-1/2}|k|+p^{-1/2}|l|)\lesssim\langle k,\eta,l\rangle$ for $k\ne0$. Hence, together \eqref{3x.3}, \eqref{4x.5} with \eqref{4x.6} gives \eqref{3.1x}-\eqref{3.3x} and  the proof of lemma~\ref{lem3.4} is finished. It is worth noting that in lemma~\ref{lem3.4} where we require $B_{\beta} > 1/4$, but in fact, obtaining the enhanced dissipation of $U_{\neq}$ and $ \Theta_{\neq}$ does not necessitate this restriction. We have provided a detailed explanation of this point in appendix~\ref{secA} later on.
\end{proof}

\subsection{Analysis of the zero modes $(u_0,\theta_0)$}\label{sec3.2}
\qquad This section is devoted to the dynamics of the zero modes in the $x$ frequency, which satisfy
\begin{equation}\label{1..7}
	\begin{cases}
		\partial_{t} u_0^1-\nu \Delta_{y, z} u_0^1+
		(1-\beta)u_0^2=0, \\
		\partial_{t} u_0^2-\nu \Delta_{y, z} u_0^2+
		\beta \partial_{z} \Delta_{y, z}^{-1} \partial_{z} u_0^1-\alpha \partial_{y} \Delta_{y, z}^{-1} \partial_{z} \theta_0=0, \\
		\partial_{t} u_0^3-\nu \Delta_{y, z} u_0^3-\beta \partial_{z} \Delta_{y, z}^{-1} \partial_{y} u_0^1+\alpha \partial_{y} \Delta_{y, z}^{-1} \partial_{y} \theta_{0}=0, \\
		\partial_{t} \theta_0-\nu \Delta_{y, z} \theta_0-
		\alpha u_0^3=0, \\
		u_0(t=0)=(u_{\mathrm{in}})_0, \quad \theta_0(t=0)=(\theta_{\mathrm{in}})_0.
	\end{cases}
\end{equation}

For the linearized system \eqref{1..7}, we express the solution in an explicit form via Fourier transform and Duhamel’s principle. One gives the following lemma.
\begin{lem}\label{lem3.5}
	$(1)$ Assume that $l=0$, then \eqref{1..7} degenerates into the following heat equations 
	\begin{equation}\label{3.2.2}
		\begin{cases}
			\partial_{t} \overline{u}_0^1-\nu \partial_{y}^2 \overline{u}_0^1=0, \\
			\overline{u}_{0}^{2}=0,\\
			\partial_{t} \overline{u}_0^3-\nu \partial_{y}^2 \overline{u}_0^3+\alpha \overline{\theta}_0=0, \\
			\partial_{t} \overline{\theta}_0-\nu \partial_{y}^2 \overline{\theta}_0-
			\alpha \overline{u}_0^3=0,\\
			\overline{u_0}(t=0)=\overline{({u}_{\mathrm{in}})_0}, \overline{\theta_0}(t=0)=\overline{({\theta}_{\mathrm{in}})_0}.
		\end{cases}
	\end{equation}
	In this case, the double zero modes $(\overline{u}_0, \overline{\theta}_0)$  of \eqref{3.2.2} is easily obtained, which is \eqref{thm1.1.1.1}--\eqref{thm1.1.3} as in theorem~\ref{thm:linear}.

	$(2)$ Assume that $l \neq 0$, the simple zero modes $\widetilde{u}_{0}$ of system \eqref{1..7} can be expressed as the formulas \eqref{3.2.4}-\eqref{3.2.7}. 
\end{lem}
\begin{proof}
First, we give the proof for the double zero modes. When $l=0$, $\overline{u}_0^1$ satisfies the standard heat equation which gives the estimate \eqref{thm1.1.1.1}. The coupling of $\overline{u}_0^3$ and $\overline{\theta}_0$  does not lead to energy dissipation or enhancement. The Fourier transform is as follows
\begin{equation*} 
	\begin{cases}
		\partial_{t} \widehat{\overline{u}_0^3}+\nu \eta^2 \widehat{\overline{u}_0^3}+\alpha \widehat{\overline{\theta}_0}=0, \\
		\partial_{t} \widehat{\overline{\theta}_0}+\nu \eta^2 \widehat{\overline{\theta}_0}-
		\alpha \widehat{\overline{u}_0^3}=0,\\
		\widehat{\overline{u^3_0}}(t=0)=\widehat{\overline{({u}^3_{\mathrm{in}})_0}}, \widehat{\overline{\theta_0}}(t=0)=\widehat{\overline{({\theta}_{\mathrm{in}})_0}}.
	\end{cases}
\end{equation*}
Let	$\phi=\widehat{\overline{u}_0^3}+i \widehat{\overline{\theta}_0}$ and $\psi=\widehat{\overline{u}_0^3}-i \widehat{\overline{\theta}_0}$ with the initial data $\phi_{\mathrm{in}}=\widehat{\overline{({u}^3_{\mathrm{in}})_0}}+i \widehat{\overline{({\theta}_{\mathrm{in}})_0}}$ and $\psi_{\mathrm{in}}=\widehat{\overline{({u}^3_{\mathrm{in}})_0}}-i \widehat{\overline{({\theta}_{\mathrm{in}})_0}}$, then we have
\begin{equation*}
	\partial_{t} \phi=(-\nu \eta^2+i \alpha) \phi, \quad \partial_{t} \psi=(-\nu \eta^2-i \alpha) \psi.
\end{equation*}
Thus, a direct calculation yields
\begin{equation*}
	\phi(t, \eta)=e^{-\nu \eta^2 t+ i \alpha t}\phi_{\mathrm{in}}, \quad \psi(t, \eta)=e^{-\nu \eta^2 t- i \alpha t}\psi_{\mathrm{in}}.
\end{equation*}
Noticing that $\widehat{\overline{u}_0^3}=\frac{\phi+\psi}{2}$ and $\widehat{\overline{\theta}_0}=\frac{\phi-\psi}{2i}$, we immediately have  \eqref{thm1.1.2}  and \eqref{thm1.1.3}.

Next, assume $l\ne0$.  Denote by $\mathsf L_0(\eta,l)$ the coefficient matrix of the Fourier-transformed
zero-mode system. By using the incompressibility condition, $\widetilde{u}_{0}^{3}$   can be expressed as $-\partial_{z}^{-1} \partial_{y} \widetilde{u}_{0}^{2}(t, y, z)$. To simplify the calculation, taking  the Fourier transform to the $\widetilde{u}_{0}^{1}$,  $\widetilde{u}_{0}^{2}$ and $\widetilde{\theta}_{0}$ equations in \eqref{1..7}, it yields
\begin{equation*}
	\partial_{t} \begin{bmatrix}
		\widehat{\widetilde{u}_{0}^{1}} \\
		\widehat{\widetilde{u}_{0}^{2}} \\
		\widehat{\widetilde{\theta}_{0}}
	\end{bmatrix}=\begin{bmatrix}
		-\nu(\eta^2+l^2)  & \beta-1 & 0 \\
		-\beta \frac{l^2}{\eta^2+l^2} & -\nu(\eta^2+l^2) & \alpha \frac{\eta l}{\eta^2+l^2} \\
		0 & -\alpha \frac{\eta}{l} & -\nu(\eta^2+l^2)
	\end{bmatrix}\begin{bmatrix}
		\widehat{\widetilde{u}_{0}^{1}} \\
		\widehat{\widetilde{u}_{0}^{2}} \\
		\widehat{\widetilde{\theta}_{0}}
	\end{bmatrix}=:\mathsf{L}_0\begin{bmatrix}
		\widehat{\widetilde{u}_{0}^{1}} \\
		\widehat{\widetilde{u}_{0}^{2}} \\
		\widehat{\widetilde{\theta}_{0}}
	\end{bmatrix}.
\end{equation*}
Then, the simple zero modes satisfy the following expression
\begin{equation}\label{7.xx}
\left(
\widehat{\widetilde u_0^1}(t,\eta,l),
\widehat{\widetilde u_0^2}(t,\eta,l),
\widehat{\widetilde\theta_0}(t,\eta,l)\right)^T
=e^{\mathsf L_0 t}
\left(
\widehat{\widetilde{(u^1_{\mathrm{in}})_0}},
\widehat{\widetilde{(u^2_{\mathrm{in}})_0}},
\widehat{\widetilde{(\theta_{\mathrm{in}})_0}}
\right)^T.
\end{equation}

We compute the fundamental matrix $e^{\mathsf{L}_0 t}$ explicitly. The characteristic polynomial associated
with $\mathsf{L}_0$ is given by
\begin{equation*}
	\left(\lambda+ \nu \left( \eta^2+l^2 \right)\right)\left[ \left(\lambda+ \nu \left( \eta^2+l^2 \right)\right)^2+\omega_0^2\right]=0,
\end{equation*}
where $\omega_0=\omega_0(\alpha, \beta, \eta, l)=\sqrt{\frac{\alpha^2 \eta^2+ B_{\beta} l^2}{\eta^2+l^2}}$ and $B_{\beta}>0$. The eigenvalues of $\mathsf L_0$ are
\begin{equation*}
	\lambda_1=-\nu |\eta, l|^2, \quad \lambda_2=-\nu |\eta, l|^2+i\, \omega_0, \quad  \quad \lambda_3=-\nu |\eta, l|^2-i\, \omega_0,
\end{equation*}
Note that the eigenvalues $\lambda_2$ and $\lambda_3$ are a pair of conjugate complex roots. The eigenvectors $\mathcal{V}_{i}$ $(i=1, 2, 3)$ associated with  $\lambda_{i}$  satisfy
\begin{align*}
	\left( \lambda_{i} \textbf{I} - \mathsf{L}_0 \right) \mathcal{V}_{i}=\begin{bmatrix}
		\lambda_{i}+\nu(\eta^2+l^2)  & 1-\beta & 0 \\
		\beta \frac{l^2}{\eta^2+l^2} & \lambda_{i}+\nu(\eta^2+l^2) & -\alpha \frac{\eta l}{\eta^2+l^2} \\
		0 & \alpha \frac{\eta}{l} & \lambda_{i}+\nu(\eta^2+l^2)
	\end{bmatrix} \begin{bmatrix}
		\mathcal{V}_{i 1}\\
		\mathcal{V}_{i 2}\\
		\mathcal{V}_{i 3}
	\end{bmatrix}=0,
\end{align*}
Let 
\begin{equation}\label{3.2.9}
	\mathsf{S}_0:=[\mathcal{V}_{1}, \mathcal{V}_{2}, \mathcal{V}_{3}]=\begin{bmatrix}
		\alpha \frac{\eta}{l} & \quad i(1-\beta) & 	\quad i(\beta-1) \\
		0  &  \quad	\omega_0 & \quad	\omega_0\\
		\beta & \quad i \alpha \frac{\eta}{l} &  \quad -i \alpha \frac{\eta}{l}
	\end{bmatrix},\quad \mathcal{J}:=\begin{bmatrix}
		e^{\lambda_1 t} &   & 	  \\
		&  e^{\lambda_2 t} &  \\
		&   &  e^{\lambda_3 t}
	\end{bmatrix},
\end{equation}
then we can calculate the adjugate matrix $\operatorname{adj}(\mathsf{S}_0)$ and determinant $\det(\mathsf{S}_0)$ of $\mathsf{S}_0$. Hence, we obtain that $e^{\mathsf L_0 t}=\frac{1}{\det(\mathsf{S}_0)} \mathsf{S}_0 \mathcal{J} \operatorname{adj}(\mathsf{S}_0)$, which combined with \eqref{7.xx} gives
\eqref{3.2.4}--\eqref{3.2.7}.
\end{proof}

Contrast to the classical Navier-Stokes equations, by means of \eqref{3.2.5} and \eqref{3.2.6} we find that although the combined effect of rotation and stratification brings about fluid oscillation, it cancels the lift-up effect. In other words, in the zero modes part, the oscillation effect caused by the rotation and stratification is significantly stronger than the lift-up effect caused by the Couette flow. 

Lemma~\ref{lem3.5} and the nonzero modes estimates in
appendix~\ref{secA} prove theorem~\ref{thm:linear}.

\section{The proof of theorem~\ref{1..1}}\label{section4}
\qquad In this section, we establish the global nonlinear stability stated in theorem~\ref{1..1} via a bootstrap argument. We begin by formulating the coupled nonlinear systems for the zero and nonzero modes. Taking the $x$-average of system  \eqref{2.11}, we obtain the zero modes system
\begin{equation}\label{eq:0modes}
\begin{cases}
(\p_t-\nu\Delta_{y,z})U_0^1+U_0\cdot\na U_0^1=(\beta-1)U_0^2 - (U_{\neq}\cdot\na_L U_{\neq}^1)_0,\\
(\p_t-\nu\Delta_{y,z})U_0^2+U_0\cdot\na U_0^2+\p_y P_0+\beta U_0^1=-(U_{\neq}\cdot\na_L U_{\neq}^2)_0,\\
(\p_t-\nu\Delta_{y,z})U_0^3+U_0\cdot\na U_0^3+\p_z P_0+\alpha\Theta_0=-(U_{\neq}\cdot\na_L U_{\neq}^3)_0,\\
(\p_t-\nu\Delta_{y,z})\Theta_0+U_0\cdot\na\Theta_0=\alpha U_0^3-(U_{\neq}\cdot\na_L\Theta_{\neq})_0,\\
\p_yU_0^2+\p_zU_0^3=0.
\end{cases}
\end{equation}
Denote $Z=(U_0^2,U_0^3,c_\beta U_0^1,\Theta_0)^T$. In terms of $Z$, the forcing from
nonzero modes on the right-hand side of \eqref{eq:0modes} can be written in divergence form as $\mathcal F_0=-\mathbb P_0\operatorname{div}_{y,z}
\mathcal S_0$, where the stress columns are given by
\begin{equation}\label{eq:stress}
\mathcal S_{0,j}=
\big(U_{\neq}^j(U_{\neq}^2,U_{\neq}^3,c_\beta U_{\neq}^1,
\Theta_{\neq})\big)_0,\qquad j=2,3,
\end{equation}
and $\mathbb P_0$ denotes the two-dimensional Leray projection acting on the first two components and the identity operator acting on the last two components.

The good unknowns \eqref{yyy} for nonzero modes satisfy the forced version of \eqref{Non0L}:
\begin{equation}\label{eq:nonzero-forced-system}
\begin{aligned}
\p_t\widehat Q_{\neq}+\nu p\widehat Q_{\neq}
={}&\left(\frac{\dot p}{4p}
 +\frac{l^2}{2p}\frac{\dot p_h}{p_h}\right)\widehat Q_{\neq}
 -\alpha\sqrt{\frac{p_h}{p}}\widehat H_{\neq}\\
&+\left(\sgn(\beta)\sqrt\B\,\frac{l}{\sqrt p}
 -2c_\beta^{-1}\frac{k^2}{p_h}\frac{l}{\sqrt p}\right)\widehat K_{\neq}
 +\widehat{\mathcal N_Q},\\
\p_t\widehat K_{\neq}+\nu p\widehat K_{\neq}
={}&-\left(\frac{\dot p}{4p}
 +\frac{l^2}{2p}\frac{\dot p_h}{p_h}\right)\widehat K_{\neq}
 -\sgn(\beta)\sqrt\B\,\frac{l}{\sqrt p}\widehat Q_{\neq}
 +\widehat{\mathcal N_K},\\
\p_t\widehat H_{\neq}+\nu p\widehat H_{\neq}
={}&\frac{\dot p}{4p}\widehat H_{\neq}
 +\alpha\sqrt{\frac{p_h}{p}}\widehat Q_{\neq}
 +\widehat{\mathcal N_H}.
\end{aligned}
\end{equation}
where 
\begin{equation}\label{Non0NL}
\begin{aligned}
\mathcal N_Q={}&|\na_L|^{3/2}|\na_{L,h}|^{-1}
                  (U\cdot\na_LU^3)_{\neq}
 +\p_z|\na_{L,h}|^{-1}|\na_L|^{-1/2}
                  (\p_i^LU^j\p_j^LU^i)_{\neq},\\
\mathcal N_K={}&ic_\beta|\na_L|^{1/2}|\na_{L,h}|^{-1}
  [\p_x(U\cdot\na_LU^2)_{\neq}
    -\p_y^L(U\cdot\na_LU^1)_{\neq}],\\
\mathcal N_H={}&|\na_L|^{1/2}(U\cdot\na_L\Theta)_{\neq}.
\end{aligned}
\end{equation}

For the zero modes, we  write $E_{0,s}=\frac{1}{2}\|Z\|_{H^s}^2$ and $F_{0,s}=\nu\|\nabla Z\|_{H^s}^2$. Let $L_s$ be the positive, time-independent weight constructed in lemma~\ref{lem:common-weight} of appendix~\ref{sec:zero-weight}. For $N\le s\le N+1$, define the weighted energy and dissipation by
\begin{equation}\label{eq:zero-energy-def}
 E_{0,s}^{w}=\frac{1}{2}\int L_s(\xi)|\widehat Z(\xi)|^2\,\mathrm d\xi,
 \qquad F_{0,s}^{w}=\nu\int|\xi|^2L_s(\xi)|\widehat Z(\xi)|^2\,\mathrm d\xi.
\end{equation}

By \eqref{eq:weight-properties}, there is $C_L=C_L(\alpha,B_\beta,N)\ge1$ such that
\begin{equation}\label{eq:energy-comparison}
 C_L^{-1}E_{0,s}\le E_{0,s}^{w}\le C_LE_{0,s},
 \qquad C_L^{-1}F_{0,s}\le F_{0,s}^{w}\le C_LF_{0,s}.
\end{equation}

Furthermore, we will use $M_7=M_7(t,k=0,\p_y)$ to obtain a equivalent energy which satisfies
\begin{equation}\label{eq:early-zero-def}
 \|M_7(t,0,\p_y)Z\|_{H^{N+1}}^2\asymp\|Z\|_{H^{N+1}}^2,
\end{equation}
with constants independent of $t$ and $\nu$.

For nonzero modes, we set
\begin{equation}\label{eq:nonzero-budget}
 \mathcal E_{\neq,N}(T)=\sup_{0\le t\le T}E_{\neq,N}(t)
                  +\int_0^T F_{\neq,N}(t)\,\mathrm dt.
\end{equation}

For fixed constant $A_*>0$ with $2\lambda A_*\ge2$ and $\lambda=\frac{1}{16}\frac{2\sqrt{B_{\beta}}-1}{2\sqrt{B_{\beta}}+1}$, we set
\begin{equation}\label{eq:switch-time}
 t_*=A_*\nu^{-1/3}(1+|\log\nu|).
\end{equation}

We now state the main bootstrap proposition as follows.

\begin{pro}\label{boots}
Fix a sufficiently large constant $\Boot$, depending only on $\alpha,\beta,a_0$ and $N$. Suppose that on $[0,T]$,  the solution satisfies the bootstrap ansatz
\begin{equation}\label{eq:full-bootstrap}
\begin{aligned}
 \mathcal E_{\neq,N}(T)&\le\Boot^2\ep^2,\\
 \sup_{t\le T}\|Z(t)\|_{H^{N+1}}^2
 +\nu\int_0^T\|\na Z\|_{H^{N+1}}^2\dd
 +\|\sqrt{\mu_7}M_7Z\|_{L^2(0,\min(T,t_*);H^{N+1})}^2&\le\Boot^2\ep^2.
\end{aligned}
\end{equation}
If $\ep\le\delta\nu^{2/3}/(1+|\log\nu|)^{1+a_0}$ with $\delta=\delta(\alpha,\beta,a_0,N,\Boot)>0$ sufficiently small and $a_0\in(0,1/2)$, then both
bounds in \eqref{eq:full-bootstrap} hold with $\Boot^2/2$ in place of $\Boot^2$.
\end{pro}

The proof relies on the following two propositions concerning the zero modes and nonzero modes. First, we state the energy estimates for the zero modes system.
\begin{pro}\label{prop:zero-forced}
Let $t_0\ge0$, $N\le s\le N+1$, and
$\sup_{N\le s\le N+1}\|Z(t_0)\|_{H^s}\lesssim\ep$.
Suppose that $\sup_{t_0\le t\le T}\|Z(t)\|_{H^{N+1}}\le\Boot\ep$, as the assumption in proposition~\ref{boots}. If $\ep^{3/2}/\nu$ is sufficiently small
with respect to the fixed bootstrap constant, then the corrected
zero modes energy gives
\begin{align}
 \sup_{t\in[t_0,T]} \|Z\|_{H^s}+\nu^{1/2}\|\na Z\|_{L^2_t([t_0,T])H^s}
 &\lesssim\ep+\|\mathcal F_0\|_{L^1_t([t_0,T])H^s},\label{eq:zero-close}\\
 \sup_{t\in[t_0,T]} \|Z\|_{H^s}+\nu^{1/2}\|\na Z\|_{L^2_t([t_0,T])H^s}
 &\lesssim\ep+\nu^{-1/2}\|\mathcal S_0\|_{L^2_t([t_0,T])H^s},
 \label{eq:zero-proposition}
\end{align}
In particular, \eqref{eq:zero-proposition} at $s=N+1$ improves the
highest-order bootstrap directly when its right-hand side is bounded
by a fixed constant times $\ep$ and $\Boot$ is chosen larger than
that constant.
\end{pro}
\begin{proof}
For a detailed proof of proposition~\ref{prop:zero-forced}, see section~\ref{sec:tame-zero}.
\end{proof}
Recall \eqref{Non0NL}, we denote the nonlinear contribution as
\begin{equation}\label{eq:source-pairing}
\begin{aligned}
\mathcal I_{\neq,N}
=\int \langle\xi\rangle^{2N}\big(|\widehat{\A Q_{\neq}}+G\widehat{\A K_{\neq}}|
       |\widehat{\A\mathcal N_Q}|+|\widehat{\A K_{\neq}}+G\widehat{\A Q_{\neq}}|
       |\widehat{\A\mathcal N_K}|
       +|\widehat{\A H_{\neq}}|\,|\widehat{\A\mathcal N_H}|\big){\rm d}\xi.
\end{aligned}
\end{equation}
In the following proposition, we give the estimate of nonzero modes.
\begin{pro}\label{prop:linear}
For any solution to \eqref{eq:nonzero-forced-system} on $[0,T]$, there holds
\begin{equation}\label{eq:linear-est}
E_{\neq,N}(T)+\int_0^TF_{\neq,N}\dd
\lesssim E_{\neq,N}(0)
 +\int_0^T\mathcal I_{\neq,N}(t)\dd.
\end{equation} 
Furthermore, under the bootstrap ansatz \eqref{eq:full-bootstrap}, we have
\begin{equation}\label{eq:nonlinear-proposition}
 \int_0^T\mathcal I_{\neq,N}(t)\dd
 \lesssim (1+|\log\nu|)^{1+a_0}\nu^{-2/3}\mathcal E_{\neq,N}^{3/2}
 +\nu^{-2/3}\sup_{t\le T}\|Z(t)\|_{H^{N+1}}\,\mathcal E_{\neq,N}.
\end{equation}
\end{pro}
\begin{proof}
Repeat the energy calculation \eqref{999} for
\eqref{eq:nonzero-forced-system}, with Sobolev weight
$\langle\xi\rangle^N$. The additional terms are precisely the
nonlinear contribution in \eqref{eq:source-pairing}. Taking the absolute values of the nonlinear contribution and absorbing the linear remainder as in \eqref{1010}, we obtain \eqref{eq:linear-est} by \eqref{00}. 
The estimate \eqref{eq:nonlinear-proposition} is proved in section~\ref{sec:nonzero-nonlinear}. 
\end{proof}

\noindent{\textit{Proof of proposition~\ref{boots}.}
We first verify the initial-energy bounds. For $k\ne0$, set
$p_0=k^2+\eta^2+l^2$ and $p_{h,0}=k^2+\eta^2$.}
The incompressibility  condition gives $l^2|\widehat{u_{\mathrm{in}}^3}|^2\le |k,\eta|^2\big(|\widehat{u_{\mathrm{in}}^1}|^2+|\widehat{u_{\mathrm{in}}^2}|^2\big)$,
and hence $|k,\eta,l|^2|\widehat{u_{\mathrm{in}}^3}|^2\le |k,\eta|^2|\widehat{u_{\mathrm{in}}}|^2$. Consequently, $|\widehat{Q_{\mathrm{in}}}|
 \le |k,\eta,l|^{1/2}|\widehat{u_{\mathrm{in}}}|$, $|\widehat{K_{\mathrm{in}}}|\le c_\beta |k,\eta,l|^{1/2}|\widehat{u_{\mathrm{in}}}|$ and $|\widehat{H_{\mathrm{in}}}|
 =|k,\eta,l|^{1/2}|\widehat{\theta_{\mathrm{in}}}|$. Since $M(0)=1$ and $m(0)=|k,\eta,l|^{1/2}$, these inequalities,
together with $L_s\asymp\langle\xi\rangle^{2s}$, imply
\[
 \|\mathcal A(0)(Q_{\mathrm{in}},K_{\mathrm{in}},
 H_{\mathrm{in}})\|_{H^N}^2+E_{0,N+1}^{w}(0)
 \lesssim
 \|u_{\mathrm{in}}\|_{H^{N+1}}^2
 +\|\theta_{\mathrm{in}}\|_{H^{N+1}}^2
 \le C_0\ep^2.
\]
Equations~\eqref{eq:linear-est} and \eqref{eq:nonlinear-proposition}
give
\begin{equation}\label{eq:nonzero-close}
 \mathcal E_{\neq,N}\le C_0\ep^2+
 C\big[(1+|\log\nu|)^{1+a_0}\nu^{-2/3}\mathcal E_{\neq,N}^{1/2}
       +\nu^{-2/3}\sup \|Z\|_{H^{N+1}}\big]\mathcal E_{\neq,N}.
\end{equation}

For zero modes, the short-time estimate \eqref{eq:early-zero-close}
has error coefficients
\[
 \ep (1+|\log\nu|)^{1+a_0}\nu^{-2/3}\le\delta,\qquad
 A_*\ep^2\nu^{-4/3}(1+|\log\nu|)\le A_*\delta^2(1+|\log\nu|)^{-1-2a_0}.
\]

Taking $\delta$ small depending on the fixed bootstrap constants
makes the early zero modes energy at most a fixed multiple of
its initial size, including the integral of
$\|\sqrt{\mu_7}M_7Z\|_{H^{N+1}}^2$.
If $T>t_*$, apply proposition~\ref{prop:zero-forced} starting at
$t_*$. The highest-order stress estimate \eqref{eq:tail-stress}
makes its forcing contribution small, while
$\ep^{3/2}/\nu\le\delta^{3/2}(1+|\log\nu|)^{-\frac{3}{2}(1+a_0)}$.
Thus \eqref{eq:late-zero-close} controls the remaining interval.
The constants from the comparison of the two energies are fixed
independently of $\nu$.

Choose $\Boot$ larger than the fixed initial and propagated energy
constants, then the normal-form cutoff constant $K$, and finally
$\delta$ small. The preceding bounds improve the zero modes
bootstrap, and \eqref{eq:nonzero-close} improves the nonzero modes
bootstrap. \hfill$\square$

\textit{Proof of theorem~\ref{1..1}.} By standard local well-posedness theory for parabolic systems, there exists a unique local strong solution on $[0,T_0]$ for some $T_0>0$. On every finite time interval, the solution is considered
in the class
\[
	(U,\Theta)\in C([0,T];H^N)\cap L^2(0,T;H^{N+1}).
\]
 The shear coordinate transformation \eqref{2.5} is bounded on these Sobolev spaces on each finite time interval. Since $N\ge4$, the $H^N$ bound controls the $W^{1,\infty}$ norm on such intervals. The usual Sobolev energy estimate therefore gives the continuation criterion needed below. The bootstrap estimates in proposition~\ref{boots} together with a standard continuity argument extend the solution globally in time to $[0,\infty)$, satisfying \eqref{eq:full-bootstrap} for all $T>0$. Applying these a priori bounds to smooth approximations and passing to the limit yield global existence at the stated regularity, while uniqueness follows from a standard $L^2$ energy estimate for the difference of two solutions.

Finally, the enhanced dissipation and inviscid damping estimate \eqref{eq:physical-decay} follows from the reconstruction bounds \eqref{3x.1}, \eqref{ll7} and \eqref{eq:full-bootstrap}, while the global stability bound \eqref{eq:main-bounds} and the higher-order zero-mode bound \eqref{eq:main-high} follow directly from \eqref{eq:zero-close}, \eqref{eq:early-zero-close}, \eqref{eq:late-zero-close}, and the viscous dissipation part of \eqref{eq:EF}. The proof of theorem~\ref{1..1} is thus complete.
\hfill$\square$

\section{Energy estimates on \texorpdfstring{$U_0$ and $\Theta_0$}{U0 and Theta0}}\label{sec:zero}
\allowdisplaybreaks[0]
\qquad In this section, we prove proposition~\ref{prop:zero-forced}. Throughout this section, we assume $0<\nu\le1$, $\alpha>0$, $B_\beta>1/4$, $B_\beta\ne\alpha^2$ and $N\ge4$, as in theorem~\ref{1..1}. In order to capture the oscillatory structure of the zero modes, we consider a linear inviscid toy model. It can be observed that the nonresonant part of the nonlinear contribution in the energy estimate can be rewritten as a cubic correction, together with the sum-of-squares energy, forming a normal form. Next, we restore the viscous term and the internal quadratic nonlinear terms for the toy model. Thus, we need to handle a near-resonance cubic term, a viscous cubic remainder, and a quartic remainder. The near-resonant cubic term requires the main structural
estimates. Finally, we restore the forcing from nonzero modes and close the estimates of zero modes system.

\subsection{Good unknowns and the cubic nonlinear contribution}
\qquad Set $Z=(U_0^2,U_0^3,c_\beta U_0^1,\Theta_0)^T$ and $v=(U_0^2,U_0^3)^T$. We temporarily set $U_{\neq}=\Theta_{\neq}=0$ in \eqref{eq:0modes}, and consider the following toy system
\begin{align}
	&(\p_t-\nu\Delta_{y,z})Z + \operatorname{diag}(\PP,I_2)\begin{pmatrix}
		0 & 0 & \operatorname{sgn}(\beta)\sqrt{B_{\beta}} & 0\\
		0 & 0 & 0 & \alpha \\
		-\operatorname{sgn}(\beta)\sqrt{B_{\beta}} & 0 & 0 & 0\\
		0 & -\alpha & 0 & 0
	\end{pmatrix}Z = \mathcal{N}(Z),\label{eq:zero-system}\\
	&\NL(Z)=-\big(\PP(v\cdot\na v),\ c_\beta v\cdot\na U_0^1,\ v\cdot\na\Theta_0\big),\label{eq:zero-vector}
\end{align}

where $\mathbb{P}$ is the two-dimensional Leray projection and the four components in $\NL(Z)$ share the divergence-free transport structure $v\cdot\na$.

For $\xi=(\eta,l)\ne0$, we recall the wave frequency in
section~\ref{sec3.2}
\begin{equation}\label{eq:omega}
\omega_0(\xi)=\sqrt{\frac{\alpha^2\eta^2+B_\beta l^2}{\eta^2+l^2}},
\qquad
\min(\alpha,\sqrt\B)\le\omega_0\le\max(\alpha,\sqrt\B).
\end{equation}

To capture the oscillatory structure of zero modes, we introduce three good unknowns.
\begin{defn}\label{def:5.1}
	We define
\begin{equation}\label{eq:normalized-amplitudes}
\begin{aligned}
\Gamma_{\sigma}(\xi)
&=\frac{1}{\sqrt2|\xi|}
 \left[l\widehat U_0^2-\eta\widehat U_0^3
 \ +\sigma\frac{i}{\omega_0(\xi)}
       (\beta l\widehat U_0^1-\alpha\eta\widehat\Theta_0)\right],\ \sigma\in\{-1,1\}\\
\Gamma_0(\xi)
&=\frac{\sgn(\beta)\alpha\eta c_\beta\widehat U_0^1
          +\sqrt\B\,l\widehat\Theta_0}{|\xi|\omega_0(\xi)},
\end{aligned}
\end{equation}
here $\Gamma_{\pm1}$ are the two wave amplitudes and $\Gamma_{0}$ is the stationary amplitude.
\end{defn}
For simplicity, we denote $\Gamma_{\pm 1}$ by $\Gamma_{\pm}$, respectively. One can see that
\begin{equation}\label{eq:zero-coefficients}
\Gamma_+(\xi)=-\overline{\Gamma_-(-\xi)},\qquad
\Gamma_0(-\xi)=-\overline{\Gamma_0(\xi)}.
\end{equation}
The following lemma shows that these amplitudes and $Z$ can be expressed with each other, here we omit its proof.
\begin{lem}\label{lem:physical-wave-expansion}
	Under the definition \ref{def:5.1}, we have 
	\begin{equation}\label{eq:physical-wave-expansion}
	\widehat Z=\sum_{\sigma=-1,0,1}\Gamma_\sigma\mathcal P_\sigma,\qquad \Gamma_\sigma=\mathcal P_\sigma^*\widehat Z
	,\ \sigma\in\{-1,0,1\},
	\end{equation}
where $*$ denotes conjugate transpose, and the reconstructed coefficient column vectors are
\begin{equation}\label{eq:polarizations}
\begin{aligned}
\mathcal P_\sigma(\xi)=\frac{1}{\sqrt2|\xi|}
\begin{pmatrix}
	l\\-\eta\\-i\sigma\beta l/(c_\beta\omega_0)\\ i\sigma\alpha\eta/\omega_0
\end{pmatrix},\ \sigma\in\{-1,1\},\quad
\mathcal P_0(\xi)=\frac{1}{|\xi|\omega_0}
\begin{pmatrix}0\\0\\\sgn(\beta)\alpha\eta\\\sqrt\B\,l\end{pmatrix},
\end{aligned}
\end{equation}
which form an orthonormal unit basis of the incompressible subspace
$\eta\widehat U_0^2+l\widehat U_0^3=0$. More precisely, the column vectors in \eqref{eq:polarizations} are odd with respect to $\xi$, and satisfy 
\[
	\mathcal P_+(\xi)=\overline{\mathcal P_-(\xi)},\qquad |\mathcal P_{\pm}(\xi)|^2=1,\qquad |\mathcal P_0(\xi)|^2=1,\qquad \mathcal{P}^{*}_{+}(\xi)\mathcal{P}_{-}(\xi)=0,\qquad \mathcal{P}^{*}_{\pm}(\xi)\mathcal{P}_{0}(\xi)=0,
\]
here $\mathcal P_{\pm 1}$ are represented by $\mathcal P_{\pm}$, respectively.
\end{lem}
According to the lemma~\ref{lem:physical-wave-expansion}, we can restore $U_0$ and $\Theta_0$ as
\begin{equation}\label{eq:zero-component-reconstruction}
\begin{aligned}
\widehat U_0^2&=\frac{l}{\sqrt2|\xi|}(\Gamma_-+\Gamma_+),&&
\widehat U_0^3=-\frac{\eta}{\sqrt2|\xi|}(\Gamma_-+\Gamma_+),\\
c_\beta\widehat U_0^1
&=\frac{i\beta l}{\sqrt2c_\beta|\xi|\omega_0}
  (\Gamma_--\Gamma_+)
 +\frac{\sgn(\beta)\alpha\eta}{|\xi|\omega_0}\Gamma_0, &&\widehat\Theta_0
=-\frac{i\alpha\eta}{\sqrt2|\xi|\omega_0}(\Gamma_--\Gamma_+)
 +\frac{\sqrt\B\,l}{|\xi|\omega_0}\Gamma_0.
\end{aligned}
\end{equation}
\begin{re}
	The double zero modes can still be expressed in terms of $\Gamma_{\pm}$ and $\Gamma_{0}$ as
\begin{equation}\label{eq:double-zero-amplitudes}
	\begin{aligned}
		c_{\beta}\widehat{\overline U_0^1}&=\operatorname{sgn}(\beta)\operatorname{sgn}(\eta)\Gamma_{0}(\eta,0), &\overline U_0^2&= 0,\\
		\widehat{\overline U_0^3}&=-\frac{\operatorname{sgn}(\eta)}{\sqrt{2}}\left(\Gamma_{-}(\eta,0)+\Gamma_{+}(\eta,0)\right), &\widehat{\overline\Theta_0}&=-i\frac{\operatorname{sgn}(\eta)}{\sqrt{2}}\left(\Gamma_{-}(\eta,0)-\Gamma_{+}(\eta,0)\right).
	\end{aligned}
\end{equation}
Note that the double zero modes system contains an oscillating
velocity--temperature pair with frequencies $\pm\alpha$ and a
stationary velocity component. Thus, we handle the entire zero modes together rather than decompose them.
\end{re}

Together with \eqref{eq:physical-wave-expansion}, we apply $\mathcal P^{*}_{\sigma}$ to the system \eqref{eq:zero-system} to obtain
\begin{equation}
\p_t\Gamma_\sigma(\xi)
=(-\nu|\xi|^2+i\sigma\omega_0(\xi))\Gamma_\sigma(\xi)
 +\mathcal P_\sigma^*\widehat{\NL(Z)},
\quad \sigma\in\{-1,0,1\}\label{eq:branch-equations},
\end{equation}
with
\begin{equation}
\mathcal P_\sigma^*\widehat{\NL(Z)}(\xi)=-\frac{i}{\sqrt2}
\sum_{\tau=\pm}\sum_{\kappa=-1,0,1}\int
 (e_{\xi'}\cdot(\xi-\xi'))
 \big[\mathcal P_\sigma(\xi)^*\mathcal P_\kappa(\xi-\xi')\big]\Gamma_\tau(\xi')\Gamma_\kappa(\xi-\xi')\,\mathrm d\xi',\label{eq:branch-quadratic-source}
\end{equation}
where $e_\xi=(l,-\eta)/|\xi|$. It is worth to mention that only wave amplitudes $\Gamma_{\pm}$ constitute the advecting $v=(U_0^2,U_0^3)^T$ by \eqref{eq:zero-component-reconstruction}, and \eqref{eq:branch-equations} is the key for the following analysis.

Write $\mathcal T_{0,s}=-\operatorname{Re}\langle Z,v\cdot\nabla Z\rangle_{H^s}$ for the usual Sobolev nonlinear contribution. We recall the energy $E_{0,N}=\frac{1}{2}\|Z\|_{H^N}^2$ and dissipation $F_{0,N}=\nu\|\nabla Z\|_{H^N}^2$. By Definition~\ref{def:5.1}, the orthonormality in lemma~\ref{lem:physical-wave-expansion} and \eqref{eq:zero-coefficients}, we have
\begin{equation}\label{eq:zero-three-branch-energy}
	\begin{aligned}
E_{0,N}&=\frac{1}{2}\sum_{\sigma=-1,0,1}
 \int\langle\xi\rangle^{2N}|\Gamma_\sigma(\xi)|^2\,\mathrm d\xi=\int\langle\xi\rangle^{2N}|\Gamma_-(\xi)|^2\,\mathrm d\xi
+\frac{1}{2}\int\langle\xi\rangle^{2N}|\Gamma_0(\xi)|^2\,\mathrm d\xi,\\ 
 F_{0,N}&=\nu\sum_{\sigma=-1,0,1}
 \int|\xi|^2\langle\xi\rangle^{2N}|\Gamma_\sigma(\xi)|^2\,\mathrm d\xi.
	\end{aligned}
\end{equation}
Differentiating \eqref{eq:zero-three-branch-energy} with
\eqref{eq:branch-equations} gives
\begin{equation}\label{eq:zero-row-energy}
\frac{{\rm d}}{{\rm d}t}E_{0,N}+F_{0,N}=-\operatorname{Re}\langle Z,v\cdot\na Z\rangle_{H^N}
=\operatorname{Re}\sum_{\sigma\in\{-1,0,1\}}\int\langle\xi\rangle^{2N}
 \overline{\Gamma_\sigma}\mathcal P_\sigma^*\widehat{\NL(Z)}\,\mathrm d\xi
=:\mathcal T_{0,N}.
\end{equation}
At $N=0$, the incompressibility condition implies $\mathcal{T}_{0,0}=0$. At $N>0$, thanks to the incompressibility condition, we can also preserve the cancellation formed by the commutator of derivatives. For simplicity, we introduce the notation $\int_{\xi_1+\xi_2+\xi_3=0}$ that contains the terms as $\sum_{l,l'}\iint_{\eta,\eta'}{\rm d}\eta{\rm d}\eta'$, which will be used in the following lemma.
\begin{lem}\label{lem:sobolev-cubic} 
	We can express $\mathcal T_{0,N}$ by the cubic terms composed of $\Gamma_{\sigma_j},\ \sigma_j\in\{-1,0,1\},\ j\in\{1,2,3\}$ and the cyclic coefficient as
\begin{equation}\label{eq:complete-cubic}
\begin{aligned}
\mathcal T_{0,N}
 &=\operatorname{Re}\sum_{\sigma_1,\sigma_2,\sigma_3}
 \int_{\xi_1+\xi_2+\xi_3=0}
 c_{N,\text{S}}\,\Gamma_{\sigma_1}(\xi_1)\Gamma_{\sigma_2}(\xi_2)\Gamma_{\sigma_3}(\xi_3),\\
c_{N,\text{S}}&=-\frac{i}{6\sqrt{2}}\sum_{\mathrm{cyc}(j,k,l)}
\sigma_j^2(e_{\xi_j}\cdot\xi_k)\,[\langle \xi_l\rangle^{2N}-\langle \xi_k\rangle^{2N}]
\mathcal P_{\sigma_k}(\xi_k)\cdot\mathcal P_{\sigma_l}(\xi_l),
\end{aligned}
\end{equation}
where the cyclic orders are $(1,2,3),(2,3,1),(3,1,2)$, and $c_{N,\text{S}}=c_{N,\text{S}}^{\sigma_1\sigma_2\sigma_3}(\xi_1,\xi_2,\xi_3)$.
\end{lem}
\begin{proof}
	Thanks to the Plancherel formula, we can rewrite $\mathcal T_{0,N}$ as
\begin{equation}\label{eq:physical-cubic-fourier}
\mathcal T_{0,N}
=-\operatorname{Re}\langle Z,v\cdot\na Z\rangle_{H^N}=\operatorname{Re}\iint
 -i \langle \xi\rangle^{2N}
 \big(\widehat v(\xi')\cdot(\xi-\xi')\big)
 \widehat Z(\xi-\xi')\cdot\overline{\widehat Z(\xi)}{\rm d}\xi'{\rm d}\xi.
\end{equation}
We denote $\xi_1=\xi'=(\eta',l')$, $\xi_2=\xi-\xi'=(\eta-\eta',l-l')$, $\xi_3=-\xi$.  Incompressible condition gives
\[
\widehat v(\xi_1)\cdot\xi_3
=-\widehat v(\xi_1)\cdot\xi_2.
\]
Interchanging positions 2 and 3 changes the wavenumber
from $\xi_2$ to $\xi_3$ yields
\begin{equation}\label{eq:paired-factor}
\begin{aligned}
\mathcal T_{0,N}
&=\operatorname{Re}\int_{\xi_1+\xi_2+\xi_3=0}
 -i(\widehat v(\xi_1)\cdot\xi_2)\langle\xi_3\rangle^{2N}
 \widehat Z(\xi_2)\cdot\widehat Z(\xi_3)\\
&=\operatorname{Re}\int_{\xi_1+\xi_2+\xi_3=0}
 +i(\widehat v(\xi_1)\cdot\xi_2)\langle\xi_2\rangle^{2N}
 \widehat Z(\xi_2)\cdot\widehat Z(\xi_3)\\
 &=\operatorname{Re}\int_{\xi_1+\xi_2+\xi_3=0}
-\frac{i}{2}(\widehat v(\xi_1)\cdot\xi_2)
 [\langle \xi_3\rangle^{2N}-\langle \xi_2\rangle^{2N}]
 \widehat Z(\xi_2)\cdot\widehat Z(\xi_3),
\end{aligned}
\end{equation}
where we have used $\overline{\widehat Z(\xi)}=\widehat Z(\xi_3)$ since $Z$ is real. Now we substitute \eqref{eq:physical-wave-expansion} into \eqref{eq:paired-factor} and rename the index pair $(1,2,3)$ cyclically with order $(1,2,3),(2,3,1),(3,1,2)$ that gives three representations of $\mathcal T_{0,N}$. Finally taking their average, we obtain \eqref{eq:complete-cubic}.
\end{proof}
Note that if $\sigma_1=\sigma_2=\sigma_3=0$, then every cyclic summand vanishes:
\[
c_{N,\text{S}}^{000}(\xi_1,\xi_2,\xi_3)=0.
\]
Therefore, \eqref{eq:complete-cubic} implies that $\mathcal T_{0,N}$ has no terms contributed by the product of three stationary amplitudes $\Gamma_{0}$.

We now bound the transport coefficient $c_{N,\text{S}}$ in \eqref{eq:complete-cubic}.
\begin{lem}\label{lem:cs}
	Let $r=\min_{j\in\{1,2,3\}}|\xi_j|$ and $R=\max_{j\in\{1,2,3\}}|\xi_j|$. Then the cyclic coefficient in \eqref{eq:complete-cubic} satisfies
\begin{equation}\label{eq:cs}
		|c_{N,\text{S}}^{\sigma_1,\sigma_2,\sigma_3}(\xi_1,\xi_2,\xi_3)|\lesssim r\langle R\rangle^{2N}.
	\end{equation}
\end{lem}
\begin{proof}
Since $\xi_1+\xi_2+\xi_3=0$, the two largest lengths lie between
$R/2$ and $R$. Recall the properties of $\mathcal P_{\sigma}$ in lemma~\ref{lem:physical-wave-expansion}, we always have $|\mathcal P_{\sigma_k}(\xi_k)\cdot\mathcal P_{\sigma_l}(\xi_l)|\leq 1$. Now we consider one order of the cyclic in \eqref{eq:complete-cubic}, which can be bounded by
\[
|e_{\xi_j}\cdot\xi_k|\,
 |\langle \xi_l\rangle^{2N}-\langle \xi_k\rangle^{2N}|.
\]
We consider the following three cases:
\\
(1) If $|\xi_j|=r$, the evenness of $\langle\xi\rangle^{2N}$ and $-\xi_l-\xi_k=\xi_j$ imply
\begin{align}\label{eq:cs-small-advector}
	|\langle \xi_l\rangle^{2N}-\langle \xi_k\rangle^{2N}|
	=\left|\int_0^1\na_\zeta\langle\zeta\rangle^{2N}\big|_{\zeta=\xi_k+\tau\xi_j}
	\cdot\xi_j\,\mathrm d\tau\right|
	\le C_N r\langle R\rangle^{2N-1}.
\end{align}
In fact, the length from origin to the point in the segment which joins $\xi_k$ to $-\xi_l$ is at most $R$.
Multiplying by $|e_{\xi_j}\cdot\xi_k|\le R$ proves the required bound.
\\
(2) If $|\xi_k|=r$, we directly have
\[
|e_{\xi_j}\cdot\xi_k|\le r,\qquad
|\langle \xi_l\rangle^{2N}-\langle \xi_k\rangle^{2N}|\le2\langle R\rangle^{2N}.
\]
(3) If $|\xi_l|=r$, the incompressibility condition $e_{\xi_j}\cdot\xi_j=0$ yields that
\begin{equation}\label{eq:cs-small-test}
e_{\xi_j}\cdot\xi_k=-e_{\xi_j}\cdot\xi_l,
\qquad |e_{\xi_j}\cdot\xi_k|\le r.
\end{equation}
Together with $|\langle \xi_l\rangle^{2N}-\langle \xi_k\rangle^{2N}|\le2\langle R\rangle^{2N}$ yield \eqref{eq:cs} and we finish the proof.
\end{proof}
Thanks to the incompressibility condition and the fact $\xi_1+\xi_2+\xi_3=0$, we obtain the following lemma which will be frequently used.
\begin{lem}\label{lem:2l}
	For $\xi_j=(\eta_j,l_j),\ j\in\{1,2,3\}$, if the coefficient $c_{N,\text{S}}$ in \eqref{eq:complete-cubic} is nonzero, then at least two of the $z$-direction wavenumbers $l_1,l_2,l_3$ will be nonzero.
\end{lem}
\begin{proof}
If $l_1=l_2=l_3=0$, then
\[
 e_{\xi_j}\cdot\xi_k
 =\frac{l_j\eta_k-\eta_jl_k}{|\xi_j|}=0,
\]
 every term in the complete cubic coefficient vanishes. The case of exactly one nonzero $z$-wavenumber is impossible because $l_1+l_2+l_3=0$. This  proves the lemma. The same argument applies when
the Sobolev weight is replaced by $L_s$.
\end{proof}
Lemma \ref{lem:2l} implies that $R=\max_{j\in\{1,2,3\}}|\xi_j|\ge1$, while $r=\min_{j\in\{1,2,3\}}|\xi_j|$ may touch zero. It also allows us to place one gradient at each of the two larger-frequency positions in \eqref{eq:complete-cubic} and obtain
\begin{equation}\label{eq:zero-direct-detailed}
	|\mathcal T_{0,N}|\lesssim
	\|Z\|_{H^N}\|\na Z\|_{H^N}^2
	\lesssim\nu^{-1}(E_{0,N})^{1/2}F_{0,N},
\end{equation}
which suggests a rough bound of threshold $O(\nu)$.
\subsection{Toy model and weighted energy}\label{sec:toy}
\qquad To capture the oscillatory structure of the cubic products in\eqref{eq:complete-cubic}, we first set $\nu=0$ and omit the nonlinear terms in
\eqref{eq:branch-equations}. The resulting linear inviscid toy system is
\begin{equation}\label{eq:inviscid-toy}
\p_t\Gamma_-=-i\omega_0\Gamma_-,\qquad
\p_t\Gamma_+=i\omega_0\Gamma_+,\qquad
\p_t\Gamma_0=0.
\end{equation} 
Consider, for instance, the product of two opposite-sign wave amplitudes and one stationary amplitude
\begin{equation}\label{eq:example-product}
\overline{\Gamma_-(\xi)}\Gamma_-(\xi')\Gamma_0(\xi-\xi').
\end{equation}
Under the toy dynamics \eqref{eq:inviscid-toy}, this product evolves according to
\begin{equation}\label{eq:toy-product-derivative}
\begin{aligned}
\p_t\big[\overline{\Gamma_-(\xi)}\Gamma_-(\xi')\Gamma_0(\xi-\xi')\big]
&=i\Psi^{(-1,0,1)}\overline{\Gamma_-(\xi)}\Gamma_-(\xi')\Gamma_0(\xi-\xi'),\\
\Psi^{(-1,0,1)}(\xi,\xi')&:=\omega_0(\xi)-\omega_0(\xi'),
\end{aligned}
\end{equation}
where $(-1,0,1)$ records the signed oscillatory branches in the three positions. More generally, for any $\sigma_1,\sigma_2,\sigma_3\in\{-1,0,1\}$, since $\omega_0$ is even, we have
\begin{equation}\label{eq:input-test-phase}
\begin{aligned}
\p_t\big[\overline{\Gamma_{-\sigma_3}(\xi)}
 \Gamma_{\sigma_1}(\xi')\Gamma_{\sigma_2}(\xi-\xi')\big]
&=i\Psi^{(\sigma_1,\sigma_2,\sigma_3)}
 \overline{\Gamma_{-\sigma_3}(\xi)}
 \Gamma_{\sigma_1}(\xi')\Gamma_{\sigma_2}(\xi-\xi'),\\
 \Psi^{(\sigma_1,\sigma_2,\sigma_3)}(\xi,\xi')
&:=\sigma_1\omega_0(\xi')+\sigma_2\omega_0(\xi-\xi')
 +\sigma_3\omega_0(\xi).
\end{aligned}
\end{equation}
Setting $(\xi_1,\xi_2,\xi_3)=(\xi',\xi-\xi',-\xi)$ and suppressing the branch indices when no confusion arises, we write the three-wave phase as
\begin{equation}\label{eq:phase}
\Psi=\sigma_1\omega_0(\xi_1)+\sigma_2\omega_0(\xi_2)+\sigma_3\omega_0(\xi_3),
\end{equation}

\begin{heur}[Construction of the adapted weight $L_s$]
For the linear system \eqref{eq:inviscid-toy}, whenever $\Psi\ne0$, the identity \eqref{eq:input-test-phase} yields
\[
 \operatorname{Re}\bigg(c_{N,\text{S}}\prod_{j=1}^3\Gamma_{\sigma_j}(\xi_j)\bigg)
 =\frac{\mathrm d}{\mathrm dt}\operatorname{Re}
 \bigg(\frac{c_{N,\text{S}}}{i\Psi}\prod_{j=1}^3\Gamma_{\sigma_j}(\xi_j)\bigg),
\]
which formally suggests dividing by $i\Psi$ to absorb nonresonant cubic interactions into a small energy correction. However, $\Psi$ may vanish, and the near-resonant terms with small $|\Psi|$ must be estimated separately. 

Among the 27 branch combinations, the interactions of three waves of nonuniform sign ($\sigma\in\{-1,1\}^3$ with $\max_{i<j}|\sigma_i-\sigma_j|>0$), require particular care. Up to relabeling and flipping all signs, such a phase takes the form
\[
 \Psi=\omega_1+\omega_2-\omega_3
 \ge2\min\{\alpha,\sqrt{B_\beta}\}-\max\{\alpha,\sqrt{B_\beta}\},
 \qquad \omega_j=\omega_0(\xi_j).
\]
Thus $\Psi$ can vanish only when $\max\{\alpha,\sqrt{B_\beta}\}\ge2\min\{\alpha,\sqrt{B_\beta}\}$, with $\alpha=2\sqrt{B_\beta}$ and $\sqrt{B_\beta}=2\alpha$ serving as the two critical endpoint ratios. To measure the size of the near-resonant set along one-dimensional frequency slices, we set $g_j=\partial_\eta\omega_0(\xi_j)$, fix $\xi_h,l_i,l_j$ for a permutation $\{i,j,h\}=\{1,2,3\}$, and view $\Psi$ as a function of $\eta_i$ via $\eta_j=-\eta_h-\eta_i$, so that $\p_{\eta_i}\Psi=\sigma_ig_i-\sigma_jg_j$. A
space-time resonance occurs when both the phase and its slice derivatives vanish simultaneously
\begin{equation}\label{eq:spacetimer}
	\Psi=0,\qquad \sigma_1g_1=\sigma_2g_2=\sigma_3g_3.
\end{equation}
The first condition in \eqref{eq:spacetimer} removes the time oscillation, while the second causes the first-derivative measure estimate for the near-resonant set to degenerate.

Although the commutator structure in \eqref{eq:paired-factor} provides the factor $\langle\xi_l\rangle^{2N}-\langle\xi_k\rangle^{2N}$, the standard Sobolev weight $\langle\xi\rangle^{2N}$ does not vanish at space-time resonances. To overcome this degeneracy, we seek an equivalent energy weight $L_s(\xi)\asymp \langle\xi\rangle^{2s}$ ($N\leq s\leq N+1$) satisfying the cancellation property
\[L_{N}(\xi_l)-L_{N}(\xi_k)\big|_{\text{space-time resonance}}=0.\] Observe that \eqref{eq:spacetimer} provides a natural invariant at space-time resonance, namely
\[
 |g_1|=|g_2|=|g_3|,\qquad \text{or equivalently,}\qquad |\p_{\eta}\omega_0(\xi_1)|=|\p_{\eta}\omega_0(\xi_2)|=|\p_{\eta}\omega_0(\xi_3)|.
\]
This motivates us to prescribe $L_s(\xi)$ as a function of $|\p_{\eta}\omega_0(\xi)|$ near the space-time resonant set. In the noncritical parameter regime ($2\alpha\neq \sqrt{B_{\beta}}$ and $2\sqrt{B_{\beta}}\neq \alpha$), lemma~\ref{lem:common-weight} shows that this common value is bounded away from zero and satisfies $|g_j|\asymp|\xi_j|^{-1}$ in a conic neighborhood $\mathcal{U}$ of the space-time resonant set. We therefore prescribe $L_s(\xi)=(1+|\partial_\eta\omega_0(\xi)|^{-2})^s$ on $\mathcal{U}$, which ensures that the weights coincide at space-time resonance while remaining comparable to $\langle\xi\rangle^{2s}$ and extend $L_s$ globally in the form
\begin{equation}\label{eq:weight-choice}
 L_s(\xi)=(1+|\xi|^2F(\eta^2/|\xi|^2))^s,\qquad L_s(0)=1,
 \qquad N\le s\le N+1,
\end{equation}
where $F\in C^{\infty}([0,1])$ is strictly positive and satisfies the bounds in \eqref{eq:weight-properties}. At the two critical endpoint ratios $\alpha=2\sqrt{B_{\beta}}$ and $\sqrt{B_{\beta}}=2\alpha$ (or whenever the space-time resonant set is empty), we simply take the standard weight $L_s(\xi)=\langle\xi\rangle^{2s}$ and estimate the near-resonant sets directly.  
\end{heur}
Setting $\omega_{*}=\min(\alpha,\sqrt{B_{\beta}})$ and fixing $0<\rho_0<\omega_{*}/4$, we  summarize the classification and treatment of all $3^3 = 27$ ordered cubic branch combinations in Table~\ref{tab:27-phases}. 
\begin{table}[H]
\centering
\small
\setlength{\tabcolsep}{5pt}
\renewcommand{\arraystretch}{1.15}
\begin{tabularx}{\linewidth}{@{}
 >{\raggedright\arraybackslash\strut}p{31mm}
 >{\raggedright\arraybackslash\strut}p{22mm}
 >{\raggedright\arraybackslash\strut}X@{}}
\toprule
Branch combination (number of cases) & Phase or coefficient & Treatment for $|\Psi|\le\rho<\rho_0$ \\
\midrule
$\Gamma_0^3$ (1)
 & $c_{N,\text{S}}^{000}=0$
 & Identically zero; no cubic correction is needed.\\[3pt]
$\Gamma_\pm\Gamma_0^2$ (6)
 & $|\Psi|\ge\omega_*$
 & \multirow{3}{=}{\centering{}Nonresonance: absorbed by the cubic correction $\mathcal{C}_{s,\vartheta}$.}\\[3pt]
$\Gamma_\pm^2\Gamma_0$ (6)
 & $|\Psi|\ge2\omega_*$
 & \\[3pt]
$\Gamma_\pm^3$ (2)
 & $|\Psi|\ge3\omega_*$
 & \\[3pt]
$\Gamma_-\Gamma_0\Gamma_+$ (6)
 & $\Psi$ may vanish
 & Nonresonance: absorbed by the cubic correction $\mathcal{C}_{s,\vartheta}$.

 Near-resonance:

 $l+l'\ne0$: the smallness of phase implies the smallness of   $c_s^{\mathrm{wws}}$.

 $l=-l'\ne0$: the smallness of phase and the cancellation between the phase and the weight difference imply the smallness of  $c_s^{\mathrm{wws}}$. \\[3pt]
$\Gamma_+^2\Gamma_-$, $\Gamma_-^2\Gamma_+$ (6)
 & $\Psi$ may vanish
 & Nonresonance: absorbed by the cubic correction $\mathcal{C}_{s,\vartheta}$.

 Near-resonance: 
 
If $\alpha=2\sqrt{B_\beta}$ or $\sqrt{B_\beta}=2\alpha$,
use the standard Sobolev weight and direct sublevel-set
estimates. If the space--time resonant set is empty,
use the standard weight and the sublevel-set estimates. Otherwise, construct an equivalent
energy whose weight differences are controlled by the
signed differences of $\partial_\eta\omega_0$.
\\
\bottomrule
\end{tabularx}
\caption{The 27 ordered cubic products, grouped by their treatment.}
\label{tab:27-phases}
\end{table}
We distinguish the standard Sobolev weight from the
nonstandard weight used near three-wave space--time resonance.
If $\sqrt{B_\beta}>\alpha$, the nonzero-frequency
three-wave space--time resonant set is empty, as shown
in appendix~\ref{sec:zero-weight}. If $\sqrt{B_\beta}<\alpha<2\sqrt{B_\beta}$, the mixed-sign
three-wave phases satisfy
\[
 |\Psi|\ge2\sqrt{B_\beta}-\alpha>0.
\] 
We use the standard Sobolev weight in the above two cases. At the critical ratios
$\alpha=2\sqrt{B_\beta}$ or $\sqrt{B_\beta}=2\alpha$,
we estimate the
near-resonant sets directly and also use the standard Sobolev weight.
Thus a nonstandard weight is needed only potentially in
the regime $\alpha>2\sqrt{B_\beta}$, when the relevant
space--time resonant set is nonempty. We now recall the weighted energy and dissipation $E_{0,s}^{w},F_{0,s}^{w}$ defined in \eqref{eq:zero-energy-def} as
\begin{equation}\label{eq:weighted-branch-energy}
 E_{0,s}^{w}=\frac{1}{2}\sum_{\sigma=-1,0,1}\int L_s(\xi)|\Gamma_\sigma(\xi)|^2\,\mathrm d\xi,\qquad F_{0,s}^{w}=\nu\sum_{\sigma=-1,0,1}\int|\xi|^2L_s(\xi)|\Gamma_\sigma(\xi)|^2\,\mathrm d\xi,
\end{equation}
here $E_{0,s}^{w}\asymp E_{0,s}$ and $F_{0,s}^{w}\asymp F_{0,s}$ as in \eqref{eq:energy-comparison}. For the system \eqref{eq:zero-system}, since $L_s$ is real, even and independent of time, the argument of lemma~\ref{lem:sobolev-cubic} gives
\begin{equation}\label{eq:weighted-cubic}
\begin{aligned}
 \frac{\mathrm d}{\mathrm dt}E_{0,s}^{w}+F_{0,s}^{w}
 &= -\operatorname{Re}\langle L_s(D)Z,v\cdot\nabla Z\rangle_{L^2}=:\mathcal T_{0,s}^{w}=\operatorname{Re}\sum_{\sigma_1,\sigma_2,\sigma_3}
 \int_{\xi_1+\xi_2+\xi_3=0}c_s\prod_{j=1}^3\Gamma_{\sigma_j}(\xi_j),\\
 c_s&=-\frac{i}{6\sqrt2}\sum_{\mathrm{cyc}(j,k,l)}
 \sigma_j^2(e_{\xi_j}\cdot\xi_k)[L_s(\xi_l)-L_s(\xi_k)]
 \mathcal P_{\sigma_k}(\xi_k)\cdot\mathcal P_{\sigma_l}(\xi_l).
\end{aligned}
\end{equation}
The proof of lemma~\ref{lem:cs} uses only evenness and the growth and gradient bounds of the weight. Hence \eqref{eq:weight-properties} gives
\begin{equation}\label{eq:weighted-cs}
 |c_s|\lesssim r\langle R\rangle^{2s}.
\end{equation}
Lemma~\ref{lem:2l} also applies to $c_s$, since its vanishing when $l_1=l_2=l_3=0$ follows from the advecting factor alone. Thus the preliminary coefficient estimates remain valid for $L_s$.

Let $P_j$, $j\in\mathbb Z$, be the projections associated to a standard homogeneous Littlewood--Paley decomposition on $\mathbb R\times\mathbb T$,
\[
f=\sum_{j\in\mathbb Z}P_jf,\qquad
\widehat{P_jf}(\xi)=\varphi_j(\xi)\widehat f(\xi),\qquad
\varphi_j(\xi)=\varphi(2^{-j}\xi),\quad \xi=(\eta,l),
\]
where $\varphi\in C_c^\infty(\mathbb R^2;[0,1])$ is radial, supported in $\{\xi:\frac{1}{2}\le|\xi|\le2\}$, and satisfies $\sum_j\varphi_j(\xi)=1$ for $\xi\ne0$. We write $Z_r=\sum_{|k-j|\le2}P_kZ$ for $r=2^j$ and $Z_{\sim R}=\sum_{2^{-6}R\le2^k\le2^6R}P_kZ$, with the same convention for other functions. Let
\[
 \varphi_{r,R}(\xi_1,\xi_2,\xi_3)
 =\sum_{\substack{\min_a2^{j_a}=r\\\max_a2^{j_a}=R}}
   \prod_{a=1}^3\varphi_{j_a}(\xi_a).
\]
Then $\sum_{r,R}\varphi_{r,R}=1$, and $\varphi_{r,R}$ is symmetric and supported where $\min_a|\xi_a|\asymp r$, $\max_a|\xi_a|\asymp R$. 
We now use \eqref{eq:input-test-phase} to absorb the nonresonant part of $\mathcal T_{0,N}^{w}$ into a cubic correction. We choose an even smooth cutoff function $0\le\chi\le1$ which satisfies
\[\chi(x)=\begin{cases}
	0,\ \text{if }|x|\leq 1,\\
	1,\ \text{if }|x|>2.
\end{cases}.\]
Fix $0<\vartheta\le1$, to be chosen later, and set $\chi_R=\chi(\Psi/(\vartheta R))$ on each block. Using \eqref{eq:toy-product-derivative} gives
\[
\begin{aligned}
 &\operatorname{Re}\left[\chi(\Psi^{(-1,0,1)}/(\vartheta R))c^{\mathrm{wws}}_N\,
 \overline{\Gamma_-(\xi)}\Gamma_{-}(\xi')\Gamma_{0}(\xi-\xi')\right]\\
 &\quad=\frac{\mathrm d}{\mathrm dt}
 \operatorname{Re}\left[\frac{\chi(\Psi^{(-1,0,1)}/(\vartheta R))c^{\mathrm{wws}}_N}{i\Psi^{(-1,0,1)}}\,
 \overline{\Gamma_-(\xi)}\Gamma_{-}(\xi')\Gamma_{0}(\xi-\xi')\right],
\end{aligned}
\]
where $c_{N}^{\mathrm{wws}}$ is the coefficient corresponding to the product \eqref{eq:example-product} in $\mathcal T_{0,N}^{w}$. For $N\le s\le N+1$, applying the same operation to every
term of \eqref{eq:weighted-cubic} and summing over the
frequency blocks gives the cubic correction
\begin{equation}\label{eq:correction}
 \mathcal C_{s,\vartheta}
 =-\operatorname{Re}\sum_{r,R}
 \sum_{\sigma_1,\sigma_2,\sigma_3}
 \int_{\xi_1+\xi_2+\xi_3=0}
 \varphi_{r,R}\frac{\chi_Rc_s}{i\Psi}
 \prod_{j=1}^3\Gamma_{\sigma_j}(\xi_j).
\end{equation}
where we have used $\Gamma_{\sigma}(-\xi)=-\overline{\Gamma_{-\sigma}(\xi)}$ by \eqref{eq:zero-coefficients}. We define $\chi_R/\Psi=0$ on $|\Psi|\le\vartheta R$. This correction absorbs the terms with coefficient $\varphi_{r,R}\chi_Rc_s$ in $\mathcal T_{0,s}^{w}$. The remaining terms, with coefficient $\varphi_{r,R}(1-\chi_R)c_s$, form the near-resonant contribution $\mathcal T_{\mathrm{near},s}$. We estimate it below together with the viscous and quartic remainders.

\subsection{Restoring viscosity and internal nonlinearity with
$U_{\neq}=\Theta_{\neq}=0$}\label{sec:normal-form}
\qquad We now restore $\nu>0$ and the nonlinear terms, then \eqref{eq:example-product} satisfies
\begin{align}
&\p_t\big[\overline{\Gamma_-(\xi)}\Gamma_-(\xi')\Gamma_0(\xi-\xi')\big]\nonumber\\
&=\big[i\Psi^{(-1,0,1)}
 -\nu(|\xi|^2+|\xi'|^2+|\xi-\xi'|^2)\big]
 \overline{\Gamma_-(\xi)}\Gamma_-(\xi')\Gamma_0(\xi-\xi')+\overline{\mathcal{P}_{-}^{*}\widehat{\mathcal N(Z)}(\xi)}\Gamma_-(\xi')\Gamma_0(\xi-\xi')\nonumber\\
&\quad+\overline{\Gamma_-(\xi)}\mathcal{P}_{-}^{*}\widehat{\mathcal N(Z)}(\xi')\Gamma_0(\xi-\xi')+\overline{\Gamma_-(\xi)}\Gamma_-(\xi')\mathcal{P}_{0}^{*}\widehat{\mathcal N(Z)}(\xi-\xi').\label{eq:worked-product}
\end{align}
Recall that
\begin{equation}\label{eq:block-cutoffs}
 \chi_R=\chi\left(\frac{\Psi}{\vartheta R}\right),
 \qquad 0<\vartheta\le1.
\end{equation}

Using the correction $\mathcal C_{s,\vartheta}$ defined in \eqref{eq:correction}, we  write the near-resonant
contribution as
\begin{equation}\label{eq:block-near-definition}
 \mathcal T_{\mathrm{near},s}=\operatorname{Re}\sum_{r,R}\sum_{\sigma_1,\sigma_2,\sigma_3}
 \int_{\sum_j\xi_j=0}\varphi_{r,R}(1-\chi_R)c_s
                 \prod_j\Gamma_{\sigma_j}(\xi_j).
\end{equation}

Applying \eqref{eq:worked-product} to each term of the correction and using \eqref{eq:weighted-cubic}, we obtain
\begin{equation}\label{eq:corrected-energy-identity}
 \frac{\mathrm d}{\mathrm dt}(E_{0,s}^{w}+\mathcal C_{s,\vartheta})+F_{0,s}^{w}
 =\mathcal T_{\mathrm{near},s}+\mathcal R_\nu+\mathcal R_4,
\end{equation}
where
\begin{align}
 \mathcal R_\nu&=\operatorname{Re}\sum_{r,R}\sum_{\sigma_1,\sigma_2,\sigma_3}
 \int_{\sum_j\xi_j=0}\varphi_{r,R}
       \nu\Big(\sum_j|\xi_j|^2\Big)
       \frac{\chi_Rc_s}{i\Psi}\prod_j\Gamma_{\sigma_j}(\xi_j),\label{eq:block-viscous-definition}\\
 \mathcal R_4&=-\operatorname{Re}\sum_{r,R}\sum_{\sigma_1,\sigma_2,\sigma_3}
 \int_{\sum_j\xi_j=0}\varphi_{r,R}\frac{\chi_Rc_s}{i\Psi}
 \sum_{j=1}^3[\mathcal P_{\sigma_j}^*\widehat{\NL(Z)}](\xi_j)\prod_{h\ne j}\Gamma_{\sigma_h}(\xi_h)\label{eq:three-quartic-insertions}.
\end{align}

The part with coefficient $\chi_Rc_s$ is absorbed by the cubic correction. For the remaining part $|\Psi|\le2\vartheta R$, if $2\vartheta R<\rho_0$, only the last two families in Table~\ref{tab:27-phases} contribute. If $2\vartheta R\ge\rho_0$, we use the general coefficient bound \eqref{eq:weighted-cs}, as shown below. We first record the coefficient of two opposite-sign waves and one stationary branch.
\begin{lem}\label{lem:2w1scoeff}
For $N\le s\le N+1$, the coefficient of
$\overline{\Gamma_-(\xi)}\Gamma_-(\xi')\Gamma_0(\xi-\xi')$ in \eqref{eq:weighted-cubic} is
\begin{equation}\label{eq:mixed-symbol}
 c_s^{\mathrm{wws}}
 =-\frac{\alpha\sqrt\B\,d^2}{2|\xi'||\xi||\xi-\xi'|\omega_0(\xi-\xi')}
 \left[\frac{L_s(\xi)-L_s(\xi-\xi')}{\omega_0(\xi)}
       -\frac{L_s(\xi')-L_s(\xi-\xi')}{\omega_0(\xi')}\right],
\end{equation}
where $d=\eta'l-\eta l'$.
The coefficient in \eqref{eq:weighted-cubic} satisfies
\[
 c_s^{\mathrm{wws}}(\xi,\xi')=-6c_s^{(-1,0,1)}(\xi',\xi-\xi',-\xi),
\]
where the factor $6$ counts the permutations and the minus sign follows from \eqref{eq:zero-coefficients}.
\end{lem}
\pf We first introduce some notations 
\begin{equation}\label{eq:mixed-geometry}
d=\eta'l-\eta l',\quad \mathfrak a(\xi,\xi')=
\frac{\alpha\sqrt\B\,d^2}{2|\xi'||\xi||\xi-\xi'|\omega_0(\xi-\xi')}.
\end{equation}
Incompressible condition gives
$e_{\xi'}\cdot (\xi-\xi')=-d/|\xi'|$ and $e_\xi\cdot (\xi-\xi')=-d/|\xi|$. And \eqref{eq:polarizations} yields
\begin{equation}\label{eq:polarization-products}
\mathcal P_\sigma(\zeta)\cdot\mathcal P_0(\xi-\xi')
=-\frac{i\sigma\alpha\sqrt\B\,(e_\zeta\cdot (\xi-\xi'))}{\sqrt2\omega_0(\zeta)|\xi-\xi'|\omega_0(\xi-\xi')},\qquad \zeta\in\{\xi,\xi'\},\ \sigma\in\{-1,1\}.
\end{equation}

\textit{Case I: Wave testing.} We test \eqref{eq:branch-equations} with the weight $L_s$ when $\sigma=\pm 1$ and $\kappa=0$. Using \eqref{eq:mixed-geometry} and  \eqref{eq:polarization-products}, the two opposite-sign waves and one stationary contribution can be written as
\begin{equation}\label{eq:wavetest}
\begin{aligned}
	-&\operatorname{Re}\iint\frac{2L_s(\xi)}{\sqrt{2}}i\big(e_{\xi'}\cdot(\xi-\xi')\big)\mathcal{P}_{-}(\xi)^*\mathcal{P}_{0}(\xi-\xi')\Gamma_{-}(\xi')\Gamma_{0}(\xi-\xi')\overline{\Gamma_{-}(\xi)}{\rm d}\xi'{\rm d}\xi\\
	&=-2\operatorname{Re}\iint\mathfrak a(\xi,\xi')
 \frac{L_s(\xi)}{\omega_0(\xi)}
 \overline{\Gamma_-(\xi)}\Gamma_-(\xi')\Gamma_0(\xi-\xi')
 \,\mathrm d\xi'\mathrm d\xi\\
 &=-\operatorname{Re}\iint\mathfrak a(\xi,\xi')
 \left(\frac{L_s(\xi)}{\omega_0(\xi)}
-\frac{L_s(\xi')}{\omega_0(\xi')}\right)
 \overline{\Gamma_-(\xi)}\Gamma_-(\xi')\Gamma_0(\xi-\xi')
 \,\mathrm d\xi'\mathrm d\xi,
\end{aligned}
\end{equation}

where in the last line we have exchanged $\xi,\xi'$ in the integral and used
\[
\overline{\Gamma_-(\xi')}\Gamma_-(\xi)\Gamma_0(\xi'-\xi)
=-\overline{\overline{\Gamma_-(\xi)}\Gamma_-(\xi')\Gamma_0(\xi-\xi')}.
\]

\textit{Case II: Stationary testing.}
We test \eqref{eq:branch-equations} with the weight $L_s$ when $\sigma=0$. The two opposite-sign waves and one stationary contribution can be written as
\begin{equation}\label{eq:0test}
	\begin{aligned}
	\operatorname{Re}&\int L_s(\xi-\xi')
\overline{\Gamma_0(\xi-\xi')}\mathcal{P}_{0}^{*}\widehat{\mathcal N}(\xi-\xi')\,\mathrm d(\xi-\xi')\\
	&=-\operatorname{Re}\frac{i}{\sqrt{2}}\sum_{\substack{\tau=\pm1\\\kappa=-\tau}}\int\int L_s(\xi-\xi')
\big(e_{\rho}\cdot(\xi-\xi'-\rho)\big)\mathcal{P}_{0}(\xi-\xi')^*\mathcal{P}_{\kappa}(\xi-\xi'-\rho)\\
&\qquad\qquad\qquad\qquad \qquad\qquad \times\overline{\Gamma_0(\xi-\xi')}\Gamma_{\tau}(\rho)\Gamma_{\kappa}(\xi-\xi'-\rho){\rm d}\rho\,\mathrm d(\xi-\xi').
	\end{aligned}
\end{equation}

If $(\tau,\kappa)=(-1,1)$, taking $\rho=\xi$ as in \eqref{eq:0test} and using $\mathcal P_+(-\xi')=-\overline{\mathcal P_-(\xi')}$ and
	\[\Gamma_{-}(\xi)\Gamma_{+}(-\xi')\overline{\Gamma_0(\xi-\xi')}=-\overline{\overline{\Gamma_{-}(\xi)}\Gamma_{-}(\xi')\Gamma_0(\xi-\xi')},
	\]
	then \eqref{eq:0test} becomes
\begin{equation}\label{eq:-11}
		-\operatorname{Re}\frac{i}{\sqrt{2}}\int\int L_s(\xi-\xi')
		\big(e_{\xi}\cdot(-\xi')\big)\mathcal{P}_0(\xi-\xi')^*\overline{\mathcal{P}_{-}(\xi')}\cdot\overline{\Gamma_{-}(\xi)}\Gamma_{-}(\xi')\Gamma_0(\xi-\xi'){\rm d}\xi\,\mathrm d(\xi-\xi').
	\end{equation}

	If $(\tau,\kappa)=(1,-1)$, taking $\rho=-\xi'$ as in \eqref{eq:0test} and using $e_{-\xi'}=-e_{\xi'}$ and
	\[\Gamma_{+}(-\xi')\Gamma_{-}(\xi)\overline{\Gamma_0(\xi-\xi')}=-\overline{\overline{\Gamma_{-}(\xi)}\Gamma_{-}(\xi')\Gamma_0(\xi-\xi')},\]
	then the coefficient in this part is  $i/\sqrt{2}\cdot \big(-e_{\xi'}\cdot \xi\big)\mathcal{P}_0(\xi-\xi')^*\mathcal{P}_{-}(\xi)$. Together with the coefficient parts in \eqref{eq:-11} gives
\begin{equation}\label{eq:wave-slow-slow-row}
\begin{aligned}
-&\frac{i}{\sqrt2}L_s(\xi-\xi')\big[
 (e_\xi\cdot(-\xi'))\mathcal P_0(\xi-\xi')^*\overline{\mathcal P_-(\xi')}+(e_{\xi'}\cdot\xi)\mathcal P_0(\xi-\xi')^*\mathcal P_-(\xi)\big]
 \\
 &=-\mathfrak a(\xi,\xi')L_s(\xi-\xi')\left[\frac{1}{\omega_0(\xi')}-\frac{1}{\omega_0(\xi)}\right],
\end{aligned}
\end{equation}
here $\mathfrak a$ is the same as \eqref{eq:mixed-geometry}.
Combining  \eqref{eq:mixed-geometry}, \eqref{eq:wavetest} and \eqref{eq:wave-slow-slow-row} gives \eqref{eq:mixed-symbol}.\hfill$\square$

We next estimate the cubic coefficient on $|\Psi|\le\rho$. The adapted weight from lemma~\ref{lem:common-weight} is used for the noncritical three-wave interactions.
The two critical ratios are handled separately by direct frequency-set estimates. The estimate in the following applies to arbitrary nonnegative amplitudes and is used
both for the near-resonant cubic term and for the quartic remainder in lemma~\ref{lem:nf-remainders}.
\begin{lem}\label{lem:phase-blocks}
Let $N\ge4$, $N\le s\le N+1$ and $\rho>0$. On a dyadic block with $|\xi_1|\asymp r$ and $|\xi_2|\asymp|\xi_3|\asymp R$, for every branch combination and nonnegative functions $f_1,f_2,f_3$ supported on the corresponding annuli, one has
\begin{equation}\label{eq:uniform-near-block}
 I_{r,R}(\rho):=\int_{\substack{\xi_1+\xi_2+\xi_3=0\\|\Psi|\le\rho}}
 \varphi_{r,R}\frac{|c_s|\prod_{j=1}^3f_j(\xi_j)}{\langle\xi_1\rangle^N\langle\xi_2\rangle^s\langle\xi_3\rangle^s|\xi_2||\xi_3|}\lesssim\left(\frac{\rho}{R}\right)^{1/2}
 \frac{\sqrt{r(1+r)}}{\langle r\rangle^N}\prod_{j=1}^3\|f_j\|_{L^2_\eta\ell^2_l}.
\end{equation}
The implicit constant depends only on the fixed parameters and $\alpha,\beta,N$.
\end{lem}
\begin{proof}
We assign the $H^N$ weight to the lowest-frequency factor after each relabeling. By lemma~\ref{lem:2l}, we may assume $R\gtrsim1$. Fix $0<\rho_0<\omega_*/4$ sufficiently small, depending only on $\alpha,B_\beta$. In the noncritical case, also require $\rho_0\le\rho_{\mathrm{sep}}(\alpha,B_\beta)$ from lemma~\ref{lem:common-weight}. For $\rho\ge\rho_0$, Young's and Cauchy--Schwarz inequalities give
\begin{align*}
 I_{r,R}(\rho)
 &\lesssim\frac{r}{R^2\langle r\rangle^N}\|f_1\|_{L^1_\eta\ell^1_l}\|f_2\|_{L^2_\eta\ell^2_l}\|f_3\|_{L^2_\eta\ell^2_l}
 \lesssim\frac{r^{3/2}\sqrt{1+r}}{R^2\langle r\rangle^N}\prod_{j=1}^3\|f_j\|_{L^2_\eta\ell^2_l}\\
 &\lesssim\rho_0^{-1/2}\left(\frac{\rho}{R}\right)^{1/2}
 \frac{\sqrt{r(1+r)}}{\langle r\rangle^N}\prod_{j=1}^3\|f_j\|_{L^2_\eta\ell^2_l},
\end{align*}
where we have used \eqref{eq:weighted-cs}, $r\lesssim R$ and $\sum_{|l|\lesssim r}\int_{|\eta|\lesssim r}\mathrm d\eta
 \lesssim r(1+r)$. Hence it remains to consider $0<\rho<\rho_0$. Since $c_s^{000}=0$ and $\omega_0\ge\omega_*>0$, only two opposite-sign waves with one stationary branch and three waves of nonuniform sign can contribute when $\rho_0<\omega_*$.

For $\xi_1+\xi_2+\xi_3=0$, the support of $\varphi_{r,R}$ is contained in
\[
 \mathcal B_{r,R}:=\left\{(\xi_1,\xi_2,\xi_3):
 \begin{gathered}
 \xi_j\in\mathbb R\times\mathbb Z,\quad \xi_1+\xi_2+\xi_3=0,\\
 r/2\le\min_j|\xi_j|\le2r,\qquad
 R/2\le\max_j|\xi_j|\le2R
 \end{gathered}\right\}.
\]
For any permutation $\{i,j,h\}=\{1,2,3\}$ and fixed $\xi_h=(\eta_h,l_h)$, $l_i,l_j\in\mathbb Z$ with $l_i+l_j+l_h=0$, set $\eta_j=-\eta_h-\eta_i$ and define the one-dimensional sets
\[
\begin{aligned}
 \mathcal S_{ij}&=\mathcal S_{ij}(\xi_h,l_i,l_j)
 :=\{\eta_i\in\mathbb R:(\xi_1,\xi_2,\xi_3)\in\mathcal B_{r,R}\},\\
 \mathcal S_{ij}(\rho)&:=\{\eta_i\in\mathcal S_{ij}:|\Psi(\xi_1,\xi_2,\xi_3)|\le\rho\}.
\end{aligned}
\]
Here $|\mathcal S_{ij}(\rho)|$ denotes one-dimensional Lebesgue measure. Let $K=K(\eta_i,\eta_h)$ be a measurable coefficient supported where $\eta_i\in\mathcal S_{ij}(\rho)$. Applying the Cauchy--Schwarz inequality first in $\eta_i$ and then in $\eta_h$, we obtain
\begin{equation}\label{eq:arbitrary-amplitudes}
\begin{aligned}
 &\iint |K|f_i(\eta_i)f_j(-\eta_h-\eta_i)f_h(\eta_h)
 \,\mathrm d\eta_i\mathrm d\eta_h\\
 &\quad\le\left(\sup_{\eta_h}\int_{\mathcal S_{ij}(\rho)}|K|^2\,\mathrm d\eta_i\right)^{1/2}
 \int f_h(\eta_h)
 \left(\int f_i(\eta_i)^2f_j(-\eta_h-\eta_i)^2\,\mathrm d\eta_i\right)^{1/2}\mathrm d\eta_h\\
 &\quad\le\left(\sup_{\eta_h}\int_{\mathcal S_{ij}(\rho)}|K|^2\,\mathrm d\eta_i\right)^{1/2}
 \prod_{j=1}^3\|f_j\|_{L^2_\eta}.
\end{aligned}
\end{equation}
The lowest frequency has $O(1+r)$ integer wavenumbers, so summation in the integer frequencies gives
\[
 \sum_{l_1+l_2+l_3=0}\prod_{j=1}^3\|f_j(\cdot,l_j)\|_{L^2_\eta}
 \le\left(\sum_{|l_1|\lesssim r}\|f_1(\cdot,l_1)\|_{L^2_\eta}\right)
 \prod_{j=2}^3\|f_j\|_{L^2_\eta\ell^2_l}
 \lesssim\sqrt{1+r}\prod_{j=1}^3\|f_j\|_{L^2_\eta\ell^2_l}.
\]

\textit{Case I: Two opposite-sign waves and one stationary branch.}
We use the variables $\xi,\xi'$ in lemma~\ref{lem:2w1scoeff} and write
\begin{equation}\label{eq:mixed-phase}
 \Psi=\omega_0(\xi)-\omega_0(\xi')
 =-\frac{(\alpha^2-\B)dj}{|\xi'|^2|\xi|^2[\omega_0(\xi')+\omega_0(\xi)]},
 \qquad j=\eta'l+\eta l'.
\end{equation}
Note that
\[
 d=\eta'l-\eta l'=\eta'(l-l')-l'(\eta-\eta')
 =-l(\eta-\eta')+\eta(l-l').
\]
Thus $|d|$ is the area of the parallelogram formed by any two of $\xi,\xi'$ and $\xi-\xi'$, and
\[
 |d|\le\min\{|\xi||\xi'|,|\xi'||\xi-\xi'|,|\xi||\xi-\xi'|\}
 \lesssim rR.
\]
Using $|\xi'||\xi||\xi-\xi'|\asymp rR^2$, $\omega_0\ge\omega_*$ and \eqref{eq:weight-properties} in \eqref{eq:mixed-symbol}, we obtain
\[
 \frac{|c_s^{\mathrm{wws}}|}{R^{2s+2}\langle r\rangle^N}
 \lesssim\frac{d^2}{rR^2}\frac{\langle R\rangle^{2s}}{R^{2s+2}\langle r\rangle^N}
 \lesssim\frac{d^2}{R^4r\langle r\rangle^N},
\]
where $R\gtrsim1$ by lemma~\ref{lem:2l}.
For fixed $q=\xi-\xi'$ and integer wavenumbers $l,l'$, set $\eta=\eta'+q_y$. With $(\xi_1,\xi_2,\xi_3)=(\xi',q,-\xi)$, the integration variable $\eta'$ belongs to $\mathcal S_{13}(\rho)$. The case $l=l'=0$ has zero coefficient by lemma~\ref{lem:2l}. If $l+l'\ne0$, then $j=(l+l')\eta'+l'q_y$. By the preceding bounds and \eqref{eq:arbitrary-amplitudes}, we have
\begin{align*}
 I_{r,R}(\rho)
 &\lesssim\sqrt{1+r}
 \left(\sup_{q_y}\int_{\mathcal S_{13}(\rho)}
 \frac{|c_s^{\mathrm{wws}}|^2}{R^{4s+4}\langle r\rangle^{2N}}\,\mathrm d\eta'\right)^{1/2}
 \prod_{j=1}^3\|f_j\|_{L^2_\eta\ell^2_l}\\
 &\lesssim\frac{\sqrt{1+r}}{R^4r\langle r\rangle^N}
 \left(\int_\R\min\{R^4r^4,\rho^4R^{16}|j|^{-4}\}\,\mathrm dj\right)^{1/2}
 \prod_{j=1}^3\|f_j\|_{L^2_\eta\ell^2_l}\\
 &\lesssim\frac{\sqrt{1+r}}{R^4r\langle r\rangle^N}
 \left(R^4r^4\int_0^{\rho R^3/r}\mathrm dj
       +\rho^4R^{16}\int_{\rho R^3/r}^{\infty}\frac{\mathrm dj}{j^4}\right)^{1/2}
 \prod_{j=1}^3\|f_j\|_{L^2_\eta\ell^2_l}\\
 &\lesssim\left(\frac{\rho}{R}\right)^{1/2}
 \frac{\sqrt{r(1+r)}}{\langle r\rangle^N}\prod_{j=1}^3\|f_j\|_{L^2_\eta\ell^2_l},
\end{align*}
where we have used $|\p_{\eta'}j|=|l+l'|\ge1$ and $|dj|\lesssim\rho R^4$.

If $l=-l'\ne0$, then $l,l'$ and $l-l'=2l$ are all nonzero,
so $r\gtrsim1$. In this case, we have $d=l(\eta+\eta')$ and $|d|\lesssim |l|R$. Equation~\eqref{eq:mixed-symbol} becomes
\begin{equation}\label{eq:mixed-split}
 c_s^{\mathrm{wws}}=\mathfrak a(\xi,\xi')
 \left[-\frac{L_s(\xi)-L_s(\xi')}{\omega_0(\xi)}
 +\Psi\frac{L_s(\xi')-L_s(\xi-\xi')}{\omega_0(\xi')\omega_0(\xi)}\right].
\end{equation}
By \eqref{eq:weight-properties}, \eqref{eq:mixed-phase}, $|d|\lesssim |l|R$ and  $|d|\lesssim rR$, we have
\begin{align*}
 |L_s(\eta,l)-L_s(\eta',-l)|
 &\le |\eta^2-\eta'^2|
 \sup_{\min\{\eta^2,\eta'^2\}\le u\le\max\{\eta^2,\eta'^2\}}
 |\partial_u L_s(\sqrt u,l)|\\
 &\lesssim\langle R\rangle^{2s-2}|\eta^2-\eta'^2|
 \lesssim\frac{|\Psi|R^{2s+2}}{l^2},\\
 \frac{|c_s^{\mathrm{wws}}|}{R^{2s+2}\langle r\rangle^N}
 &\lesssim\min\left\{|\Psi|
 \left(\frac{d^2}{l^2R^2r\langle r\rangle^N}
       +\frac{d^2}{R^4r\langle r\rangle^N}\right),
 \frac{d^2}{R^4r\langle r\rangle^N}\right\}\\
 &\lesssim\min\{\rho r^{-N-1},R^{-2}r^{1-N}\}.
\end{align*}
Since $\mathcal S_{13}(\rho)\subset[-2R,2R]$, we have $|\mathcal S_{13}(\rho)|\le4R$. Applying \eqref{eq:arbitrary-amplitudes}, we obtain
\begin{align*}
 I_{r,R}(\rho)
 &\lesssim R^{1/2}r^{1/2}
 \min\{\rho r^{-N-1},R^{-2}r^{1-N}\}\prod_{j=1}^3\|f_j\|_{L^2_\eta\ell^2_l}\\
 &\lesssim\left(\rho R^{1/2}r^{-N-1/2}\cdot R^{-3/2}r^{3/2-N}\right)^{1/2}
 \prod_{j=1}^3\|f_j\|_{L^2_\eta\ell^2_l}\\
 &=\left(\frac{\rho}{R}\right)^{1/2}r^{1/2-N}\prod_{j=1}^3\|f_j\|_{L^2_\eta\ell^2_l}
 \lesssim\left(\frac{\rho}{R}\right)^{1/2}
 \frac{\sqrt{r(1+r)}}{\langle r\rangle^N}\prod_{j=1}^3\|f_j\|_{L^2_\eta\ell^2_l}.
\end{align*}

\textit{Case II: Three waves of nonuniform sign, $\alpha\ne2\sqrt\B$ and $\sqrt\B\ne2\alpha$.}
Set $g_j=\p_\eta\omega_0(\xi_j)$ and $d_{ij}=\sigma_i g_i-\sigma_j g_j$. On the neighborhood $\mathcal U$ constructed in lemma~\ref{lem:common-weight}, we have $|\xi_j|\asymp R$, $|g_j|\asymp R^{-1}$ and $L_s(\xi_j)=(1+|g_j|^{-2})^s$. Hence
\begin{align*}
 |L_s(\xi_i)-L_s(\xi_j)|
 &\le\sup_{\min\{|g_i|,|g_j|\}\le u\le\max\{|g_i|,|g_j|\}}
 \left|\frac{\mathrm d}{\mathrm du}(1+u^{-2})^s\right|
 \big||g_i|-|g_j|\big|
 \lesssim R\langle R\rangle^{2s}|d_{ij}|,\\
 |c_s|&\lesssim rR\langle R\rangle^{2s}\sum_{i<j}|d_{ij}|,
\end{align*}
where we have used \eqref{eq:weighted-cubic} and $r\asymp R$ on $\mathcal U$. On the complement of $\mathcal U$, \eqref{eq:phase-gap-outside-weight} and \eqref{eq:weighted-cs} give
\[
 |c_s|\lesssim r\langle R\rangle^{2s}
 \lesssim c_{\mathrm{sep}}(\alpha,B_\beta)^{-1}rR\langle R\rangle^{2s}\sum_{i<j}|d_{ij}|.
\]
Combining the two regions with \eqref{eq:weighted-cs}, we obtain
\begin{equation}\label{eq:three-null}
 |c_s|\lesssim r\langle R\rangle^{2s}
 \min\bigg\{1,R\sum_{i<j}|d_{ij}|\bigg\},\qquad |\Psi|\le\rho_0.
\end{equation}
Writing $c_{ij}=c_s|d_{ij}|/\sum_{a<b}|d_{ab}|$, with $c_{ij}=0$ when the denominator vanishes, we have $c_s=\sum_{i<j}c_{ij}$ and
\[
 \frac{|c_{ij}|^2}{r^2\langle R\rangle^{4s}}
 \lesssim\min\{1,R|d_{ij}|\}^2\le R|d_{ij}|.
\]
To estimate the integral, we fix $\xi_h,l_i,l_j$, where $\{i,j,h\}=\{1,2,3\}$, and set $\eta_j=-\eta_h-\eta_i$. Along this one-dimensional set parametrized by $\eta_i$, the phase function
\[
 \Psi(\eta_i)=\sigma_i\omega_0(\eta_i,l_i)
 +\sigma_j\omega_0(-\eta_h-\eta_i,l_j)+\sigma_h\omega_0(\xi_h)
\]
satisfies $\partial_{\eta_i}\Psi=d_{ij}$ by the chain rule. First, we show that $\mathcal S_{ij}$ can be decomposed, up to finitely many endpoints, into a uniformly bounded number of intervals on which $\Psi$ is either constant or strictly monotone. When $l_i,l_j\ne0$, we write
\[
 \partial_{\eta_i}\Psi=A_i^{(1)}/\sqrt{B_i^{(1)}}-A_j^{(1)}/\sqrt{B_j^{(1)}},
\]
where, for $\kappa\in\{i,j\}$, $A_\kappa^{(1)}$ and $B_\kappa^{(1)}$ are polynomials in $\eta_i$ given by
\[
 A_\kappa^{(1)}=\sigma_\kappa(\alpha^2-B_\beta)\eta_\kappa l_\kappa^2,
 \qquad B_\kappa^{(1)}=(\alpha^2\eta_\kappa^2+B_\beta l_\kappa^2)(\eta_\kappa^2+l_\kappa^2)^3.
\]
The equation $\partial_{\eta_i}\Psi=0$ implies
\begin{equation}\label{eq:phase-slice-polynomial}
 (A_i^{(1)})^2B_j^{(1)}-(A_j^{(1)})^2B_i^{(1)}=0,
 \qquad \deg A_\kappa^{(1)}\le1,\quad\deg B_\kappa^{(1)}=8.
\end{equation}
If the polynomial in \eqref{eq:phase-slice-polynomial} is nontrivial, it has degree at most ten and hence at most ten real roots. If it vanishes identically, the two quotients have equal absolute values. Between the zeros of $A_i^{(1)}$ and $A_j^{(1)}$, their signs are fixed, so $\partial_{\eta_i}\Psi$ is either identically zero or twice one of the quotients, with a strict sign. If one of $l_i,l_j$ vanishes, the corresponding derivative vanishes and the remaining quotient has at most one zero. If both vanish, $\partial_{\eta_i}\Psi=0$ away from the zero vectors.

The boundary of $\mathcal S_{ij}$ is contained in the level sets $|\xi_\kappa|\in\{r/2,2r,R/2,2R\}$, $\kappa\in\{i,j\}$, each of which has at most two points unless it is empty. Partitioning $\mathcal S_{ij}$ at these boundary points and the real roots of the nontrivial polynomials just described, we obtain intervals $I_1,\ldots,I_J$, with $J\le C$ independent of the frequencies and $r,R$. On the interior of each interval, $\partial_{\eta_i}\Psi$ is smooth and either vanishes identically or maintains a strict sign. Neglecting finitely many endpoints and applying $u=\Psi(\eta_i)$ on each interval where $\partial_{\eta_i}\Psi\ne0$, we obtain
\[
 \int_{\mathcal S_{ij}(\rho)}
 R|\partial_{\eta_i}\Psi|\,\mathrm d\eta_i
 =R\sum_{m=1}^J
 \int_{I_m\cap\{|\Psi|\le\rho\}}|\partial_{\eta_i}\Psi|\,\mathrm d\eta_i=R\sum_{m=1}^J
 \int_{\Psi(I_m)\cap[-\rho,\rho]}\mathrm du
 \le2JR\rho\lesssim R\rho.
\]
If $\partial_{\eta_i}\Psi\equiv0$ on $I_m$, both integrals over $I_m\cap\{|\Psi|\le\rho\}$ and $\Psi(I_m)\cap[-\rho,\rho]$ vanish. Combining this bound with \eqref{eq:arbitrary-amplitudes} and the summation in the integer frequencies, we obtain
\begin{align*}
 I_{r,R}(\rho)
 &\lesssim\sum_{i<j}
 \left(\sup_{\eta_h}\int_{\mathcal S_{ij}(\rho)}
 \frac{|c_{ij}|^2}{R^{4s+4}\langle r\rangle^{2N}}\,\mathrm d\eta_i\right)^{1/2}
 \sqrt{1+r}\prod_{j=1}^3\|f_j\|_{L^2_\eta\ell^2_l}\\
 &\lesssim\frac{r\sqrt{1+r}}{R^2\langle r\rangle^N}
 \sum_{i<j}\left(\sup_{\eta_h}\int_{\mathcal S_{ij}(\rho)}R|\p_{\eta_i}\Psi|\,\mathrm d\eta_i\right)^{1/2}
 \prod_{j=1}^3\|f_j\|_{L^2_\eta\ell^2_l}\\
 &\lesssim\frac{r\sqrt{1+r}}{R^2\langle r\rangle^N}
 \left(R\int_{-\rho}^{\rho}\mathrm du\right)^{1/2}\prod_{j=1}^3\|f_j\|_{L^2_\eta\ell^2_l}\\
 &\lesssim\frac{\rho^{1/2}r\sqrt{1+r}}{R^{3/2}\langle r\rangle^N}\prod_{j=1}^3\|f_j\|_{L^2_\eta\ell^2_l}
 \lesssim\left(\frac{\rho}{R}\right)^{1/2}
 \frac{\sqrt{r(1+r)}}{\langle r\rangle^N}\prod_{j=1}^3\|f_j\|_{L^2_\eta\ell^2_l}.
\end{align*}

\textit{Case III: Three waves of nonuniform sign, $\alpha=2\sqrt\B$.}
After reordering and changing all signs if necessary, we have
\begin{align}
 \Psi&=\omega_0(\xi_1)+\omega_0(\xi_2)-\omega_0(\xi_3)\nonumber\\
 &=(\omega_0(\xi_1)-\sqrt\B)+(\omega_0(\xi_2)-\sqrt\B)
 +(2\sqrt\B-\omega_0(\xi_3))\nonumber\\
 &\ge\frac{3\sqrt\B}{4}
 \left(\frac{\eta_1^2}{|\xi_1|^2}
 +\frac{\eta_2^2}{|\xi_2|^2}+\frac{l_3^2}{|\xi_3|^2}\right)\ge0.
 \label{eq:critical-three-a}
\end{align}
For $(\xi_1,\xi_2,\xi_3)\in\mathcal B_{r,R}$ with $|\Psi|\le\rho$, \eqref{eq:critical-three-a} gives
\[
 |\eta_1|+|\eta_2|\lesssim R\sqrt\rho,\qquad |l_3|\lesssim\sqrt\rho\,|\xi_3|,\qquad |\xi_3|\le|\eta_1|+|\eta_2|+|l_3|
 \lesssim R\sqrt\rho+\sqrt\rho\,|\xi_3|.
\]
Reducing $\rho_0$ if necessary, we absorb the last term and obtain $|\xi_3|\lesssim R\sqrt\rho$. Hence $\xi_3$ is the lowest-frequency factor and $|\xi_3|\asymp r$.
For fixed $\xi_3,l_1,l_2$, every $\eta_1\in\mathcal S_{12}(\rho)$ satisfies $|\eta_1|\le CR\sqrt\rho$. Therefore,
\[
 \mathcal S_{12}(\rho)\subset[-CR\sqrt\rho,CR\sqrt\rho],
 \qquad |\mathcal S_{12}(\rho)|\le2CR\sqrt\rho.
\]
Using this set estimate in \eqref{eq:arbitrary-amplitudes}, together with \eqref{eq:weighted-cs}, $r\lesssim R\sqrt{\rho}$ and $R\geq 1$, we obtain
\begin{align*}
 I_{r,R}(\rho)
 \lesssim\frac{r}{R^2\langle r\rangle^N}(R\sqrt\rho)^{1/2}\sqrt{1+r}\prod_{j=1}^3\|f_j\|_{L^2_\eta\ell^2_l}\lesssim\left(\frac{\rho}{R}\right)^{1/2}\frac{\sqrt{r(1+r)}}{\langle r\rangle^N}\prod_{j=1}^3\|f_j\|_{L^2_\eta\ell^2_l}.
\end{align*}

\textit{Case IV: Three waves of nonuniform sign, $\sqrt\B=2\alpha$.}
Using \eqref{eq:omega} as in Case III, we obtain
\begin{equation}\label{eq:critical-three-b}
 \Psi\ge\frac{3\alpha}{4}
 \left(\frac{l_1^2}{|\xi_1|^2}+\frac{l_2^2}{|\xi_2|^2}
 +\frac{\eta_3^2}{|\xi_3|^2}\right)\ge0.
\end{equation}
For $(\xi_1,\xi_2,\xi_3)\in\mathcal B_{r,R}$ with $|\Psi|\le\rho$, \eqref{eq:critical-three-b} and $l_3=-l_1-l_2$ give
\[
 |l_1|+|l_2|\lesssim R\sqrt\rho,\qquad |\eta_3|\lesssim\sqrt\rho\,|\xi_3|,\qquad |\xi_3|\le|l_1|+|l_2|+|\eta_3|
 \lesssim R\sqrt\rho+\sqrt\rho\,|\xi_3|.
\]
As in Case III, reducing $\rho_0$ yields $|\xi_3|\asymp r\lesssim R\sqrt\rho$. Moreover, $|\eta_3|\lesssim\sqrt\rho\,|\xi_3|$ and $\xi_3\ne0$ imply $l_3\ne0$ and $r\gtrsim1$. For fixed $\xi_1,l_3,l_2$, with $\eta_2=-\eta_1-\eta_3$, every $\eta_3\in\mathcal S_{32}(\rho)$ satisfies $|\eta_3|\le Cr\sqrt\rho$. Thus
\[
 \mathcal S_{32}(\rho)\subset[-Cr\sqrt\rho,Cr\sqrt\rho],
 \qquad |\mathcal S_{32}(\rho)|\le2Cr\sqrt\rho.
\]
Using this set estimate in \eqref{eq:arbitrary-amplitudes}, $R\geq 1$ and $r\lesssim R\sqrt{\rho}$, we obtain
\begin{align*}
 I_{r,R}(\rho)
 \lesssim\frac{r}{R^2\langle r\rangle^N}(r\sqrt\rho)^{1/2}\sqrt{1+r}\prod_{j=1}^3\|f_j\|_{L^2_\eta\ell^2_l}\lesssim\left(\frac{\rho}{R}\right)^{1/2}\frac{\sqrt{r(1+r)}}{\langle r\rangle^N}\prod_{j=1}^3\|f_j\|_{L^2_\eta\ell^2_l}.
\end{align*}
Combining the above estimates proves \eqref{eq:uniform-near-block}.
\end{proof}

We now estimate the cubic correction, the viscous and quartic remainders, and the near-resonant term.
\begin{lem}[Cubic correction and remainder estimates]\label{lem:nf-remainders}
For $N\ge4$, $N\le s\le N+1$ and $0<\vartheta\le1$, one has
\begin{align}
 |\mathcal C_{s,\vartheta}|&\lesssim\vartheta^{-1}(E_{0,N}^{w})^{1/2}E_{0,s}^{w},\label{eq:C-size}\\
 |\mathcal R_\nu|&\lesssim\vartheta^{-1}(E_{0,N}^{w})^{1/2}F_{0,s}^{w},\label{eq:visc-rem}\\
 |\mathcal R_4|&\lesssim\frac{1}{\nu\sqrt\vartheta}
 \big((E_{0,N}^{w})^{1/2}(F_{0,s}^{w})^{1/2}+(F_{0,N}^{w})^{1/2}(E_{0,s}^{w})^{1/2}\big)^2,\label{eq:quartic-rem}\\
 |\mathcal T_{\mathrm{near},s}|&\lesssim\frac{\sqrt\vartheta}{\nu}(E_{0,N}^{w})^{1/2}F_{0,s}^{w}.\label{eq:near-final}
\end{align}
\end{lem}
\begin{proof}
Without loss of generality, we assume $|\xi_1|\asymp r$ and
$|\xi_2|\asymp|\xi_3|\asymp R$. We have
\begin{align*}
 |\mathcal C_{s,\vartheta}|
 &\le\sum_{r,R}\sum_{\sigma_1,\sigma_2,\sigma_3}
 \int_{\sum_j\xi_j=0}\varphi_{r,R}
 \left|\frac{\chi_Rc_s}{\Psi}\right|\prod_{j=1}^3|\Gamma_{\sigma_j}(\xi_j)|\\
 &\lesssim\vartheta^{-1}\sum_r\sum_{R\gtrsim\max(1,r)}
 \frac{r}{R\langle r\rangle^N}
 \left(\sum_{|l_1|\lesssim r}\int_{|\eta_1|\lesssim r}\mathrm d\eta_1\right)^{1/2}
 \|Z_r\|_{H^N}\|Z_R\|_{H^s}\|Z_{\sim R}\|_{H^s}\\
 &\lesssim\vartheta^{-1}
 \bigg(\sum_{r\le1}r^{3/2}+\sum_{r\ge1}r^{1-N}\bigg)\|Z\|_{H^N}
 \left(\sum_R\|Z_R\|_{H^s}^2\right)^{1/2}
 \left(\sum_R\|Z_{\sim R}\|_{H^s}^2\right)^{1/2}\\
 &\lesssim\vartheta^{-1}\|Z\|_{H^N}\|Z\|_{H^s}^2
 \lesssim\vartheta^{-1}(E_{0,N}^{w})^{1/2}E_{0,s}^{w},
\end{align*}
where we have used \eqref{eq:weighted-cs}, $|\chi_R/\Psi|\le(\vartheta R)^{-1}$, $R\gtrsim1$ by lemma~\ref{lem:2l}, and the finite overlap of the Littlewood--Paley projections. We also used $R\gtrsim\max\{1,r\}$, which implies
\[
 \frac{r^{3/2}\sqrt{1+r}}{R\langle r\rangle^N}
 \lesssim
 \begin{cases}
 r^{3/2},&r\le1,\\
 r^{1-N},&r\ge1.
 \end{cases}
\]
Similarly, since $|\xi_1|^2+|\xi_2|^2+|\xi_3|^2\lesssim R^2\asymp|\xi_2||\xi_3|$, we obtain
\begin{align*}
 |\mathcal R_\nu|
 &\le\sum_{r,R}\sum_{\sigma_1,\sigma_2,\sigma_3}
 \int_{\sum_j\xi_j=0}\varphi_{r,R}
 \frac{\nu\sum_j|\xi_j|^2}{|\xi_2||\xi_3|}
 \left|\frac{\chi_Rc_s}{\Psi}\right|
 |\Gamma_{\sigma_1}(\xi_1)|\prod_{j=2}^3|\xi_j|\,|\Gamma_{\sigma_j}(\xi_j)|\\
 &\lesssim\nu\vartheta^{-1}\sum_r\sum_{R\gtrsim\max(1,r)}
 \frac{r^{3/2}\sqrt{1+r}}{R\langle r\rangle^N}
 \|Z_r\|_{H^N}\|\nabla Z_R\|_{H^s}\|\nabla Z_{\sim R}\|_{H^s}\\
 &\lesssim\nu\vartheta^{-1}
 \left(\sum_{r\le1}r^{3/2}+\sum_{r\ge1}r^{1-N}\right)\|Z\|_{H^N}
 \left(\sum_R\|\nabla Z_R\|_{H^s}^2\right)^{1/2}
 \left(\sum_R\|\nabla Z_{\sim R}\|_{H^s}^2\right)^{1/2}\\
 &\lesssim\nu\vartheta^{-1}\|Z\|_{H^N}\|\nabla Z\|_{H^s}^2
 \lesssim\vartheta^{-1}(E_{0,N}^{w})^{1/2}F_{0,s}^{w},
\end{align*}
here we used the same dyadic bound. For $\mathcal R_4$, Young's inequality gives
\begin{align*}
 \|\NL(Z)\|_{H^s}
 \lesssim\|\widehat Z\|_{L^1_\eta\ell^1_l}\|\nabla Z\|_{H^s}
       +\|Z\|_{H^s}\||\xi|\widehat Z\|_{L^1_\eta\ell^1_l}\lesssim\|Z\|_{H^N}\|\nabla Z\|_{H^s}
       +\|\nabla Z\|_{H^N}\|Z\|_{H^s}.
\end{align*}
In every insertion, at least one unchanged factor has frequency comparable to $R$. We place one gradient on this factor and bound the other high-frequency length in \eqref{eq:uniform-near-block} by~$CR$. Multiplying \eqref{eq:uniform-near-block} by $R$, and using $|\mathcal P_\sigma|\le1$, $\int_{\vartheta R}^{\infty}\frac{R}{\rho^2}\left(\frac{\rho}{R}\right)^{1/2}\mathrm d\rho
 =2\vartheta^{-1/2}$ and $|\chi_R/\Psi|\le\int_{\vartheta R}^{\infty}\mathbf1_{\{|\Psi|\le\rho\}}\rho^{-2}\,\mathrm d\rho$, we obtain
\begin{align*}
 |\mathcal R_4|
 &\le\sum_{r,R}\sum_{\sigma_1,\sigma_2,\sigma_3}
 \int_{\sum_j\xi_j=0}\varphi_{r,R}\left|\frac{\chi_Rc_s}{\Psi}\right|
 \sum_{j=1}^3|\widehat{\NL(Z)}(\xi_j)|\prod_{h\ne j}|\Gamma_{\sigma_h}(\xi_h)|\\
 &\lesssim\sum_r\frac{\sqrt{r(1+r)}}{\langle r\rangle^N}
 \sum_{R\gtrsim\max(1,r)}
 \int_{\vartheta R}^{\infty}\frac{R}{\rho^2}\left(\frac{\rho}{R}\right)^{1/2}\mathrm d\rho\\
 &\qquad\times\Big[\|Z_r\|_{H^N}\|\nabla Z_R\|_{H^s}\|\NL(Z)_{\sim R}\|_{H^s}+\|\NL(Z)_r\|_{H^N}\|Z_R\|_{H^s}\|\nabla Z_{\sim R}\|_{H^s}\Big]\\
 &\lesssim\vartheta^{-1/2}
 \left(\sum_{r\le1}r^{1/2}+\sum_{r\ge1}r^{1-N}\right)
 \Big[\|Z\|_{H^N}\|\nabla Z\|_{H^s}\|\NL(Z)\|_{H^s}+\|\NL(Z)\|_{H^N}\|Z\|_{H^s}\|\nabla Z\|_{H^s}\Big]\\
 &\lesssim\vartheta^{-1/2}
 \big(\|Z\|_{H^N}\|\nabla Z\|_{H^s}
                  +\|\nabla Z\|_{H^N}\|Z\|_{H^s}\big)^2\\
 &\lesssim\frac{1}{\nu\sqrt\vartheta}
       \big((E_{0,N}^{w})^{1/2}(F_{0,s}^{w})^{1/2}
                       +(F_{0,N}^{w})^{1/2}(E_{0,s}^{w})^{1/2}\big)^2.
\end{align*}

Finally, $\operatorname{supp}(1-\chi_R)\subset\{|\Psi|\le2\vartheta R\}$. Taking $\rho=2\vartheta R$ in \eqref{eq:uniform-near-block}, which already includes the one-dimensional integrals and the summation in the integer frequencies, we obtain
\begin{align*}
 |\mathcal T_{\mathrm{near},s}|
 &\le\sum_{r,R}\sum_{\sigma_1,\sigma_2,\sigma_3}
 \int_{\sum_j\xi_j=0}\varphi_{r,R}|(1-\chi_R)c_s|
 \prod_{j=1}^3|\Gamma_{\sigma_j}(\xi_j)|\\
 &\lesssim\sqrt\vartheta\sum_r\frac{\sqrt{r(1+r)}}{\langle r\rangle^N}\|Z_r\|_{H^N}
 \sum_{R\gtrsim\max(1,r)}\|\nabla Z_R\|_{H^s}\|\nabla Z_{\sim R}\|_{H^s}\\
 &\lesssim\sqrt\vartheta
 \left(\sum_{r\le1}r^{1/2}+\sum_{r\ge1}r^{1-N}\right)\|Z\|_{H^N}
 \left(\sum_R\|\nabla Z_R\|_{H^s}^2\right)^{1/2}
 \left(\sum_R\|\nabla Z_{\sim R}\|_{H^s}^2\right)^{1/2}\\
 &\lesssim\sqrt\vartheta\|Z\|_{H^N}\|\nabla Z\|_{H^s}^2
 \lesssim\frac{\sqrt\vartheta}{\nu}(E_{0,N}^{w})^{1/2}F_{0,s}^{w}.
\end{align*}
This finishes the proof.
\end{proof}

\begin{re}
If $(E_{0,N}^{w})^{1/2}\le\Boot\ep$ and $C\Boot\ep/\vartheta\le1/4$,
then \eqref{eq:C-size} gives
\begin{equation}\label{eq:normal-form-equivalence}
 \frac{3}{4}E_{0,s}^{w}\le E_{0,s}^{w}+\mathcal C_{s,\vartheta}\le\frac{5}{4}E_{0,s}^{w}.
\end{equation}
\end{re}
Combining \eqref{eq:C-size}, \eqref{eq:visc-rem}, \eqref{eq:near-final}, \eqref{eq:quartic-rem} and \eqref{eq:corrected-energy-identity}, for $N\le s\le N+1$, we have
\begin{equation}\label{eq:zero-smallness}
\begin{aligned}
 \frac{\mathrm d}{\mathrm dt}(E_{0,s}^{w}+\mathcal C_{s,\vartheta})+F_{0,s}^{w}
 &\lesssim\left(\frac{\sqrt\vartheta (E_{0,N}^{w})^{1/2}}{\nu}
       +\frac{(E_{0,N}^{w})^{1/2}}{\vartheta}
       +\frac{E_{0,s}^{w}}{\nu\sqrt\vartheta}\right)F_{0,s}^{w}\\
 &\lesssim_{\Boot}\left(\frac{\ep\sqrt\vartheta}{\nu}
       +\frac{\ep}{\vartheta}+\frac{\ep^2}{\nu\sqrt\vartheta}\right)F_{0,s}^{w},
\end{aligned}
\end{equation}
where we have used $E_{0,N}^{w}\le E_{0,s}^{w}$, $F_{0,N}^{w}\le F_{0,s}^{w}$ and the bootstrap ansatz. Choosing
\begin{equation}\label{eq:cutoff-choice}
 \vartheta=K\ep,
\end{equation}
here $\ep$ is the size of the initial perturbation in \eqref{1.23}. Thus $\vartheta$ is fixed throughout the evolution and is independent of time. We obtain
\begin{equation}\label{eq:zero-final}
 \frac{\ep}{\vartheta}+\frac{\ep\sqrt\vartheta}{\nu}
       +\frac{\ep^2}{\nu\sqrt\vartheta}
 =\frac{1}K+\left(\sqrt K+\frac{1}{\sqrt K}\right)\frac{\ep^{3/2}}{\nu}\ll1.
\end{equation}
First fix $K$ sufficiently large, and then choose $c_K>0$ sufficiently small so that $\ep\le c_K\nu^{2/3}$ implies \eqref{eq:zero-final} and $\vartheta\le1$. Since $|\mathcal C_{s,\vartheta}|/E_{0,s}^{w}\lesssim_{\Boot}K^{-1}$, the corrected energy is equivalent to $E_{0,s}^{w}$ and the right-hand side of \eqref{eq:zero-smallness} is absorbed into $F_{0,s}^{w}$.

\subsection{Restoring the nonzero modes contribution}\label{sec:zero-forcing}
\qquad We now restore the forcing terms from the nonzero modes in \eqref{eq:0modes} and estimate their contribution via the Reynolds stress $\mathcal S_0$ defined in \eqref{eq:stress}, which requires the condition $B_\beta>1/4$ from lemma~\ref{lem3.4}. By the incompressibility condition, the nonzero modes forcing in \eqref{eq:0modes} can be written in divergence form as 
\[
\mathcal F_0=-\mathbb P_0(\partial_y\mathcal S_{0,2}+\partial_z\mathcal S_{0,3}),
\]
where $\mathbb P_0$ denotes the two-dimensional Leray projection on the first two components of $Z$ and the identity on the last two components. Here the stress columns $\mathcal S_{0,2}$ and $\mathcal S_{0,3}$ form the $4\times 2$ matrix
\begin{equation}\label{eq:eight-stresses}
\begin{pmatrix}
(U_{\neq}^2U_{\neq}^2)_0 & (U_{\neq}^3U_{\neq}^2)_0\\
(U_{\neq}^2U_{\neq}^3)_0 & (U_{\neq}^3U_{\neq}^3)_0\\
c_\beta(U_{\neq}^2U_{\neq}^1)_0 & c_\beta(U_{\neq}^3U_{\neq}^1)_0\\
(U_{\neq}^2\Theta_{\neq})_0 & (U_{\neq}^3\Theta_{\neq})_0
\end{pmatrix}.
\end{equation}
Consequently, the branch equation \eqref{eq:branch-equations} acquires the additional nonlinear term 
$\mathcal P_\sigma^*\widehat{\mathcal F_0}$. The double zero modes satisfy
\[
	\begin{cases}
		(\p_t-\nu\p_{yy})\overline U_0^1+\p_y(\overline{U^2_{\neq}U_{\neq}^1})_0=-\p_y(\overline{U^2_0\widetilde U_0^1}),\\
		(\p_t-\nu\p_{yy})\overline U_0^3+\alpha\overline\Theta_0
		+\p_y(\overline{U^2_{\neq}U_{\neq}^3})_0=-\p_y(\overline{U^2_0\widetilde U_0^3}),\\
		(\p_t-\nu\p_{yy})\overline\Theta_0-\alpha\overline U_0^3
		+\p_y(\overline{U^2_{\neq}\Theta_{\neq}})_0=-\p_y(\overline{U^2_0\widetilde \Theta_0}).
		\end{cases}
\]
Note that taking the $z$-average retains only the $\p_y \mathcal S_{0,2}$ derivative, which recovers the double zero modes forcing. 

Integration by parts yields
\begin{equation}\label{eq:IBPS}
	\operatorname{Re}\langle L_N(D)Z,\mathcal F_0\rangle_{L^2}
 =\sum_{j=2,3}\operatorname{Re}
   \langle L_N(D)\partial_jZ,\mathcal S_{0,j}\rangle_{L^2}.
\end{equation}
Here $L_N(D)$ commutes with derivatives and the Leray projection. Expanding the inner product for each $j\in\{2,3\}$ gives the four nonlinear contributions
\begin{equation}\label{eq:stress-row-pairings}
\begin{aligned}
&c_\beta^2\langle L_N(D)\p_jU_0^1,(U_{\neq}^jU_{\neq}^1)_0\rangle_{L^2},
 &&\langle L_N(D)\p_jU_0^2,(U_{\neq}^jU_{\neq}^2)_0\rangle_{L^2},\\
&\langle L_N(D)\p_jU_0^3,(U_{\neq}^jU_{\neq}^3)_0\rangle_{L^2},
 &&\langle L_N(D)\p_j\Theta_0,(U_{\neq}^j\Theta_{\neq})_0\rangle_{L^2},
\end{aligned}
\end{equation}
whose sum satisfies
$|\operatorname{Re}\langle L_N(D)Z,\mathcal F_0\rangle_{L^2}|\leq C\|\nabla Z\|_{H^N}\|\mathcal S_0\|_{H^N}
\lesssim (F_{0,N}^{w}/\nu)^{1/2}\|\mathcal S_0\|_{H^N}$.

For the cubic correction $\mathcal C_{N,\vartheta}$, differentiating in time produces the additional forcing remainder analogous to \eqref{eq:three-quartic-insertions}
\[
\mathcal{R}_{4}^{\neq}
 =-\operatorname{Re}\sum_{\sigma_1,\sigma_2,\sigma_3}
 \sum_{r,R}\int_{\xi_1+\xi_2+\xi_3=0}\varphi_{r,R}\frac{\chi_Rc_N}{i\Psi}
 \sum_{j=1}^3[\mathcal P_{\sigma_j}(\xi_j)^*\widehat{\mathcal F_0}(\xi_j)]
 \prod_{k\ne j}\Gamma_{\sigma_k}(\xi_k).
\]
In position $j$, the Fourier transform of the projected forcing reads
\[
\mathcal{P}_{\sigma_j}(\xi_j)^{*}\widehat{\mathcal{F}}_0(\xi_j)=-i\mathcal P_{\sigma_j}(\xi_j)^*\mathbb P_0(\xi_j)
 [\eta_j\widehat{\mathcal S}_{0,2}(\xi_j)
  +l_j\widehat{\mathcal S}_{0,3}(\xi_j)].
\]
If $k,m$ denote the two unchanged positions, then $\xi_j=-\xi_k-\xi_m$ implies $|\xi_j|\le|\xi_k|+|\xi_m|$, which allows us to transfer one derivative onto one of the two unchanged factors $\Gamma_{\sigma_k}$ or $\Gamma_{\sigma_m}$. Without loss of generality, assume that $|\xi_1|=\min_{a\in\{1,2,3\}}|\xi_a|\asymp r$. Assigning the $H^N$ weights to all three factors and using the bound $|c_N|/\prod_{a=1}^{3}\langle\xi_a\rangle^N\lesssim |\xi_1|/\langle\xi_1\rangle^N$ as in lemma~\ref{lem:nf-remainders}, we deduce from Young's convolution inequality that
\begin{equation}\label{eq:forced-C}
|\mathcal{R}_{4}^{\neq}|\lesssim \vartheta^{-1}\|\mathcal S_0\|_{H^N}\|\nabla Z\|_{H^N}\|Z\|_{H^N} \lesssim\vartheta^{-1}(E_{0,N}^{w})^{1/2}(F_{0,N}^{w}/\nu)^{1/2}\|\mathcal S_0\|_{H^N}.
\end{equation}
Combining this with the direct energy forcing bound and using 
$\vartheta^{-1}(E_{0,N}^{w})^{1/2}\ll1$ and Young's inequality, we obtain
\[
|\operatorname{Re}\langle L_N(D)Z,\mathcal F_0\rangle_{L^2}|+|\mathcal R_4^{\neq}|\leq C(1+\vartheta^{-1}(E_{0,N}^{w})^{1/2})(F_{0,N}^{w}/\nu)^{1/2}
 \|\mathcal S_0\|_{H^N}
 \le c_1F_{0,N}^{w}+C\nu^{-1}\|\mathcal S_0\|_{H^N}^2.
\]
Absorbing $c_1F_{0,N}^{w}$ and the terms $\mathcal T_{\mathrm{near},N}$, $\mathcal R_\nu$, $\mathcal R_4$ by \eqref{eq:zero-final} yields
\begin{equation}\label{eq:forced-low-diff}
\frac{\mathrm d}{\mathrm dt}(E_{0,N}^{w}+\mathcal C_{N,\vartheta})
+cF_{0,N}^{w}\lesssim\nu^{-1}\|\mathcal S_0\|_{H^N}^2 .
\end{equation}
Writing $X_{\neq}=(Q_{\neq},K_{\neq},H_{\neq})$ and recalling \eqref{eq:ghost-bounds}, \eqref{eq:M7-bounds}, \eqref{3x.1-linear}--\eqref{3x.2}, \eqref{ll6}--\eqref{ll7-linear}, and $mp^{1/4}(t,\xi)\gtrsim1$, we have the pointwise Fourier bound
\[
|\widehat{\mathcal{Z}}_{\neq}(\xi)|\lesssim
 \sqrt{1-(2\sqrt{B_\beta})^{-1}}\,
 \mathfrak L_{\alpha,\beta}|\widehat{\mathcal AX_{\neq}}(\xi)|,
 \qquad \mathcal{Z}\in\{U^1,U^2,U^3,\Theta,c_\beta U^1\},
\]
where $\kappa_{\beta}=(1-1/(2\sqrt{B_{\beta}}))^{-1}$ and 
\begin{equation}\label{eq:zero-forcing-parameter-factor}
	\mathfrak L_{\alpha,\beta}=\sqrt{\kappa_{\beta}}(1+c_\beta+c_\beta^{-1})\exp\left[\kappa_{\beta}\cdot\left(\frac{\pi(1+\alpha)+2}{\sqrt{B_\beta}}+\frac{4}{|\beta|}+\frac{2\pi}{c_\beta}\right)\right].
\end{equation}
The dependence on $a_0$ arising from $M_7^{-1}$ is included
in the implicit constants, whereas $\mathfrak L_{\alpha,\beta}$
records the dependence on the physical parameters.
Hence the nonzero modes estimates imply
\[
\begin{aligned}
\|\mathcal{Z}_{\neq}\|_{L^\infty H^N}
 &\lesssim\mathfrak L_{\alpha,\beta}\mathcal E_{\neq,N}^{1/2},\\
\|\mathcal{Z}_{\neq}\|_{L^2_tH^N}
 &\lesssim\mathfrak L_{\alpha,\beta}\nu^{-1/6}\mathcal E_{\neq,N}^{1/2},\\
\|\na_L\mathcal{Z}_{\neq}\|_{L^2_tH^N}
 &\lesssim\mathfrak L_{\alpha,\beta}\nu^{-1/2}\mathcal E_{\neq,N}^{1/2}.
\end{aligned}
\]
The $H^N$  estimate gives
\begin{align}\label{eq:stress-entry-low}
	\|(U_{\neq}^j \mathcal{Z}_{\neq})_0\|_{L^2_tH^N} \lesssim \|U_{\neq}^j\|_{L^\infty H^N}\|\mathcal{Z}_{\neq}\|_{L^2_tH^N}+\|U_{\neq}^j\|_{L^2_tH^N}\|\mathcal{Z}_{\neq}\|_{L^\infty H^N}\lesssim \mathfrak L_{\alpha,\beta}^2\nu^{-1/6}\mathcal E_{\neq,N}.
\end{align}
Summing over the entries of $\mathcal{S}_0$ yields
\begin{equation}\label{eq:stress-low}
\|\mathcal S_0\|_{L^2(0,T;H^N)}
\lesssim\mathfrak L_{\alpha,\beta}^2\nu^{-1/6}\mathcal E_{\neq,N}(T).
\end{equation}
Integrating \eqref{eq:forced-low-diff} and using energy equivalence
at the endpoints, we have
\[
\sup_{t\le T}E_{0,N}^{w}(t)+\int_0^TF_{0,N}^{w}\dd
 \lesssim E_{0,N}^{w}(0)+\nu^{-1}\|\mathcal S_0\|_{L^2_tH^N}^2.
\]
\eqref{1.23} gives
$E_{0,N}^{w}(0)\lesssim b_\beta^2\varepsilon^2,\ b_\beta=\max\{1,c_\beta\}$, and we obtain
\begin{equation}\label{eq:low-feedback}
\sup_{t\le T}E_{0,N}^{w}(t)+\int_0^TF_{0,N}^{w}\dd
\lesssim b_\beta^2\ep^2
 +\mathfrak L_{\alpha,\beta}^4\nu^{-4/3}\mathcal E_{\neq,N}(T)^2,
\end{equation}

\subsubsection{High-order zero modes estimates}\label{sec:tame-zero}
\begin{lem}\label{lem:forced-high}
Let $N\le s\le N+1$. Under the bootstrap assumptions \eqref{eq:full-bootstrap} and the parameter choice \eqref{eq:cutoff-choice}--\eqref{eq:zero-final}, for every $0\le t_0\le T$,
\begin{equation}\label{eq:tame-result}
 \sup_{t_0\le t\le T}E_{0,s}^{w}(t)+\int_{t_0}^TF_{0,s}^{w}\dd
 \lesssim E_{0,s}^{w}(t_0)+\nu^{-1}\|\mathcal S_0\|_{L^2(t_0,T;H^s)}^2.
\end{equation}
\begin{equation}\label{eq:direct-high-result}
 \sup_{t_0\le t\le T}\|Z(t)\|_{H^s}+\left(\int_{t_0}^TF_{0,s}^{w}\dd\right)^{1/2}
 \lesssim \|Z(t_0)\|_{H^s}+\|\mathcal F_0\|_{L^1(t_0,T;H^s)}.
\end{equation}
\end{lem}
\begin{proof}
We write $\mathcal R_{4,s}$ and $\mathcal R_{4,s}^{\neq}$ for the quartic and forcing remainders at order $s$. By $L_N\le L_s$ and \eqref{eq:quartic-rem},
\begin{equation}\label{eq:tame-quartic}
 |\mathcal R_{4,s}|\lesssim\frac{1}{\nu\sqrt\vartheta} E_{0,s}^{w}F_{0,s}^{w}
 \lesssim_{\Boot}\frac{\ep^2}{\nu\sqrt\vartheta}F_{0,s}^{w}.
\end{equation}
Together with \eqref{eq:visc-rem} and \eqref{eq:near-final}, the smallness condition \eqref{eq:zero-final} allows us to absorb $\mathcal T_{\mathrm{near},s}$, $\mathcal R_\nu$ and $\mathcal R_{4,s}$ into the dissipation $F_{0,s}^{w}$.

Repeating the estimate of $\mathcal C_{s,\vartheta}$ in lemma~\ref{lem:nf-remainders} with arbitrary functions $f_j$, and using $\langle\xi\rangle^N\le\langle\xi\rangle^s$, we obtain
\[
 \sum_{r,R}\sum_{\sigma_1,\sigma_2,\sigma_3}\int_{\sum_j\xi_j=0}
 \varphi_{r,R}\left|\frac{\chi_Rc_s}{\Psi}\right|
 \prod_{j=1}^3|\widehat f_j(\xi_j)|
 \lesssim \vartheta^{-1}\prod_{j=1}^3\|f_j\|_{H^s}.
\]
In each of the three forcing terms in $\mathcal R_{4,s}^{\neq}$, we use $\widehat{\mathcal F_0}(\xi_j)=-i\mathbb P_0(\xi_j)
 [\eta_j\widehat{\mathcal S}_{0,2}(\xi_j)+l_j\widehat{\mathcal S}_{0,3}(\xi_j)]$.
Using $|\xi_j|\le\sum_{h\ne j}|\xi_h|$, we obtain
\begin{equation}\label{eq:tame-force}
\begin{aligned}
 |\mathcal R_{4,s}^{\neq}|
 &\lesssim\sum_{r,R}\sum_{\sigma_1,\sigma_2,\sigma_3}\sum_{j=1}^3
 \int_{\sum_i\xi_i=0}\varphi_{r,R}
 \left|\frac{\chi_Rc_s}{\Psi}\right|
 |\widehat{\mathcal S_0}(\xi_j)|
 \sum_{h\ne j}|\xi_h|\prod_{i\ne j}|\Gamma_{\sigma_i}(\xi_i)|\\
 &\lesssim \vartheta^{-1}\|Z\|_{H^s}\|\na Z\|_{H^s}\|\mathcal S_0\|_{H^s}.
\end{aligned}
\end{equation}
Combining \eqref{eq:tame-force} with \eqref{eq:IBPS} at order $s$, and using $\vartheta^{-1}\|Z\|_{H^s}\ll1$, we have
\[
 C(1+\vartheta^{-1}\|Z\|_{H^s})\|\na Z\|_{H^s}\|\mathcal S_0\|_{H^s}
 \le c_1F_{0,s}^{w}+C\nu^{-1}\|\mathcal S_0\|_{H^s}^2.
\]
Integrating on $[t_0,T]$ and using \eqref{eq:normal-form-equivalence} proves \eqref{eq:tame-result}.

Alternatively, applying the same trilinear estimate to $\mathcal F_0$ without transferring its divergence derivative gives
\begin{equation}\label{eq:force-insertions}
 |\mathcal R_{4,s}^{\neq}|\lesssim \vartheta^{-1} \|Z\|_{H^s}^2\|\mathcal F_0\|_{H^s}.
\end{equation}
Together with the direct energy estimate, this gives
\[
 \frac{\mathrm d}{\mathrm dt}(E_{0,s}^{w}+\mathcal C_{s,\vartheta})+cF_{0,s}^{w}
 \lesssim (1+\vartheta^{-1}\|Z\|_{H^s})\|Z\|_{H^s}\|\mathcal F_0\|_{H^s}
 \lesssim \|Z\|_{H^s}\|\mathcal F_0\|_{H^s}.
\]
By \eqref{eq:normal-form-equivalence},
\begin{align*}
 \frac{\mathrm d}{\mathrm dt}(E_{0,s}^{w}+\mathcal C_{s,\vartheta})^{1/2}
 &\lesssim\|\mathcal F_0\|_{H^s},\\
 \sup_{t_0\le t\le T}\|Z(t)\|_{H^s}
 &\lesssim \|Z(t_0)\|_{H^s}+\|\mathcal F_0\|_{L^1(t_0,T;H^s)},\\
 \int_{t_0}^TF_{0,s}^{w}\dd
 &\lesssim \|Z(t_0)\|_{H^s}^2+\sup_{t_0\le t\le T}\|Z(t)\|_{H^s}
           \|\mathcal F_0\|_{L^1(t_0,T;H^s)}\\
 &\lesssim\big(\|Z(t_0)\|_{H^s}+\|\mathcal F_0\|_{L^1(t_0,T;H^s)}\big)^2.
\end{align*}
This proves \eqref{eq:direct-high-result}. Then we know that \eqref{eq:direct-high-result} and \eqref{eq:tame-result} give \eqref{eq:zero-close} and \eqref{eq:zero-proposition}, respectively. Thus proposition~\ref{prop:zero-forced} follows.
\end{proof}
\begin{lem}\label{lem:stress-high}
Under the assumptions of lemma~\ref{lem:forced-high} and \eqref{1.23}, for $0\le\sigma\le1$,
\begin{align}
	\|\mathcal S_0\|_{L^2(0,T;H^{N+\sigma})}&\lesssim\mathfrak L_{\alpha,\beta}^2\nu^{-1/6-\sigma/3}\mathcal E_{\neq,N}(T),\label{eq:stress-high}\\
	\sup_{t\le T}(E_{0,N+\sigma}^{w})^{1/2} + \left(\int_0^TF_{0,N+\sigma}^{w}(t){\rm d}t\right)^{1/2}
&\lesssim b_\beta\ep
 +\mathfrak L_{\alpha,\beta}^2\nu^{-(2+\sigma)/3}\mathcal E_{\neq,N}(T).\label{eq:high-budget}
\end{align}
\end{lem}

\begin{proof}
For $q=\xi'+\xi''$ with $k'+k''=0$, one has
\[
 \langle q\rangle^N\lesssim\langle\xi'\rangle^N+\langle\xi''\rangle^N,
 \qquad \langle q\rangle^\sigma\lesssim (p')^{\sigma/2}+(p'')^{\sigma/2}.
\]
For $\mathcal Z\in\{U^1,U^2,U^3,\Theta,c_\beta U^1\}$, Cauchy--Schwarz and interpolation give
\begin{align*}
 \|p^{\sigma/2}\widehat{\mathcal Z}_{\neq}\|_{L^1_\eta\ell^1_{k,l}}
 &\le\left(\sum_{k,l}\int_\R\langle\xi\rangle^{-2N}\,\mathrm d\eta\right)^{1/2}
       \||\na_L|^\sigma \mathcal Z_{\neq}\|_{H^N}
 \lesssim\||\na_L|^\sigma \mathcal Z_{\neq}\|_{H^N},\\
 \||\na_L|^\sigma \mathcal Z_{\neq}\|_{L^2_tH^N}
 &\le\|\mathcal Z_{\neq}\|_{L^2_tH^N}^{1-\sigma}
        \|\na_L\mathcal Z_{\neq}\|_{L^2_tH^N}^{\sigma}
 \lesssim\mathfrak L_{\alpha,\beta}\nu^{-1/6-\sigma/3}\mathcal E_{\neq,N}^{1/2}.
\end{align*}
Thus, we have 
\begin{align}\label{eq:stress-entry-high}
	\|(U_{\neq}^j \mathcal Z_{\neq})_0\|_{L^2H^{N+\sigma}}\lesssim &\||\na_L|^\sigma U_{\neq}^j\|_{L^2_tH^N}\|\mathcal Z_{\neq}\|_{L^\infty H^N}+\|U_{\neq}^j\|_{L^\infty H^N}\||\na_L|^\sigma \mathcal Z_{\neq}\|_{L^2_tH^N}\nonumber\\
	\lesssim &\mathfrak L_{\alpha,\beta}^2\nu^{-1/6-\sigma/3}\mathcal E_{\neq,N}.
\end{align}
Summing the stress entries proves \eqref{eq:stress-high}. Substitution in \eqref{eq:tame-result} proves \eqref{eq:high-budget}.
\end{proof}
Note that at the top regularity level $\sigma=1$ (i.e., $H^{N+1}$), the stress bound \eqref{eq:high-budget} yields a forcing contribution of order $\nu^{-1}\mathcal E_{\neq,N}\lesssim \ep^2\nu^{-1}$, which exceeds $O(\ep)$ under the threshold $\ep\lesssim \nu^{2/3}(1+|\log\nu|)^{-1-a_0}$ in \eqref{1.23} and thus cannot close the $H^{N+1}$ bootstrap on $[0,T]$ directly. On the other hand, the direct forcing bound \eqref{eq:direct-high-result} avoids transferring a derivative onto $Z$, but bounding $\|\mathcal{F}_0\|_{L^1_tH^{N+1}}$ when the advecting nonzero mode has high frequency near its critical time requires the cross-frequency transfer property $\eqref{eq:single-transfer}_3$ of the multiplier $M_7$.

To overcome this obstacle, we adopt a two-stage strategy divided at the switching time $t_{*}=A_{*}\nu^{-1/3}(1+|\log\nu|)$ defined in \eqref{eq:switch-time}:
\\
\textit{(1) Short-time stage $[0,\min(T,t_{*})]$.} In section~\ref{sec:zero-early}, we estimate the $M_7(t,0,\p_y)$-weighted zero modes energy $\frac12\|M_7Z\|_{H^{N+1}}^2$. The dissipation $\mu_7(t,0,\p_y):=-\dot{M_7}/M_7$ on the zero modes pairs with that of the nonzero modes via $\eqref{eq:single-transfer}_3$ to control the top-order forcing $\langle M_7Z,M_7\mathcal{F}_0\rangle_{H^{N+1}}$, while over the finite time horizon $t_{*}$, the internal zero-mode nonlinearity 
$\mathcal{N}(Z)$ can be absorbed directly by viscous dissipation without normal-form corrections.
\\
\textit{(2) Long-time stage $[t_{*},T]$.} In section~\ref{sec:zero-late}, for $t\geq t_{*}$, the exponential decay in time yields the factor $e^{-2\lambda\nu^{1/3}t_{*}}\leq\nu^2$, which makes the remaining Reynolds-stress contribution small. We can then switch back to the modified energy $E_{0,N+1}^w+\mathcal{C}_{N+1,\vartheta}$ and apply lemma~\ref{lem:forced-high} from $t_0=t_{*}$ to $T$.

\subsubsection{The highest-order forcing estimate on $[0,\min(T,t_{*})]$}
\label{sec:zero-early}
\qquad Recall $X_{\neq}=(Q_{\neq},K_{\neq},H_{\neq})$. By \eqref{eq:EF}--\eqref{eq:Equiv2} and \eqref{eq:enhanced-rate}, the total nonzero-mode energy functional $\mathcal{E}_{\neq,N}(T)$ controls
\begin{equation}\label{eq:nonzero-norm-budget}
\begin{aligned}
 &\|\mathcal AX_{\neq}\|_{L^\infty(0,T;H^N)}^2
 +\nu\|\na_L\mathcal AX_{\neq}\|_{L^2(0,T;H^N)}^2\\
 &\quad+\nu^{1/3}\|\mathcal AX_{\neq}\|_{L^2(0,T;H^N)}^2
 +\big\|\sqrt{\muEight}\,\mathcal AX_{\neq}\big\|_{L^2(0,T;H^N)}^2
 \lesssim\mathcal E_{\neq,N}(T).
\end{aligned}
\end{equation}

By \eqref{3.20} and \eqref{4.11-2}, we have
\[
 |(\widehat U_{\neq},\widehat\Theta_{\neq})|
 \lesssim e^{-\lambda\nu^{1/3}t}\frac{|\widehat{\mathcal AX_{\neq}}|}{mp^{1/4}}.
\]
Consequently
\begin{equation}\label{eq:low-force-point}
\begin{aligned}
 \|\mathcal F_0(t)\|_{H^N}
 &\lesssim e^{-2\lambda\nu^{1/3}t}
       \|\na_L\mathcal AX_{\neq}\|_{H^N}\|\mathcal AX_{\neq}\|_{H^N},\\
 \|\mathcal F_0\|_{L^1(0,T;H^N)}
 &\lesssim\nu^{-2/3}\mathcal E_{\neq,N}(T).
\end{aligned}
\end{equation}

We now estimate the highest-order forcing in the zero modes energy.
\begin{lem}\label{lem:top-forcing}
The highest-order forcing satisfies
\begin{equation}\label{eq:top-force-point}
\begin{aligned}
 \big|\operatorname{Re}\langle M_7Z,M_7\mathcal F_0\rangle_{H^{N+1}}\big|\lesssim{}&e^{-2\lambda\nu^{1/3}t}\|Z\|_{H^{N+1}}
    \|\mathcal AX_{\neq}\|_{H^N}
    \big[\|\na_L\mathcal AX_{\neq}\|_{H^N}
           +(1+t)\|\mathcal AX_{\neq}\|_{H^N}\big]\\
 &+e^{-2\lambda\nu^{1/3}t}(1+t)^2[\log(e+1+t)]^{1+a_0}
     \|\mathcal AX_{\neq}\|_{H^N}\\
 &\qquad\times\big\|\sqrt{\muEight}\,\mathcal AX_{\neq}\big\|_{H^N}
     \big\|\sqrt{\mu_7}M_7(t,0,D_y)Z\big\|_{H^{N+1}}.
\end{aligned}
\end{equation}
For $0<t_1\le T$, one also has
\begin{equation}\label{eq:top-force-integral}
\begin{aligned}
 \int_0^{t_1}\big|\operatorname{Re}\langle M_7Z,M_7\mathcal F_0\rangle_{H^{N+1}}\big|\dd
 \lesssim{}&(1+|\log\nu|)^{1+a_0}\nu^{-2/3}\mathcal E_{\neq,N}(t_1)
 \Big[\|Z\|_{L^\infty(0,t_1;H^{N+1})}\\
 &\qquad+\big\|\sqrt{\mu_7}M_7(t,0,D_y)Z
                 \big\|_{L^2(0,t_1;H^{N+1})}\Big].
\end{aligned}
\end{equation}
\end{lem}

\begin{proof}
Let $\xi'=(k,\eta',l')$ be the velocity frequency and
$\xi''=(-k,\eta'',l'')$ the transported frequency, with $k\ne0$.
Their sum $q=\xi'+\xi''$ is a zero mode. We write $p'=p(t,\xi')$, $m'=m(t,\xi')$, $p''=p(t,\xi'')$ and $m''=m(t,\xi'')$. Since $q=\xi'_L+\xi''_L$,
\[
 \langle q\rangle\lesssim\sqrt{p'}+\sqrt{p''},\qquad
 \sqrt{p'}\le(1+t)\langle\xi'\rangle,\qquad
 1+t\lesssim\langle\xi'\rangle\sqrt{p'}.
\]
The last inequality follows from $|k|\ge1$ and
\[
 |k|t\le|\eta'|+|\eta'-kt|
 \le2\langle\xi'\rangle\sqrt{p'}.
\]
The incompressibility condition $\xi'_L\cdot\widehat U(\xi')=0$ gives
\[
 |\widehat U^2(\xi')|\le\frac{|k,l'|}{\sqrt{p'}}|\widehat U(\xi')|,
\]
and hence
\[
 |\widehat U(\xi')\cdot\xi''_L|
 \le\left(|(k,l'')|
    +\frac{|(k,l')|}{\sqrt{p'}}|\eta''+kt|\right)
       |\widehat U(\xi')|.
\]
Since the Leray projection $\mathbb P_0$ is bounded, we have
\begin{align*}
 &\big|\operatorname{Re}\langle M_7Z,M_7\mathcal F_0\rangle_{H^{N+1}}\big|\lesssim\int_{q=(0,q_y,q_z)}\langle q\rangle^{N+1}|M_7(q)\widehat Z(q)|
       \int_{\xi'+\xi''=q,\ k\ne0}
       \langle q\rangle^{N+1}M_7(q)\\
 &\hspace{50mm}\times
       |\widehat U(\xi')\cdot\xi''_L|
       |(\widehat U,\widehat\Theta)(\xi'')|.
\end{align*}
Using $\langle q\rangle^{N+1}\lesssim\langle q\rangle(\langle\xi'\rangle^N+\langle\xi''\rangle^N)$ and the reconstruction formula \eqref{3.20}, we estimate the resulting multiplier
\[
 \frac{\langle q\rangle}{m'(p')^{1/4}m''(p'')^{1/4}}
 \left(|k,l''|+\frac{|k,l'|}{\sqrt{p'}}|\eta''+kt|\right).
\]

\smallskip\noindent\textit{Case I: $\langle\xi'\rangle\ge\langle\xi''\rangle$.}
By \eqref{4.11-2},
\begin{align*}
 \frac{\langle q\rangle |(k,l'')|}{m'(p')^{1/4}m''(p'')^{1/4}}
 &\lesssim\frac{\sqrt{p'}+\sqrt{p''}}{m'(p')^{1/4}}
 \lesssim(1+t)\langle\xi''\rangle,\\
 \frac{\langle q\rangle |(k,l')|\,|\eta''+kt|}{m'(p')^{1/4}m''(p'')^{1/4}\sqrt{p'}}
 &\lesssim\left(1+\frac{\sqrt{p''}}{\sqrt{p'}}\right)
                \frac{|\eta''+kt|}{m''(p'')^{1/4}}\lesssim(1+t)+\frac{(1+t)^2\langle\xi''\rangle}{\sqrt{p'}},
\end{align*}
where we have used \eqref{4.11-2}. If $\sqrt{p'}\ge1+t$, the last term is bounded by $(1+t)\langle\xi''\rangle$. If $\sqrt{p'}<1+t$, we have
\[
 \frac{1}{\sqrt{p'}}
 \le[\log(e+1+t)]^{1+a_0}\frac{1}{\langle t-\eta'/k\rangle[\log(e+\langle t-\eta'/k\rangle)]^{1+a_0}}.
\]
Combining the two cases and using \eqref{eq:single-transfer}, we obtain
\begin{align*}
 &\frac{\langle q\rangle}{m'(p')^{1/4}m''(p'')^{1/4}}
 \left(|(k,l'')|+\frac{|(k,l')|\,|\eta''+kt|}{\sqrt{p'}}\right)\\
 &\quad\lesssim(1+t)\langle\xi''\rangle
 +(1+t)^2[\log(e+1+t)]^{1+a_0}
   \langle\xi''\rangle^{2+a_0}\sqrt{\mu_7(t,\xi')\mu_7(t,q)}.
\end{align*}

\smallskip\noindent\textit{Case II: $\langle\xi'\rangle<\langle\xi''\rangle$.}
By \eqref{4.11-2} and $1+t\lesssim\langle\xi'\rangle\sqrt{p'}$, we have
\begin{align*}
 \frac{\langle q\rangle |(k,l'')|}{m'(p')^{1/4}m''(p'')^{1/4}}
 &\lesssim\sqrt{p'}+\sqrt{p''}
 \lesssim\langle\xi'\rangle(1+t)+\sqrt{p''},\\
 \frac{\langle q\rangle |(k,l')|\,|\eta''+kt|}{m'(p')^{1/4}m''(p'')^{1/4}\sqrt{p'}}
 &\lesssim(1+t)\left(1+\frac{\sqrt{p''}}{\sqrt{p'}}\right)
 \lesssim(1+t)+\langle\xi'\rangle\sqrt{p''}.
\end{align*}
Using Young's inequality in the two cases gives
\begin{align*}
 &e^{2\lambda\nu^{1/3}t}
 \big|\operatorname{Re}\langle M_7Z,M_7\mathcal F_0\rangle_{H^{N+1}}\big|\\
 &\quad\lesssim\|Z\|_{H^{N+1}}
 \big[\|\na_L\mathcal AX_{\neq}\|_{H^N}
              +(1+t)\|\mathcal AX_{\neq}\|_{H^N}\big]
       \|\langle\xi\rangle\widehat{\mathcal AX_{\neq}}\|_{L^1_\eta\ell^1_{k,l}}\\
 &\qquad+(1+t)^2[\log(e+1+t)]^{1+a_0}
 \big\|\sqrt{\mu_7}M_7Z\big\|_{H^{N+1}}
 \big\|\sqrt{\mu_7}\mathcal AX_{\neq}\big\|_{H^N}
       \|\langle\xi\rangle^{2+a_0}\widehat{\mathcal AX_{\neq}}\|_{L^1_\eta\ell^1_{k,l}}.
\end{align*}
For $d=1$ and $d=2+a_0$, Cauchy--Schwarz gives
\[
 \|\langle\xi\rangle^d\widehat f\|_{L^1_\eta\ell^1_{k,l}}
 \le\left(\sum_{k,l}\int_\R
             \langle\xi\rangle^{-2(N-d)}\,\mathrm d\eta\right)^{1/2}
             \|f\|_{H^N}\lesssim_{N,a_0}\|f\|_{H^N},
\]
where $N-d>3/2$. The $k=0$ projection also satisfies
\[
 \|(f*g)(0,\cdot,\cdot)\|_{L^2_\eta\ell^2_l}
 \le\|f*g\|_{\ell^2_kL^2_\eta\ell^2_l}\le\|f\|_{L^2_\eta\ell^2_{k,l}}\|g\|_{L^1_\eta\ell^1_{k,l}}.
\]
This proves \eqref{eq:top-force-point}.

By \eqref{eq:nonzero-norm-budget} and Cauchy--Schwarz in time,
\begin{align*}
 \|\na_L\mathcal AX_{\neq}\|_{L^2H^N}
 \|\mathcal AX_{\neq}\|_{L^2H^N}
 &\lesssim\nu^{-1/2}\nu^{-1/6}\mathcal E_{\neq,N}
 =\nu^{-2/3}\mathcal E_{\neq,N},\\
 \sup_{t\ge0}(1+t)e^{-2\lambda\nu^{1/3}t}
       \|\mathcal AX_{\neq}\|_{L^2H^N}^2
 &\lesssim\nu^{-1/3}\nu^{-1/3}\mathcal E_{\neq,N}
 =\nu^{-2/3}\mathcal E_{\neq,N}.
\end{align*}

For the last time weight put $x=\nu^{1/3}t$. Since $0<\nu\le1$,
\[
 1+t\le\nu^{-1/3}(1+x),\qquad
 \log(e+1+t)\lesssim(1+|\log\nu|)[1+\log(e+x)].
\]
It follows that
\begin{align*}
 &\sup_{t\ge0}(1+t)^2[\log(e+1+t)]^{1+a_0}e^{-2\lambda\nu^{1/3}t}\\
 &\quad\lesssim\nu^{-2/3}(1+|\log\nu|)^{1+a_0}
       \sup_{x\ge0}(1+x)^2[1+\log(e+x)]^{1+a_0}e^{-2\lambda x}
 \lesssim_{a_0,\lambda}\nu^{-2/3}(1+|\log\nu|)^{1+a_0}.
\end{align*}
Moreover,
\begin{align*}
 &\|\mathcal AX_{\neq}\|_{L^\infty(0,t_1;H^N)}
   \big\|\sqrt{\muEight}\mathcal AX_{\neq}\big\|_{L^2(0,t_1;H^N)}
   \big\|\sqrt{\mu_7}M_7Z\big\|_{L^2(0,t_1;H^{N+1})}\\
 &\qquad\lesssim\mathcal E_{\neq,N}(t_1)
   \big\|\sqrt{\mu_7}M_7Z\big\|_{L^2(0,t_1;H^{N+1})}.
\end{align*}
Combining these estimates gives \eqref{eq:top-force-integral} and finishes the proof.
\end{proof}

Combining the forcing bound in lemma~\ref{lem:top-forcing} with a direct $H^{N+1}$ energy estimate for $M_7(t,0,\p_y)Z$, we obtain the short-time zero-mode closure on $[0,\min(T,t_{*})]$.
\begin{lem}\label{lem:early-zero}
Under the bootstrap assumptions \eqref{eq:full-bootstrap}, for any $0<t_1\le\min(T,t_*)$, there holds
\begin{equation}\label{eq:early-zero-close}
\begin{aligned}
 &\|Z\|_{L^\infty(0,t_1;H^{N+1})}^2
       +\nu\|\na Z\|_{L^2(0,t_1;H^{N+1})}^2
       +\big\|\sqrt{\mu_7}M_7(t,0,D_y)Z\big\|_{L^2(0,t_1;H^{N+1})}^2\\
 &\qquad\le C_0\ep^2+
 C_{\Boot}\big[A_*\ep^2\nu^{-4/3}(1+|\log\nu|)+
                 \ep (1+|\log\nu|)^{1+a_0}\nu^{-2/3}\big]\ep^2.
\end{aligned}
\end{equation}
\end{lem}
\begin{proof}
Since the real, even multiplier $M_7(t,0,\eta)$ depends only on $(t,\eta)$, it commutes with the skew-adjoint linear operator in \eqref{eq:zero-system}. On the finite time interval $[0,t_{*}]$, we do not need transport cancellation for the internal nonlinearity $\mathcal{N}(Z)$. Instead, using the algebra property of $H^N$ and the uniform bounds $e^{-Ca_0^{-2}}\leq M_7\leq 1$ from \eqref{eq:M7-bounds}, we directly bound
\[
 \|\NL(Z)\|_{H^{N+1}}\lesssim
       \|Z\|_{H^{N+1}}\|\na Z\|_{H^{N+1}},\quad 
 |\langle M_7Z,M_7\NL(Z)\rangle_{H^{N+1}}|
 \lesssim\|Z\|_{H^{N+1}}^2\|\na Z\|_{H^{N+1}}.
\]
Absorbing $\|\na Z\|_{H^{N+1}}$ into the viscous dissipation via Young's inequality yields
\[
\begin{aligned}
 &\frac{1}{2}\frac{\mathrm d}{\mathrm dt}\|M_7Z\|_{H^{N+1}}^2
 +c\nu\|\na Z\|_{H^{N+1}}^2+\|\sqrt{\mu_7}M_7Z\|_{H^{N+1}}^2\\
 &\qquad\le C\nu^{-1}\|Z\|_{H^{N+1}}^4
       +\big|\operatorname{Re}\langle M_7Z,M_7\mathcal F_0\rangle_{H^{N+1}}\big|.
\end{aligned}
\]
Under the bootstrap ansatz \eqref{eq:full-bootstrap} on $[0,t_1]$ with $t_1\le\min(T,t_*)$, we have
\[
	\nu^{-1}\int_0^{t_1}\|Z\|_{H^{N+1}}^4\dd
 \lesssim_{\Boot}\nu^{-1}t_*\ep^4
 =A_*\ep^4\nu^{-4/3}(1+|\log\nu|),
\]
while \eqref{eq:top-force-integral} gives 
\[
	\int_0^{t_1}\big|\operatorname{Re}\langle M_7Z,M_7\mathcal F_0\rangle_{H^{N+1}}\big|\dd\lesssim_{\Boot}(1+|\log\nu|)^{1+a_0}\nu^{-2/3}\ep^3.
\]
Integrating the differential inequality over $[0,t_1]$ and using the equivalence $\|M_7 Z\|_{H^{N+1}}\asymp\|Z\|_{H^{N+1}}$ from \eqref{eq:early-zero-def} yields \eqref{eq:early-zero-close}.
\end{proof}

\subsubsection{Zero modes estimates after enhanced dissipation}
\label{sec:zero-late}
When $T> t_*$, we turn to the remaining interval $[t_{*},T]$. Since $2\lambda A_{*}\geq 2$ in \eqref{eq:switch-time}, the nonzero modes have decayed by a factor of $e^{-\lambda\nu^{1/3}t_{*}}\leq \nu$ by time $t_{*}$. This exponential decay suppresses the polynomial growth $(1+t)$ in the $H^{N+1}$ estimate of the Reynolds stress $\mathcal{S}_0$,  allowing us to close the modified energy estimate on $[t_{*},T]$ via lemma~\ref{lem:forced-high}.
\begin{lem}\label{lem:late-zero}
Let $T\ge t_*$ and assume the bootstrap ansatz \eqref{eq:full-bootstrap} holds. Then 
\begin{equation}\label{eq:tail-stress}
 \nu^{-1/2}\|\mathcal S_0\|_{L^2(t_*,T;H^{N+1})}
 \lesssim
 \nu^{-1}(1+\nu^{1/3}t_*)e^{-2\lambda\nu^{1/3}t_*}
 \mathcal E_{\neq,N}\lesssim C_{A_*}\nu(1+|\log\nu|)
 \mathcal E_{\neq,N}.
\end{equation}
If \eqref{1.23} holds with sufficiently small $\delta$, then
\begin{equation}\label{eq:late-zero-close}
 \|Z\|_{L^\infty(t_*,T;H^{N+1})}
       +\sqrt\nu\,\|\na Z\|_{L^2(t_*,T;H^{N+1})}
 \lesssim\ep+C_{A_*}\nu(1+|\log\nu|)\mathcal E_{\neq,N}
 \lesssim\ep.
\end{equation}
\end{lem}
\begin{proof}
For the same frequencies $q=\xi'+\xi''=\xi'_L+\xi''_L$,
\eqref{4.11-2} gives
\[
 \frac{\langle q\rangle}{m'(p')^{1/4}m''(p'')^{1/4}}
 \lesssim\frac{\sqrt{p'}+\sqrt{p''}}{m'(p')^{1/4}m''(p'')^{1/4}}\lesssim1+t.
\]
Using Young's inequality and placing the $H^N$ weight on the larger input, we obtain
\[
 \|\mathcal S_0(t)\|_{H^{N+1}}
 \lesssim(1+t)e^{-2\lambda\nu^{1/3}t}
                    \|\mathcal AX_{\neq}\|_{H^N}^2.
\]
With $x=\nu^{1/3}(t-t_*)$ and $0<\nu\le1$, we have
\begin{align*}
 \int_{t_*}^\infty(1+t)^2e^{-4\lambda\nu^{1/3}t}\,\mathrm dt&\le\nu^{-1}(1+\nu^{1/3}t_*)^2e^{-4\lambda\nu^{1/3}t_*}
       \int_0^\infty(1+x)^2e^{-4\lambda x}\,\mathrm dx\\
 &\lesssim_\lambda\nu^{-1}(1+\nu^{1/3}t_*)^2
       e^{-4\lambda\nu^{1/3}t_*}.
\end{align*}
This proves \eqref{eq:tail-stress}.
By \eqref{eq:normal-form-equivalence}, \eqref{eq:M7-bounds} and \eqref{eq:early-zero-close},
\[
 E_{0,N+1}^{w}(t_*)+\mathcal C_{N+1,\vartheta}(t_*)
 \asymp\|Z(t_*)\|_{H^{N+1}}^2
 \asymp\|M_7(t_*,0,D_y)Z(t_*)\|_{H^{N+1}}^2\lesssim\ep^2.
\]
Applying \eqref{eq:tame-result} with $s=N+1$ and $t_0=t_*$ gives
\begin{align*}
 \sup_{t_*\le t\le T}\|Z(t)\|_{H^{N+1}}
 +\sqrt\nu\|\nabla Z\|_{L^2(t_*,T;H^{N+1})}&\lesssim\|Z(t_*)\|_{H^{N+1}}
 +\nu^{-1/2}\|\mathcal S_0\|_{L^2(t_*,T;H^{N+1})}\\
 &\lesssim\ep+C_{A_*}\nu(1+|\log\nu|)\mathcal E_{\neq,N}
 \lesssim\ep,
\end{align*}
which proves \eqref{eq:late-zero-close}. Together, lemma~\ref{lem:early-zero} and lemma~\ref{lem:late-zero} establish the global $H^{N+1}$ zero modes bounds over the entire interval $[0,T]$.
\end{proof}

\section{Energy estimates on \texorpdfstring{$Q_{\neq},K_{\neq},H_{\neq}$}{Q, K, H}}
\label{sec:nonzero-nonlinear}

\qquad In this section we prove proposition~\ref{prop:linear}
by estimating the nonlinear terms in \eqref{eq:source-pairing}. We use the Fourier multiplier $M_7$ defined in \eqref{eq:M8} to handle the case when the frequencies in a product have different critical times. We first give the convolution estimates, and then estimate the nonlinear terms from \eqref{Non0NL} in turn. The two orders of interaction with
zero modes are treated separately at the end.

Recall the three good unknowns $Q,K,H$ defined in \eqref{yyy}, we  denote $X_{\neq}=(Q_{\neq},K_{\neq},H_{\neq})$ and $\mu_7(t,k,\eta,l):=-\dot{M_7}/M_7(t,k,\eta,l)$. We will frequently use the decay bound $(1+t)^q e^{-\lambda\nu^{1/3}t}\lesssim_q\nu^{-q/3},\ q\geq0$, and
\begin{equation}\label{eq:detailed-time-norms}
\begin{aligned}
\|\mathcal AX_{\neq}\|_{L^\infty_t H^N}
 +\nu^{1/6}\|\mathcal AX_{\neq}\|_{L^2_tH^N}
 +\nu^{1/2}\|\na_L\mathcal AX_{\neq}\|_{L^2_tH^N}
 +\|\sqrt{\mu_7}\mathcal AX_{\neq}\|_{L^2_tH^N}
 \lesssim\mathcal E_{\neq,N}^{1/2},
\end{aligned}
\end{equation}
and \eqref{4.11-1}--\eqref{4.11-2} give 
\begin{equation}\label{eq:a-bounds}
 1\le |(k,l)|\le mp^{1/4},\qquad \sqrt p\le(1+t)mp^{1/4},\qquad
 mp^{1/4}\le\left(1+\frac{t}{\tau_\nu}\right)\sqrt p.
\end{equation}
where the notation $p,\ m,\ \mathcal{A}$ are the same as in \eqref{2.8}, \eqref{eq:m-formula}, \eqref{888}, respectively.

We show some relationship between good unknowns $X_{\neq}$ and physical variables $(U_{\neq},\Theta_{\neq})$. From \eqref{3.20}, we have
\begin{equation}\label{eq:reconstruction}
|\widehat U_{\neq}|+|\widehat\Theta_{\neq}|
\lesssim p^{-1/4}|\widehat{X_{\neq}}|,
\qquad
\widehat U_{\neq}^3=-\sqrt{p_h}\,p^{-3/4}\widehat Q_{\neq}.
\end{equation}
\eqref{eq:reconstruction}
and \eqref{4.11-2} give
\begin{equation}\label{3x.1}
|\widehat{(U_{\neq},\Theta_{\neq})}|
\lesssim e^{-\lambda\nu^{1/3}t}|\widehat{\mathcal AX_{\neq}}|.
\end{equation}

For the interactions of nonzero modes, we write
\[
\xi=(k,\eta,l),\qquad \xi_L=(k,\eta-kt,l),\qquad 
\xi'=(k',\eta',l'),\qquad
\xi''=\xi-\xi'.
\]
For simplicity, we denote $p'=p(t,\xi'),\ m'=m(t,\xi')$, $p''=p(t,\xi'')$ and $m''=m(t,\xi'')$.

Every nonlinear term in \eqref{Non0NL} has a transport structure.
For a nonzero transported component $f_{\neq}$, its Fourier formula is
\begin{equation}\label{eq:transport-origin}
\widehat{(U_{\neq}\cdot\na_Lf_{\neq})}(\xi)
=i\int_{\xi'+\xi''-\xi=0}
 (\widehat U_{\neq}(\xi')\cdot\xi''_L)
 \widehat f_{\neq}(\xi'').
\end{equation}
Incompressible condition $\widehat U_{\neq}(\xi')\cdot\xi'_L=0$ also gives
\begin{equation}\label{eq:pp''}
	\widehat U_{\neq}(\xi')\cdot\xi''_L
=\widehat U_{\neq}(\xi')\cdot\xi_L,
\qquad |\widehat U_{\neq}(\xi')\cdot\xi''_L|\le\min(\sqrt{p''},\sqrt p)|\widehat U_{\neq}(\xi')|.
\end{equation}
$U_{\neq}^2$ has the bound
\begin{equation}\label{ll7}
|\widehat U_{\neq}^2(\xi')|\le (p')^{-1/2}|k',l'||\widehat U_{\neq}(\xi')|\lesssim
\langle t\rangle^{-1}\langle \xi'\rangle^2|\widehat U_{\neq}(\xi')|.
\end{equation}
The last bound follows from $1+t\lesssim\langle \xi'\rangle\sqrt{p'}$ and
$|(k',l')|\le\langle \xi'\rangle$.

Note that $(fg)_{\neq}=(f_{\neq}g_{\neq})_{\neq}+(f_0g_{\neq})_{\neq}+(f_{\neq}g_0)_{\neq}$, we set  $\mathscr P=\{(\neq,\neq),(0,\neq),(\neq,0)\}$. For each $(a,b)\in\mathscr P$, where the first index marks the advecting velocity and the second marks the transported quantity, we decompose the nonlinear terms in \eqref{Non0NL} as
\begin{equation}\label{eq:NL-sum}
\mathcal N_Q=\sum_{\mathscr P}(\mathcal Q^T_{ab}+\mathcal Q^P_{ab}),
\quad \mathcal N_K=\sum_{\mathscr P}\mathcal K_{ab},
\quad \mathcal N_H=\sum_{\mathscr P}\mathcal H_{ab},
\end{equation}
where the two nonlinear terms of $Q$ are given by
\[
\mathcal Q^T_{ab}=|\na_L|^{3/2}|\na_{L,h}|^{-1}(U_a\cdot\na_LU_b^3)_{\neq},\quad \mathcal Q^P_{ab}=\p_z|\na_{L,h}|^{-1}|\na_L|^{-1/2}(\p_i^LU_a^j\p_j^LU_b^i)_{\neq},
\]
and $\mathcal K_{ab},\mathcal H_{ab}$ denote the corresponding
ordered nonlinear terms in the $\mathcal N_K$ and $\mathcal N_H$ equations of \eqref{Non0NL}. 
For $\mathcal K_{ab}$, we expand the curl of the transport term via the exact identity
\begin{equation}\label{eq:K-product}
\begin{aligned}
\p_x(U_a\cdot\na_LU_b^2)-\p_y^L(U_a\cdot\na_LU_b^1)
={}\p_xU_a\cdot\na_LU_b^2-\p_y^LU_a\cdot\na_LU_b^1+U_a\cdot\na_LW_b^3,
\end{aligned}
\end{equation}
with $W_0^3=-\p_yU_0^1$.

We denote $I_1,I_2,I_3$ for nonzero modes self-interactions $(\neq,\neq)$ and
$J_1,J_2,J_3$ for mixed interactions $(0,\neq)$ and $(\neq,0)$, in the order
$\mathcal Q^T+\mathcal Q^P,\mathcal K,\mathcal H$. More precisely, the self-interaction contributions are
\[
\begin{aligned}
I_1&=\int_\xi\langle\xi\rangle^{2N}
 |\widehat{\A Q_{\neq}}+G\widehat{\A K_{\neq}}|
 |\widehat{\A\mathcal Q^T_{\neq\neq}}+\widehat{\A\mathcal Q^P_{\neq\neq}}|\,\mathrm d\xi,\\
I_2&=\int_\xi\langle\xi\rangle^{2N}
 |\widehat{\A K_{\neq}}+G\widehat{\A Q_{\neq}}|
 |\widehat{\A\mathcal K_{\neq\neq}}|\,\mathrm d\xi,\\
 I_3&=\int_\xi\langle\xi\rangle^{2N}
 |\widehat{\A H_{\neq}}|\,
 |\widehat{\A\mathcal H_{\neq\neq}}|\,\mathrm d\xi.
\end{aligned}
\]
For either mixed pair $(a,b)\in\{(\neq,0),(0,\neq)\}$, $J_1$ records the transport
and pressure contributions in the $Q$ equation. We estimate them
together and use the grouped notation
\[
\begin{aligned}
J_{1,ab}&:=\int_\xi\langle\xi\rangle^{2N}
 |\widehat{\A Q_{\neq}}+G\widehat{\A K_{\neq}}|
 |\widehat{\A(\mathcal Q^T_{ab}+\mathcal Q^P_{ab})}|\,\mathrm d\xi,\\
J_{2,ab}&:=\int_\xi\langle\xi\rangle^{2N}
 |\widehat{\A K_{\neq}}+G\widehat{\A Q_{\neq}}|
 |\widehat{\A\mathcal K_{ab}}|\,\mathrm d\xi,\\
J_{3,ab}&:=\int_\xi\langle\xi\rangle^{2N}
 |\widehat{\A H_{\neq}}|\,
 |\widehat{\A\mathcal H_{ab}}|\,\mathrm d\xi.
\end{aligned}
\]
By the triangle inequality,
\[
\mathcal I_{\neq,N}\le I_1+I_2+I_3
 +\sum_{(a,b)\in\{(0,\neq),(\neq,0)\}}
 \big(J_{1,ab}+J_{2,ab}+J_{3,ab}\big).
\]

\subsection{Nonzero modes self-interactions}
\subsubsection{The projected transport terms}\label{sec:ordinary}
\qquad As we will verify in section~\ref{sec:I1} below, after exploiting the Leray projection structure in $\mathcal Q_{\neq\neq}^T+\mathcal Q_{\neq\neq}^P$, and the curl cancellation in $\mathcal K_{\neq\neq}$, the Fourier multipliers acting on the nonlinear transport terms $(U_{\neq}\cdot\na_L \mathcal{Z}_{\neq})_{\neq}$ for $\mathcal{Z}\in\{U^1,U^2,U^3,\Theta\}$  all factor out the common weight $p(t,\xi)^{1/4}$ up to bounded coefficients. This motivates us to establish a unified estimate for the model transport pairing weighted by $\mathcal{A}(t,\xi)p(t,\xi)^{1/4}$.
\begin{lem}\label{lem:TYZ}
	Under the assumption of proposition~\ref{boots}, for any $\mathcal Z\in \{U^1,U^2,U^3,\Theta\}$ and $Y\in\{Q+GK,K+GQ,H\}$, we define
	\[
	\mathcal T_Y[\mathcal Z](t):=\int_\xi\langle\xi\rangle^{2N}|\widehat{\mathcal AY_{\neq}}(\xi)|\,\mathcal A(t,\xi)p(t,\xi)^{1/4}|\widehat{(U_{\neq}\cdot\na_L\mathcal{Z}_{\neq})_{\neq}}(\xi)|\,\mathrm d\xi.
	\]
	Then there holds 
\[
		\int_0^T \mathcal T_Y[\mathcal{Z}](t){\rm d}t\lesssim  (1+|\log\nu|)^{1+a_0}\nu^{-2/3}\mathcal E_{\neq,N}^{3/2}(T),
	\]
where the implicit constant depends only on $N,\alpha,\beta,a_0$,
and is independent of the bootstrap constant $\Boot$.
\end{lem}
\begin{proof}
	Using $|\widehat{\mathcal{Z}}_{\neq}(\xi'')|\lesssim (p'')^{-1/4}|\widehat X_{\neq}(\xi'')|$ from \eqref{eq:reconstruction}, we have
	\begin{equation}\label{eq:IYlow}
	\begin{aligned}
\int_0^T\mathcal T_Y[\mathcal{Z}]\dd
\lesssim\int_0^T\!\int_{\xi,\xi'}
 &\langle\xi\rangle^{2N}|\widehat{\mathcal AY_{\neq}}(\xi)|
 \frac{mp^{1/4}}{m''(p'')^{1/4}}
 |\widehat U_{\neq}(\xi')\cdot\xi_L''|
 |\widehat{\mathcal AX_{\neq}}(\xi'')|\,
 \mathrm d\xi'\mathrm d\xi\mathrm dt .
\end{aligned}
\end{equation}
We split $\mathcal T_{Y}[\mathcal{Z}]\le\mathcal T_Y^{\rm low}[\mathcal{Z}]+\mathcal T_{Y}^{\rm high}[\mathcal{Z}]$ by restricting the convolution to the regions $\langle\xi'\rangle\leq\langle\xi''\rangle$ and $\langle\xi'\rangle>\langle\xi''\rangle$, respectively. In both regimes, we estimate the rescaled kernel 
\begin{equation}\label{eq:ordinary-kernel}
\frac{m(t,\xi)p(t,\xi)^{1/4}}{m(t,\xi')p(t,\xi')^{1/4}m(t,\xi'')p(t,\xi'')^{1/4}}\,|\widehat U_{\neq}(\xi')\cdot\xi_L''|.
\end{equation}

By \eqref{eq:reconstruction} and the bounds on $M$, we have
\[
 m'(p')^{1/4}|\widehat U_{\neq}(\xi')|
 \lesssim m'|\widehat X_{\neq}(\xi')|
 \lesssim e^{-\lambda\nu^{1/3}t}
          |\widehat{\mathcal AX_{\neq}}(\xi')|.
\]
We use this bound to reconstruct the velocity after estimating the multiplier in \eqref{eq:ordinary-kernel}. The test factor satisfies $|\widehat{\mathcal AY}|\lesssim|\widehat{\mathcal AX_{\neq}}|$ since $|G|\lesssim1$.

\textit{Case 1 ($\langle \xi'\rangle\le\langle \xi''\rangle$):} The estimate \eqref{ll7} and \eqref{eq:a-bounds} give
\[
 |\widehat U(\xi')\cdot\xi_L''|
 \lesssim\big(|(k'',l'')|+\langle\xi'\rangle^2\sqrt{p''}/(1+t)\big)|\widehat U(\xi')|
 \lesssim\langle\xi'\rangle^2m''(p'')^{1/4}|\widehat U(\xi')|.
\]
Since $m'(p')^{1/4}\ge1$ and $mp^{1/4}\le(1+t/\tau_\nu)\sqrt p$, it follows that
\begin{equation}\label{eq:low-ordinary}
 \frac{mp^{1/4}}{m'(p')^{1/4}m''(p'')^{1/4}}|\widehat U(\xi')\cdot\xi_L''|
 \lesssim\left(1+\frac{t}{\tau_\nu}\right)\langle\xi'\rangle^2\sqrt p\,|\widehat U(\xi')|.
\end{equation}

For $q=0,1,2$, taking $x=\nu^{1/3}t$, we have
\begin{align}
 &\sup_{t\ge0}\left(1+\frac{t}{\tau_\nu}\right)(1+t)^q
 e^{-\lambda\nu^{1/3}t}
 \le\nu^{-q/3}\sup_{x\ge0}(1+x/1000)(1+x)^q e^{-\lambda x}
 \lesssim\nu^{-q/3},\label{eq:cutoff-time-weight}\\
 &\sup_{t\ge0}\left(1+\frac{t}{\tau_\nu}\right)(1+t)^2
 [\log(e+1+t)]^{1+a_0}e^{-\lambda\nu^{1/3}t}
 \lesssim\nu^{-2/3}(1+|\log\nu|)^{1+a_0}.\label{eq:cutoff-log-time-weight}
\end{align}

Multiplying \eqref{eq:low-ordinary} by $m'(p')^{1/4}$ and reconstructing the velocity as above, we deduce from \eqref{eq:cutoff-time-weight} and Young's convolution inequality that
\begin{equation}\label{eq:ordinary-low-norm}
 \int_0^T\mathcal T_Y^{\rm low}[\mathcal Z]\dd
 \lesssim\|\na_L\mathcal AY\|_{L^2H^N}
 \|\langle\xi\rangle^2\widehat{\mathcal AX_{\neq}}\|_{L^\infty L^1}
 \|\mathcal AX_{\neq}\|_{L^2H^N}\lesssim\nu^{-2/3}\mathcal E_{\neq,N}^{3/2}.
\end{equation}

\textit{Case 2 ($\langle \xi''\rangle<\langle \xi'\rangle$): }
We write $\mathcal T_Y^{\rm high,x}[\mathcal Z],\ \mathcal T_Y^{\rm high,y}[\mathcal Z],\ \mathcal T_Y^{\rm high,z}[\mathcal Z]$ for the part with the $x,y,z$ derivative on $\mathcal Z$, respectively. Using $|k'',l''|\le m''(p'')^{1/4}$ and $m'(p')^{1/4}\ge 1$ from \eqref{eq:a-bounds}, the $x,z$ derivatives satisfy
\begin{align*}
 \frac{mp^{1/4}}{m'(p')^{1/4}m''(p'')^{1/4}}|k''\widehat U_{\neq}^1(\xi')+l''\widehat U_{\neq}^3(\xi')|
 &\le\frac{mp^{1/4}}{m'(p')^{1/4}m''(p'')^{1/4}}|(k'',l'')|\,|\widehat U_{\neq}(\xi')|\\
 &\le\left(1+\frac{t}{\tau_\nu}\right)\sqrt p\,|\widehat U_{\neq}(\xi')|.
\end{align*}

The corresponding $x,z$ contributions are estimated as in Case 1. For the $y$ derivative, using $|(k',l')|\le m'(p')^{1/4}$, $mp^{1/4}\le(1+t/\tau_\nu)\sqrt p$,
$\sqrt p\le\sqrt{p'}+\sqrt{p''}$, and
$\frac{|\eta''-k''t|}{m''(p'')^{1/4}}\le1+t$, we have
\begin{align}\label{eq:high-y-kernel}
 \frac{mp^{1/4}}{m'(p')^{1/4}m''(p'')^{1/4}}|\widehat U_{\neq}^2(\xi')|\,|\eta''-k''t|
 &\le\frac{mp^{1/4}}{m'(p')^{1/4}m''(p'')^{1/4}}\frac{|(k',l')|}{\sqrt{p'}}|\eta''-k''t|
       |\widehat U_{\neq}(\xi')|\nonumber\\
 &\le\left(1+\frac{t}{\tau_\nu}\right)\frac{\sqrt p}{\sqrt{p'}}\frac{|\eta''-k''t|}{m''(p'')^{1/4}}
       |\widehat U_{\neq}(\xi')|\nonumber\\
 &\le\left(1+\frac{t}{\tau_\nu}\right)(1+t)\left(1+\frac{\sqrt{p''}}{\sqrt{p'}}\right)
       |\widehat U_{\neq}(\xi')|\nonumber\\
 &\le\left(1+\frac{t}{\tau_\nu}\right)\left((1+t)+\frac{(1+t)^2\langle\xi''\rangle}{\sqrt{p'}}\right)
       |\widehat U_{\neq}(\xi')|.
\end{align}

If $\sqrt{p'}\ge1+t$, the second term is bounded by $(1+t)\langle\xi''\rangle$. If $\sqrt{p'}<1+t$, using $p'\ge1+(t-\eta'/k')^2$, we have
\[
 \frac{1}{\sqrt{p'}}\le[\log(e+1+t)]^{1+a_0}\frac{1}{\langle t-\eta'/k'\rangle[\log(e+\langle t-\eta'/k'\rangle)]^{1+a_0}}.
\]

By \eqref{eq:single-transfer}, we further obtain
\[
 \frac{(1+t)^2\langle\xi''\rangle}{\sqrt{p'}}
 \lesssim_{a_0}(1+t)^2[\log(e+1+t)]^{1+a_0}
 \langle\xi''\rangle^{2+a_0}\sqrt{\mu_7(\xi')\mu_7(\xi)}.
\]

Combining the two regimes gives
\begin{equation}\label{eq:high-y-final}
\begin{aligned}
 &\frac{mp^{1/4}}{m'(p')^{1/4}m''(p'')^{1/4}}|\widehat U_{\neq}^2(\xi')|\,|\eta''-k''t|\\
 &\lesssim\left(1+\frac{t}{\tau_\nu}\right)\Big[(1+t)\langle\xi''\rangle+(1+t)^2[\log(e+1+t)]^{1+a_0}\langle\xi''\rangle^{2+a_0}
 \sqrt{\mu_7(\xi')\mu_7(\xi)}\Big]|\widehat U_{\neq}(\xi')|.
\end{aligned}
\end{equation}

For the lower-frequency factor, Cauchy--Schwarz gives

\begin{equation}\label{eq:low-weight-L1}
 \|\langle\xi\rangle^{2+a_0}\widehat f\|_{L^1_\eta\ell^1_{k,l}}
 \le\left(\sum_{k,l}\int_\R
       \langle\xi\rangle^{-2(N-2-a_0)}\,\mathrm d\eta\right)^{1/2}
       \|f\|_{H^N}\lesssim_{N,a_0}\|f\|_{H^N}.
\end{equation}
 Here $N-2-a_0>3/2$, which holds for every integer $N\ge4$ because $0<a_0<1/2$. Multiplying \eqref{eq:high-y-final} by $m'(p')^{1/4}$ and reconstructing the velocity, we place the lower-frequency factor in Fourier $L^1$. Then \eqref{eq:cutoff-time-weight}, \eqref{eq:cutoff-log-time-weight} and \eqref{eq:low-weight-L1} give
\begin{align}\label{eq:ordinary-high-norm}
 \int_0^T\mathcal T_Y^{{\rm high},y}[\mathcal Z]\dd
 &\lesssim \|\mathcal AX_{\neq}\|_{L^\infty_tH^N}\big[\nu^{-1/3}\|\mathcal AY\|_{L^2H^N}\|\mathcal AX_{\neq}\|_{L^2H^N}\nonumber\\
 &\qquad+(1+|\log\nu|)^{1+a_0}\nu^{-2/3}\|\sqrt{\mu_7}\mathcal AY\|_{L^2H^N}\|\sqrt{\mu_7}\mathcal AX_{\neq}\|_{L^2H^N}\big]\nonumber\\
 &\lesssim (1+|\log\nu|)^{1+a_0}\nu^{-2/3}\mathcal E_{\neq,N}^{3/2},
\end{align}
where we have used \eqref{eq:detailed-time-norms} and
$|\widehat{\mathcal AY}|
\lesssim|\widehat{\mathcal AX_{\neq}}|$. Together with the estimates of $\mathcal T_Y^{\rm high,x}[\mathcal Z]$ and $\mathcal T_Y^{\rm high,z}[\mathcal Z]$, this finishes the proof of lemma~\ref{lem:TYZ}.
\end{proof}

\subsubsection{The estimate of  \texorpdfstring{$I_1$}{I1}, \texorpdfstring{$I_2$}{I2} and \texorpdfstring{$I_3$}{I3}}\label{sec:I1}
\qquad Viewed in isolation, the transport term $\mathcal Q_{\neq\neq}^T=|\na_L|^{3/2}|\na_{L,h}|^{-1}(U_{\neq}\cdot\na_LU_{\neq}^3)_{\neq}$ carries the symbol $p^{3/4}p_h^{-1/2}$ and the pressure term $\mathcal Q_{\neq\neq}^P=\p_z|\na_{L,h}|^{-1}|\na_L|^{-1/2}(\p_i^LU_{\neq}^j\p_j^LU_{\neq}^i)_{\neq}$ involves the multiplier $lp^{1/4}p_h^{-1/2}$. When $|l|>p_h^{1/2}$, they both contain an additional factor $|l|^{3/2}p_h^{-1/2}$. Fortunately, we notice that the combination of $\mathcal Q_{\neq\neq}^P+\mathcal Q_{\neq\neq}^T$  forms a Leray projection and shows the cancellation due to the incompressibility condition, and reduces the corresponding multiplier to $p^{1/4}$. Specifically, let $F=(F^1,F^2,F^3)^T=(U_{\neq}\cdot\na_L U_{\neq})_{\neq}^T$. Expanding the divergence and using the incompressibility condition $\na_L\cdot U_{\neq}=0$, we have 
\[
\na_L\cdot F
=\partial_i^LU_{\neq}^j\partial_j^LU_{\neq}^i
 +U_{\neq}^j\partial_j^L(\partial_i^LU_{\neq}^i)
=\partial_i^LU_{\neq}^j\partial_j^LU_{\neq}^i.
\]
Using $\widehat{\partial_z\Delta_L^{-1}}=-il/p$, we can express the Fourier sum of the two $Q$-nonlinear terms as
\begin{equation}\label{eq:projected-Q}
\begin{aligned}
\widehat{\mathcal Q^T_{\neq\neq}+\mathcal Q^P_{\neq\neq}}
=\frac{p^{3/4}}{\sqrt{p_h}}
 \left(\widehat F^3-\frac{l}{p}\xi_L\cdot\widehat F\right)=p^{1/4}
 \frac{-kl\widehat F^1-(\eta-kt)l\widehat F^2+p_h\widehat F^3}{\sqrt p\sqrt{p_h}}
 =p^{1/4}e_2\cdot\widehat F,
\end{aligned}
\end{equation}
where 
\[
	e_2:=\frac{(-kl,-(\eta-kt)l,p_h)}{\sqrt p\sqrt{p_h}}
\]
is a unit vector perpendicular to $\xi_L$. Multiplying \eqref{eq:projected-Q} by $\mathcal{A}(t,\xi)$, and taking $Y=Q+G K$ and $ \mathcal{Z}\in\{U^1,U^2,U^3\}$ in  lemma~\ref{lem:TYZ} and using $|e_2(\xi)|=1$, we have
\begin{equation}\label{eq:I2-detailed}
\begin{aligned}
\int_0^T I_1\dd
&\le\int_0^T\big(\mathcal{T}_{Q+GK}[U^1](t)+\mathcal{T}_{Q+GK}[U^2](t)+\mathcal{T}_{Q+GK}[U^3](t)\big)\dd\\
&\lesssim (1+|\log\nu|)^{1+a_0}\nu^{-2/3}\mathcal E_{\neq,N}^{3/2}.
\end{aligned}
\end{equation}
Thus, we finish the estimate of $I_1$.

A similar cancellation structure holds for the curl term $\mathcal K_{\neq\neq}=ic_{\beta}|\na_L|^{1/2}|\na_{L,h}|^{-1}(\p_x(U_{\neq}\cdot\na_L U_{\neq}^2)_{\neq}-\p_y^L(U_{\neq}\cdot\na_L U_{\neq}^1)_{\neq})$. Indeed, \eqref{eq:K-product} gives 
\begin{equation}\label{eq:K-three-symbols}
\begin{aligned}
\widehat{\mathcal K_{\neq\neq}}(\xi)
={}&-ic_{\beta}\frac{p^{1/4}}{\sqrt{p_h}}\int_{\xi'+\xi''=\xi}
 \big(\widehat U(\xi')\cdot\xi''_L\big)\\
&\quad\times\left(
 k'\widehat U^2(\xi'') -(\eta'-k't)\widehat U^1(\xi'')+
 k''\widehat U^2(\xi'')-(\eta''-k''t)\widehat U^1(\xi'')
 \right)\,\mathrm d\xi'.
\end{aligned}
\end{equation}
Using $k'+k''=k$ and $(\eta'-k't)+(\eta''-k''t)=\eta-kt$, we obtain the exact cancellation 
\begin{align}
k'\widehat U^2(\xi'')-(\eta'-k't)\widehat U^1(\xi'')
 +k''\widehat U^2(\xi'')-(\eta''-k''t)\widehat U^1(\xi'')
&=k\widehat U^2(\xi'')-(\eta-kt)\widehat U^1(\xi'')\nonumber\\
&=\sqrt{p_h}\,e_1(\xi)\cdot\widehat U(\xi''),\label{eq:K-pair-cancel}
\end{align}
where
\[
e_1:=\frac{(-\eta+kt,k,0)}{\sqrt{p_h}}
\]
is a horizontal unit vector perpendicular to $\xi_L$. Consequently, the factor $p_h^{-1/2}$ in \eqref{eq:K-three-symbols} is canceled by $\sqrt{p_h}$ in \eqref{eq:K-pair-cancel}, and $|e_1(\xi)|=1$ implies
\begin{equation}\label{eq:K-integer-ratio}
\frac{|k\widehat U^2(\xi'')-(\eta-kt)\widehat U^1(\xi'')|}{\sqrt{p_h}}
\le|\widehat U(\xi'')|.
\end{equation}
Using \eqref{eq:K-pair-cancel}, \eqref{eq:K-integer-ratio} and $\widehat{\mathcal{K}_{\neq\neq}}(\xi)=-c_{\beta}p(t,\xi)^{1/4}e_1(\xi)\cdot\widehat{F}(\xi)$, we have
\begin{equation}\label{eq:I22}
	|\widehat{\mathcal{K}_{\neq\neq}}(\xi)|\leq c_{\beta}p(t,\xi)^{1/4}(|(\widehat{U_{\neq}\cdot\na_L U_{\neq}^1})_{\neq}(\xi)|+|(\widehat{U_{\neq}\cdot\na_LU_{\neq}^2})_{\neq}(\xi)|)
\end{equation}
Multiplying \eqref{eq:I22} by $\mathcal{A}(t,\xi)$, and taking  $Y=K +GQ$ and $\mathcal Z\in\{U^1,U^2\}$ in lemma~\ref{lem:TYZ}, we have
\begin{equation}\label{eq:I3-detailed}
\int_0^T I_2\dd
\le c_{\beta}\int_0^T\big(\mathcal{T}_{K+GQ}[U^1](t)+\mathcal{T}_{K+GQ}[U^2](t)\big)\dd
\lesssim c_{\beta}(1+|\log\nu|)^{1+a_0}\nu^{-2/3}\mathcal E_{\neq,N}^{3/2}.
\end{equation}

Similarly, since $\mathcal{H}_{\neq\neq}=|\na_L|^{1/2}(U_{\neq}\cdot\na_L\Theta_{\neq})_{\neq}$, taking $Y=H$ and  $\mathcal{Z}=\Theta$ in lemma~\ref{lem:TYZ} directly gives
\begin{equation}\label{eq:I4}
	\int_0^T I_3\dd=\int_0^T \mathcal{T}_{H}[\Theta](t){\rm d}t\lesssim (1+|\log\nu|)^{1+a_0}\nu^{-2/3}\mathcal E_{\neq,N}^{3/2}.
\end{equation}
Combining  \eqref{eq:I2-detailed}, \eqref{eq:I3-detailed}, and \eqref{eq:I4}, we conclude that the self-interactions satisfy
\begin{equation}\label{eq:self-bound}
\int_0^T\big(I_1+I_2+I_3\big)\dd
\lesssim b_{\beta}(1+|\log\nu|)^{1+a_0}\nu^{-2/3}\mathcal E_{\neq,N}^{3/2}(T).
\end{equation}
\subsection{Mixed interactions}\label{sec:regularity}
\qquad For each ordered pair $(a,b)=(0,\neq)$ or $(\neq,0)$, let
$F_{ab}^j=(U_a\cdot\na_LU_b^j)_{\neq}$ and
$F_{ab}^\theta=(U_a\cdot\na_L\Theta_b)_{\neq}$. Since $\na_L\cdot U_{\neq}=0$ and $\p_x U_0=0$, we apply the Leray projection identity \eqref{eq:projected-Q} and the curl cancellation \eqref{eq:K-pair-cancel}   to obtain
\begin{equation}\label{eq:mixed-rows}
\begin{aligned}
\widehat{\mathcal Q^T_{ab}+\mathcal Q^P_{ab}}
 &=p^{1/4}\left[
 -\frac{kl}{\sqrt p\sqrt{p_h}}\widehat F_{ab}^1
 -\frac{(\eta-kt)l}{\sqrt p\sqrt{p_h}}\widehat F_{ab}^2
 +\frac{\sqrt{p_h}}{\sqrt p}\widehat F_{ab}^3\right],\\
\widehat{\mathcal K_{ab}}
 &=-c_\beta p^{1/4}\left[
 -\frac{\eta-kt}{\sqrt{p_h}}\widehat F_{ab}^1
 +\frac{k}{\sqrt{p_h}}\widehat F_{ab}^2\right],\\
\widehat{\mathcal H_{ab}}&=p^{1/4}\widehat F_{ab}^\theta,
\end{aligned}
\end{equation}
\qquad Throughout this subsection, we set $\xi=(k,\eta,l),\ \xi_L=(k,\eta-kt,l),\ \zeta=(k,\eta',l'),\ \zeta_L=(k,\eta'-kt,l')$ and $(\eta'',l'')=(\eta-\eta',l-l')$ for $k\ne0$. As before, we write $p=p(t,\xi),\ m=m(t,\xi),\ p'=p(t,\zeta)$ and $m'=m(t,\zeta)$, and split each $J_{i,ab}\leq J_{i,ab}^{\text{low}}+J_{i,ab}^{\text{high}}$, $i=1,2,3$, by restricting the convolution to $\langle\eta'',l''\rangle\le8\langle\zeta\rangle$
and $\langle\eta'',l''\rangle>8\langle\zeta\rangle$, respectively.
\subsubsection{Zero modes transporting nonzero modes}
\qquad Let $F_{0\neq}^j=(U_0\cdot\nabla_L U_{\neq}^j)_{\neq}$ and $F_{0\neq}^{\theta}=(U_{0}\cdot\na_L\Theta_{\neq})_{\neq}$. Similar to the rescaled kernel \eqref{eq:ordinary-kernel}, we focus on
\begin{equation}\label{eq:forward-kernel}
\frac{mp^{1/4}}{m(t,\zeta)p^{1/4}(t,\zeta)}
\,|\widehat U_0(\eta-\eta',l-l')\cdot\zeta_L|.
\end{equation}

\textit{Case I ($\langle\eta'',l''\rangle\le 8\langle\zeta\rangle$): }By \eqref{4.11-3}, we have 
\[
\frac{mp^{1/4}}{m'(p')^{1/4}}
\lesssim \langle\eta'',l''\rangle.
\]
Moreover, the incompressibility condition $\widehat{U}_0(\eta'',l'')\cdot(\eta'',l'')=0$ implies
\[\widehat U_0(\eta-\eta',l-l')\cdot\zeta_L
=\widehat U_0(\eta-\eta',l-l')\cdot\xi_L.\]
Combining these bounds with \eqref{eq:mixed-rows} and $|\widehat{\mathcal{A}Y}(\xi)|\lesssim |\widehat{\mathcal{A}X_{\neq}}(\xi)|$ yields
\begin{align}\label{eq:forward-low-norm}
\int_0^T \big(J_{1,0\neq}^{\rm low}+J_{2,0\neq}^{\rm low}+J_{3,0\neq}^{\rm low}\big)\dd
&\lesssim \|\na_L\mathcal AX_{\neq}\|_{L^2_tH^N}\|\langle\eta,l\rangle\widehat U_0\|_{L^\infty_t L^1}\|\mathcal AX_{\neq}\|_{L^2_tH^N}\nonumber\\
&\lesssim b_{\beta}\nu^{-1/2}\nu^{-1/6}
 \sup (E_{0,N}^{w})^{1/2}\mathcal E_{\neq,N}.
\end{align}

\textit{Case II ($\langle\eta'',l''\rangle>8\langle\zeta\rangle$): }By \eqref{4.11-3} and the incompressibility condition, we have
\begin{equation}\label{eq:forward-high}
 \frac{mp^{1/4}}{m'(p')^{1/4}}|\widehat U_0\cdot\zeta_L|
 =\frac{mp^{1/4}}{m'(p')^{1/4}}|\widehat U_0\cdot\xi_L|
 \lesssim\langle\eta'',l''\rangle\sqrt p\,|\widehat U_0|.
\end{equation}

Assigning $\langle\xi\rangle^N\lesssim\langle\eta'',l''\rangle^N$ to the higher-frequency zero modes and using Young's inequality, we obtain
\begin{align}\label{eq:forward-high-norm}
 \int_0^T \big(J_{1,0\neq}^{\rm high}+J_{2,0\neq}^{\rm high}+J_{3,0\neq}^{\rm high}\big)\dd
 &\lesssim\|\na_L\mathcal AX_{\neq}\|_{L^2H^N}
 \|U_0\|_{L^\infty H^{N+1}}\|\widehat{\mathcal AX_{\neq}}\|_{L^2L^1}\nonumber\\
 &\lesssim b_{\beta}\nu^{-2/3}\sup \|Z\|_{H^{N+1}}\mathcal E_{\neq,N}.
\end{align}
Summing \eqref{eq:forward-low-norm} and \eqref{eq:forward-high-norm} over the three pairings yields
\begin{equation}\label{eq:forward-bound}
\int_0^T\big(J_{1,0\neq}+J_{2,0\neq}+J_{3,0\neq}\big)\dd
\lesssim b_{\beta}\nu^{-2/3}
\sup\|Z\|_{H^{N+1}}
\mathcal E_{\neq,N}(T).
\end{equation}

\subsubsection{Nonzero modes transporting zero modes}
\qquad For $U_{\neq}\cdot\na \mathcal{Z}_0$ with $\mathcal{Z}_0\in\{U_0^1,U_0^2,U_0^3,\Theta_0\}$, the incompressibility condition $\widehat{U}_{\neq}(\zeta)\cdot\zeta_L=0$ gives
\[
(\eta-\eta')\widehat U_{\neq}^2(\zeta)
 +(l-l')\widehat U_{\neq}^3(\zeta)
=\widehat U_{\neq}(\zeta)\cdot\xi_L.
\]
Combining this with \eqref{4.11-3} and \eqref{eq:reconstruction}, we obtain
\begin{equation}\label{eq:reverse-kernel}
 p^{1/4}|\widehat{\mathcal A(U_{\neq}\cdot\na \mathcal{Z}_0)}(\xi)|
 \lesssim\sqrt p\sum_{l'}\int_\R\langle\eta-\eta',l-l'\rangle
 |\widehat{\mathcal AX_{\neq}}(\zeta)|\,|\widehat{\mathcal{Z}}_0(\eta-\eta',l-l')|\,\mathrm d\eta'.
\end{equation}
Applying \eqref{eq:reverse-kernel} to \eqref{eq:mixed-rows} and  $\langle\xi\rangle^N\lesssim
\langle\eta-\eta',l-l'\rangle^N+\langle\zeta\rangle^N$ give
\begin{align*}
J_{1,\neq0}+J_{2,\neq0}+J_{3,\neq0}
&\lesssim b_{\beta}\|\na_L\mathcal AX_{\neq}\|_{H^N}
 \|(U_0,\Theta_0)\|_{H^{N+1}}
 \|\mathcal AX_{\neq}\|_{H^N}.
\end{align*}
Integrating in time and using \eqref{eq:detailed-time-norms}, we conclude that
\begin{equation}\label{eq:reverse-bound}
\begin{aligned}
\int_0^T \big(J_{1,\neq0}+J_{2,\neq0}+J_{3,\neq0}\big)\dd
&\lesssim\|\na_L\mathcal AX_{\neq}\|_{L^2_tH^N}
 \|(U_0,\Theta_0)\|_{L^\infty_tH^{N+1}}
 \|\mathcal AX_{\neq}\|_{L^2_tH^N}\\
 &\lesssim b_{\beta}\nu^{-2/3}\sup\|Z\|_{H^{N+1}}\mathcal E_{\neq,N}(T).
\end{aligned}
\end{equation}
Combining \eqref{eq:forward-bound} and \eqref{eq:reverse-bound} yields
\begin{equation}\label{eq:mixed-total}
\int_0^T\sum_{(a,b)\in\{(0,\neq),(\neq,0)\}}
\big(J_{1,ab}+J_{2,ab}+J_{3,ab}\big)\dd\lesssim b_{\beta}\nu^{-2/3}\sup\|Z\|_{H^{N+1}}\mathcal E_{\neq,N}(T),
\end{equation}
By the bootstrap assumption \eqref{eq:full-bootstrap}, the coefficient in \eqref{eq:mixed-total} is at most $C_{\Boot}\ep\nu^{-2/3}$, which is small for $\ep\le\delta\nu^{2/3}/(1+|\log\nu|)^{1+a_0}$.
Combining \eqref{eq:self-bound} and \eqref{eq:mixed-total} with the nonlinear terms decomposition above proves proposition~\ref{prop:linear}.

\allowdisplaybreaks

\vskip .2in 
\subsection*{Declaration of competing interest} The authors declare that they have no conflict of interest.

\vspace{4mm} \noindent\textbf{Acknowledgements.} {The authors would like to thank F. Shao, D. Wei and Z. Zhang} for illuminating discussions. Huang was partially supported by China Postdoctoral Science Foundation (No. 2026M793341). Sun was partially supported by the National Natural Science Foundation of China (No. 12501296) and the Fundamental Research Funds for the Central Universities. Xu was partially supported by the National Natural Science Foundation of China (No. 12571244; 12171040) and the National Key R\&D Program of China (No. 2020YFA0712900).  

\appendix

\section{Dynamics of nonzero modes when \texorpdfstring{$B_{\beta}>0$}{B beta > 0}} 
\label{secA}
\qquad In fact, when $B_{\beta}>0$,    we can also prove the linear stability of the system \eqref{1.6}. Let us talk about the differences in the proof process here. We introduce a new group of good unknowns
\[
	Q^*_{\neq}:=-|\na_L|^2|\na_{L,h}|^{-1}U^3_{\neq},\ K^*_{\neq}:=-ic_{\beta}|\na_L||\na_{L,h}|^{-1}W_{\neq}^3,\ H_{\neq}^*:=-|\na_L|\Theta_{\neq},
\]
where $c_{\beta}=\sqrt{\beta/(\beta-1)}$, and then we rewrite \eqref{Non0L} as 
\begin{align}\label{Non0L1}
	\begin{cases}
		\p_t \widehat{Q_{\neq}^*} + \nu p \widehat{Q_{\neq}^*} - \frac{1}{2} \frac{\p_t p_{h}}{p_{h}}\widehat{Q_{\neq}^*} - \sgn(\beta)\sqrt{B_{\beta}}\frac{l}{p^{1/2}}\widehat{K_{\neq}^*} + \alpha \frac{p_h^{1/2}}{p^{1/2}}\widehat{H_{\neq}^*} + \frac{2}{c_{\beta}}\frac{k^2}{p_h}\frac{l}{p^{1/2}}\widehat{K_{\neq}^*}  = 0, \\
		\p_t\widehat{K_{\neq}^*} + \nu p\widehat{K_{\neq}^*} +\frac{1}{2}\frac{l^{2}}{p}\frac{\p_t p_h}{p_h}\widehat{K_{\neq}^*}+ \sgn(\beta)\sqrt{B_{\beta}}\frac{l}{p^{1/2}}\widehat{Q_{\neq}^*} = 0,\\
		\p_t\widehat{H_{\neq}^*} + \nu p\widehat{H_{\neq}^*}  -\frac{1}{2}\frac{\partial_{t}p}{p} \widehat{H_{\neq}^*}- \alpha \frac{p_h^{1/2}}{p^{1/2}}\widehat{Q_{\neq}^*} = 0.
	\end{cases}
\end{align}
We use a suitably constructed Fourier multiplier $m^*$ with symbol $m^*=m^*(t,k,\eta,l)$, defined by 
\begin{align}\label{A.....4}	
	-\frac{\dot{m^*}}{m^*}=\left\{\begin{array}{ll}
		\frac{| k \left( \eta- kt \right)|}{k^2+\left(\eta- kt\right)^2}, & \text{ if }  |t-\frac{\eta}{k}|\leq \tau_\nu, \\
		0,  & \text{ if } |t-\frac{\eta}{k}|> \tau_\nu
	\end{array}\right. ,\qquad \tau_\nu=1000\nu^{-1/3},
\end{align}
with the initial data $m^*(0,k,\eta,l)=\left(1+\min\{\max(\eta/k,0)^2,\tau_\nu^2\}\right)^{1/2}$ and $m^*(t,0,\eta,l)=1$. We also define the  multiplier $M^{*}=M_{1}M_{5}$ with $M_{1}$ and $M_{5}$ as in  \eqref{4..5}.

\begin{re}
	$(1)$ Similar to the lemma~\ref{lem3.1}, one can directly obtain the formula of $m^*(t, k, \eta,l)$, here we omit it. In particular, we have $\nu^{\frac{1}{3}}\lesssim m^*\left(t, k, \eta,l\right) \lesssim \nu^{-\frac{1}{3}}$, which is worse than \eqref{4.11-1}. 

	$(2)$  $m^*(t, k, \eta,l)$  and the frequency have the following relationship
\begin{equation}\label{A.11-2}
	|k|/\sqrt{p_h}\lesssim m^*\left(t, k, \eta,l \right) \lesssim\sqrt{p_h}/|k|.
\end{equation}
\end{re} 

The proof of \eqref{eq:linear-decay}--\eqref{A.33x} as in  theorem~\ref{thm:linear}  can be directly obtained from the following lemma, which gives the linear enhanced dissipation and inviscid damping on $U_{\neq}$.  
\begin{lem}\label{lemA.4}
	Assume that $0<\nu\le1$, $\alpha>0$ and $\beta \in \mathbb{R}$ with $B_\beta>0$, then there holds
	\begin{align}
		&\left\| U^1_{\neq} (t) \right\|_{L^{2}} \lesssim_{\alpha,\beta}  \, b_{\beta}^* \, e^{-\lambda^*  \nu^{\frac{1}{3}} t} \left(\left\|  u_{\mathrm{in}}\right\|_{H^{2}}+\left\|  \theta_{\mathrm{in}}\right\|_{H^{2}} \right) , \label{A.1x}  \\
		&{\left\langle {t} \right\rangle} \left\|  U^2_{\neq} (t) \right\|_{L^{2}} \lesssim_{\alpha,\beta}  \, b_{\beta}^* \, e^{-\lambda^* \nu^{\frac{1}{3}} t} \left(\left\|  u_{\mathrm{in}}\right\|_{H^{3}}+\left\|  \theta_{\mathrm{in}}\right\|_{H^{3}} \right) , \label{A.2x}  \\
		&\left\| (U^3_{\neq}, \Theta_{\neq}) (t) \right\|_{L^2} \lesssim_{\alpha,\beta}   e^{-\lambda^* \nu^{\frac{1}{3}} t} \left(\left\|  u^3_{\mathrm{in}}\right\|_{H^{2}} + c_{\beta} \left\|  \p_xu^2_{\mathrm{in}}-\p_yu^1_{\mathrm{in}}\right\|_{H^{1}}+\left\|  \theta_{\mathrm{in}}\right\|_{H^{2}} \right),\label{A.3x}  
	\end{align}
	where  $\lambda^*=1/16$ and $b_{\beta}^*=\max \left\{c_{\beta}, c_{\beta}^{-1} \right\}$.
\end{lem}

\begin{proof}
To control the evolution of the nonzero modes, we introduce the Fourier multiplier
\begin{equation}\label{A88}
	\mathcal{A}^*:=m^*M^{*}e^{\lambda^* \nu^{\frac{1}{3}}t},
\end{equation} 
which incorporates the multiplier $m^*$ already used in the linear analysis,   additional ghost weights ${M^*}$ as well as the time weight $e^{\lambda^*\nu^{\frac{1}{3}}t}$, which $\lambda^*$ can be confirmed later. Recall the system  \eqref{Non0L1}, energy estimates yield that 
\begin{equation}
\begin{aligned}\label{A99}
    \frac{1}{2}&\frac{d}{dt}|\mathcal{A}^*(\widehat{Q_{\neq}^*},\widehat{K_{\neq}^*},\widehat{H_{\neq}^*})|^2 +\left[ \nu p+\left(-\frac{\dot{M^*}}{M^*}\right)\right]|\mathcal{A}^*(\widehat{Q_{\neq}^*},\widehat{K_{\neq}^*},\widehat{H_{\neq}^*})|^2  \\
    &=\lambda^*\nu^{\frac{1}{3}}|\mathcal{A}^*(\widehat{Q_{\neq}^*},\widehat{K_{\neq}^*},\widehat{H_{\neq}^*})|^2-\left[(-\frac{\dot{m^*}}{m^*})-\frac{1}{2} \frac{\p_t p_{h}}{p_{h}}\right]|\widehat{\mathcal{A}^*Q_{\neq}^*}|^2-\left[(-\frac{\dot{m^*}}{m^*})+\frac{1}{2}\frac{l^{2}}{p}\frac{\p_t p_h}{p_h}\right]|\widehat{\mathcal{A}^*K_{\neq}^*}|^2\\
	&\quad-\left[(-\frac{\dot{m^*}}{m^*})-\frac{1}{2} \frac{\p_t p}{p}\right]|\widehat{\mathcal{A}^*H_{\neq}^*}|^2 -\frac{2}{c_{\beta}}\frac{k^2}{p_h}\frac{l}{\sqrt{p}}\Rea\big(\widehat{\mathcal{A}^*K_{\neq}^*}\overline{\widehat{\mathcal{A}^*Q_{\neq}^*}}\big).
\end{aligned}
\end{equation}

Denote $e_{\neq}^*(t):=|\mathcal{A}^*(\widehat{Q_{\neq}^*},\widehat{K_{\neq}^*},\widehat{H_{\neq}^*})|^2$. Notice that the definition of $m^*$ as in \eqref{A.....4}, we have
 $ 
 -\frac{\dot{m^*}}{m^*}+\frac{\nu p}{2}-\frac{1}{2} \big|\frac{\p_t p}{p_{h}}\big|\gtrsim 0$.
The three diagonal stretching coefficients have absolute value
at most $|k(\eta-kt)|/p_h$. The last term is bounded by
$(-\dot M_5/M_5)e_{\neq}^*/2$.
\eqref{eq:enhanced-rate} and $\lambda^*=1/16$ gives
\begin{equation}\label{A010}
\frac{\mathrm d}{\mathrm dt}e_{\neq}^*
+c\left(\nu p-\dot M^*/M^*\right)e_{\neq}^*\le0.
\end{equation}
Integration in frequency and time yields, for $r\ge0$,
\begin{equation}\label{Ax.0}
\begin{aligned}
&\|\mathcal A^*(Q_{\neq}^*,K_{\neq}^*,H_{\neq}^*)(t)\|_{H^r}^2
 +\nu\|\mathcal A^*\nabla_L(Q_{\neq}^*,K_{\neq}^*,H_{\neq}^*)\|_{L_t^2H^r}^2\\
&\quad+\big\|\mathcal A^*\sqrt{-\dot{ M^*}/M^*}
       (Q_{\neq}^*,K_{\neq}^*,H_{\neq}^*)\big\|_{L_t^2H^r}^2
\lesssim\|\mathcal A^*(0)(Q_{\mathrm{in}}^*,K_{\mathrm{in}}^*, H_{\mathrm{in}}^*)\|_{H^r}^2.
\end{aligned}
\end{equation}
Note that $m^*(0)\le |k|^{-1}|k,\eta|$, we have $m^*(0)|\widehat{Q_{\mathrm{in}}^*}|\le |k|^{-1}|k,\eta,l|^2|\widehat{u_{\mathrm{in}}^3}|$, $m^*(0)|\widehat{K_{\mathrm{in}}^*}|
 \le c_\beta |k|^{-1}|k,\eta,l|
 |\widehat{\partial_xu_{\mathrm{in}}^2-\partial_yu_{\mathrm{in}}^1}|$, $m^*(0)|\widehat{H_{\mathrm{in}}^*}|
 \le |k|^{-1}|k,\eta,l|^2|\widehat{\theta_{\mathrm{in}}}|$. Consequently,
\[
 \|\mathcal A^*(0)
 (Q_{\mathrm{in}}^*,K_{\mathrm{in}}^*,H_{\mathrm{in}}^*)\|_{H^r}
 \lesssim
 \|u_{\mathrm{in}}^3\|_{H^{r+2}}
 +c_\beta\|\p_xu_{\rm in}^2-\p_yu_{\rm in}^1\|_{H^{r+1}}
 +\|\theta_{\mathrm{in}}\|_{H^{r+2}}.
\]
Next, recalling the definition of $Q_{\neq}^{*}$, $K_{\neq}^{*}$ and  $H_{\neq}^{*}$ and combining with \eqref{A.11-2}, we have 
\begin{align}\label{Ax.1}
	|\widehat{U^3_{\neq}}|+|\widehat{\Theta_{\neq}}|= p_h^{\frac{1}{2}}p^{-1}|\widehat{Q_{\neq}^*}|+p^{-\frac{1}{2}}|\widehat{H_{\neq}^*}| \lesssim |k|^{-1}\left(|m^*M^*\widehat{Q_{\neq}^*}|+|m^*M^*\widehat{H_{\neq}^*}|\right).
\end{align}
Then we use \eqref{Ax.0} and \eqref{Ax.1} to obtain
\begin{equation}
\begin{aligned}\label{Ax.3}
	\left\|({U^3_{\neq}},\Theta_{\neq})\right\|_{L^2}^2
	\lesssim e^{-2\lambda^* \nu^\frac{1}{3} t} \left(\| u^3_{\mathrm{in}}\|_{H^2}^2+c_{\beta}^2 \left\|  \p_xu^2_{\mathrm{in}}-\p_yu^1_{\mathrm{in}}\right\|_{H^1}^2+\| \theta_{\mathrm{in}}\|_{H^2}^2 \right),
\end{aligned}
\end{equation}
Using  \eqref{3.20} and \eqref{A.11-2}, we get $\widehat{U_{\neq}^1}$ and $\widehat{U_{\neq}^2}$ which satisfy
\begin{equation}\label{Al6}
\begin{aligned}
	|\widehat{U_{\neq}^1}|
	\lesssim \frac{1}{|k|} \left(c_{\beta}^{-1}|m^*M^*\widehat{K_{\neq}^*}| + |m^*M^*\widehat{Q_{\neq}^*}|\right),\ |\widehat{U_{\neq}^2}|
	\lesssim   \frac{1}{p^{\frac{1}2}}\left(c_{\beta}^{-1} |m^*M^*\widehat{K_{\neq}^*}| +\frac{|l|}{|k|} |m^*M^*\widehat{Q_{\neq}^*}|\right).
\end{aligned}
\end{equation}
Together with \eqref{Ax.0} and $\langle t\rangle p^{-1/2}(1+|l|/|k|)\lesssim \langle k,\eta,l\rangle$, this gives
\begin{equation}\label{Ax.5}
\begin{aligned}
	\left\|{U^1_{\neq}}\right\|_{L^2}^2&\lesssim 
	\max \left\{ c_{\beta}^2,  c_{\beta}^{-2} \right\} e^{-2\lambda^* \nu^{\frac{1}{3}} t} \left(\| u_{\mathrm{in}}\|_{H^2}^2+\| \theta_{\mathrm{in}}\|_{H^2}^2 \right),\\
	\left\|{U^2_{\neq}}\right\|_{L^2}^2&\lesssim 
	{\left\langle {t} \right\rangle}^{-2} \,\max \left\{ c_{\beta}^2,  c_{\beta}^{-2} \right\} e^{-2\lambda^* \nu^{\frac{1}{3}} t} \left(\| u_{\mathrm{in}}\|_{H^{3}}^2+\| \theta_{\mathrm{in}}\|_{H^{3}}^2 \right),
\end{aligned}
\end{equation}
where we have used  the fact $p^{-1} \lesssim {\left\langle {t} \right\rangle}^{-2} |k, \eta, l|^2$. Hence, together \eqref{Ax.3} with \eqref{Ax.5} gives \eqref{A.1x}-\eqref{A.3x} and  the proof of lemma~\ref{lemA.4} is finished.  
\end{proof}

\section{Construction of the zero-mode energy weight}\label{sec:zero-weight}
\allowdisplaybreaks[0]
Recall \eqref{eq:phase} and \eqref{eq:omega}. For $\sigma\in\{-1,1\}^3$, $\xi_j\ne0$ and $\sum_j\xi_j=0$, set $\omega_j=\omega_0(\xi_j)$ and
\[
 \Psi(\sigma,\xi_1,\xi_2,\xi_3)=\sum_{j=1}^3\sigma_j\omega_j,\qquad
 g(\xi)=\p_\eta\omega_0(\xi)
 =\frac{(\alpha^2-B_\beta)\eta l^2}{(\alpha^2\eta^2+B_\beta l^2)^{1/2}(\eta^2+l^2)^{3/2}}.
\]
Define the normalized set
\[
 \mathfrak S:=\left\{(\sigma,\zeta_1,\zeta_2,\zeta_3):
 \begin{gathered}
 \sigma\in\{-1,1\}^3,\quad \max_{i<j}|\sigma_i-\sigma_j|>0,\quad
 \zeta_j\in\mathbb R^2\setminus\{0\},\\
 \sum_{j=1}^3\zeta_j=0,\qquad\max_{1\le j\le3}|\zeta_j|=1
 \end{gathered}\right\},
\]
and its space--time resonant subset
\[
 \mathfrak R:=\left\{(\sigma,\zeta_1,\zeta_2,\zeta_3)\in\mathfrak S:
 \Psi=0,\quad\sigma_1g(\zeta_1)=\sigma_2g(\zeta_2)=\sigma_3g(\zeta_3)\right\}.
\]
Here both components of $\zeta_j$ are allowed to be real, since normalization by the largest frequency does not preserve integer wavenumbers.

\begin{lem}\label{lem:common-weight}
Assume $\alpha^2\ne\B$. For $N\le s\le N+1$, there is a positive smooth function $F$ on $[0,1]$, depending only on $\alpha,B_\beta$, such that, with $L_s(0)=1$,
\begin{equation}\label{eq:weight-properties}
\begin{aligned}
 L_s(\xi)&=(1+|\xi|^2F(\eta^2/|\xi|^2))^s\asymp\langle\xi\rangle^{2s},\\
 |\na L_s(\xi)|&\lesssim\langle\xi\rangle^{2s-1},\\
 |\partial_q L_s(\sqrt q,l)|&\lesssim\langle\sqrt q,l\rangle^{2s-2},\qquad q\ge0,\quad q+l^2>0.
\end{aligned}
\end{equation}
If $\alpha\ne2\sqrt\B$ and $\sqrt\B\ne2\alpha$, there is a relatively open neighborhood $\mathcal U_1$ of $\mathfrak R$ in $\mathfrak S$, chosen only in terms of $\alpha,B_\beta$. Let $\mathcal U$ be its extension under simultaneous positive scaling of the three frequencies. On $\mathcal U$, we have
\begin{equation}\label{eq:common-weight}
 L_s(\xi_j)=(1+|g(\xi_j)|^{-2})^s,\qquad j=1,2,3.
\end{equation}
Moreover, the constant
\begin{equation}\label{eq:phase-separation-constant}
 c_{\mathrm{sep}}(\alpha,B_\beta):=\frac{1}{2}\min\left\{1,
 \inf_{(\sigma,\zeta_1,\zeta_2,\zeta_3)\in\mathfrak S\setminus\mathcal U_1}
 \left(|\Psi|+\sum_{i<j}|\sigma_i g(\zeta_i)-\sigma_jg(\zeta_j)|\right)\right\}
\end{equation}
is positive, where $\inf\varnothing=+\infty$. Set
\[
 \rho_{\mathrm{sep}}(\alpha,B_\beta)
 :=\min\left\{\tfrac14\min(\alpha,\sqrt{B_\beta}),c_{\mathrm{sep}}(\alpha,B_\beta)\right\}.
\]
For $\sigma\in\{-1,1\}^3$ of nonuniform sign, $\sum_j\xi_j=0$ and $\xi_j\ne0$, one has
\begin{equation}\label{eq:phase-gap-outside-weight}
 \begin{aligned}
 |\Psi|\le\rho_{\mathrm{sep}}(\alpha,B_\beta),\quad
 (\sigma,\xi_1,\xi_2,\xi_3)\notin\mathcal U
 \quad\Longrightarrow\quad
 R\sum_{i<j}|\sigma_i g(\xi_i)-\sigma_j g(\xi_j)|\ge c_{\mathrm{sep}}(\alpha,B_\beta),
 \end{aligned}
\end{equation}
where $R=\max_j|\xi_j|$. If $\mathfrak R=\varnothing$, we take $\mathcal U_1=\mathcal U=\varnothing$ and $L_s(\xi)=\langle\xi\rangle^{2s}$. At either critical ratio we also take this standard weight, without asserting \eqref{eq:phase-gap-outside-weight}. The implicit constants in \eqref{eq:weight-properties} depend only on $\alpha,B_\beta,N$.
\end{lem}
\begin{proof}
We first assume $\alpha\ne2\sqrt\B$ and $\sqrt\B\ne2\alpha$. By the homogeneity of $\omega_0$, we have
\[
 \omega_0(\xi/R)=\omega_0(\xi),\qquad g(\xi/R)=Rg(\xi).
\]
We begin with the compactness argument needed for both the weight and the constant in \eqref{eq:phase-separation-constant}. Let
$(\sigma^{(n)},\zeta_1^{(n)},\zeta_2^{(n)},\zeta_3^{(n)})\in\mathfrak S$ be a sequence satisfying
\begin{equation}\label{eq:approximate-joint-resonance}
 \Psi_n\to0,\qquad
 \sum_{i<j}|\sigma_i^{(n)}g(\zeta_i^{(n)})-\sigma_j^{(n)}g(\zeta_j^{(n)})|\to0,
 \qquad \Psi_n:=\sum_j\sigma_j^{(n)}\omega_0(\zeta_j^{(n)}).
\end{equation}
Since there are finitely many sign triples, we may take a subsequence for which $\sigma^{(n)}=\sigma$ is fixed. Since $\max_j|\zeta_j^{(n)}|=1$, we take a further subsequence and relabel the indices so that
\[
 \zeta_j^{(n)}\to\zeta_j^*,\qquad |\zeta_1^*|=1,\qquad
 g_j^{(n)}:=\sigma_jg(\zeta_j^{(n)})\to g^*\in\mathbb R,
 \quad j=1,2,3.
\]
Here and below we keep the same notation for subsequences. The convergence of $g_j^{(n)}$ follows from the smoothness of $\omega_0$ near $\zeta_1^*\ne0$ and $g_i^{(n)}-g_j^{(n)}\to0$. We consider the following two cases.

\textit{Case 1 ($\zeta_j^*\ne0$ for all $j=1,2,3$):}
Passing to the limit gives
\[
 \sum_j\zeta_j^*=0,\qquad\max_j|\zeta_j^*|=1,\qquad
 \sum_j\sigma_j\omega_0(\zeta_j^*)=0,\qquad
 \sigma_jg(\zeta_j^*)=g^*,\quad j=1,2,3.
\]
Thus $(\sigma,\zeta_1^*,\zeta_2^*,\zeta_3^*)\in\mathfrak R$.

\textit{Case 2 ($\min_j|\zeta_j^*|=0$):}
Since $\sum_j\zeta_j^*=0$ and $|\zeta_1^*|=1$, only one limiting vector can vanish. Without loss of generality, assume $\zeta_3^*=0$, so that $\zeta_2^*=-\zeta_1^*$. If $\sigma_1=-\sigma_2$, the evenness of $\omega_0$ gives
\[
 \sigma_1\omega_0(\zeta_1^{(n)})+\sigma_2\omega_0(\zeta_2^{(n)})\to0,
 \qquad \liminf_{n\to\infty}|\Psi_n|\ge\min\{\alpha,\sqrt\B\}>0,
\]
contrary to \eqref{eq:approximate-joint-resonance}. Hence $\sigma_1=\sigma_2=-\sigma_3$. Since $g$ is odd,
\[
 g^*=\sigma_1g(\zeta_1^*)=\sigma_2g(-\zeta_1^*)=-g^*,
 \qquad g_j^{(n)}\to0,\quad j=1,2,3,
\]
and the phase relation implies
\[
 \omega_0(\zeta_3^{(n)})\to2\omega_0(\zeta_1^*).
\]
Expressing $\zeta_3^{(n)}=r_n(u_n,v_n)$, where $r_n\to0$ and $u_n^2+v_n^2=1$, we take a further subsequence so that $(u_n,v_n)\to(u_*,v_*)\in\mathbb S^1$. The formula for $g$ gives
\[
 \frac{(\alpha^2-B_\beta)u_nv_n^2}{\sqrt{\alpha^2u_n^2+B_\beta v_n^2}}
 =r_ng(\zeta_3^{(n)})=\sigma_3r_ng_3^{(n)}\to0.
\]
Thus $u_*v_*^2=0$. Also $g(\zeta_1^*)=0$, so both limiting directions lie on coordinate axes. Consequently,
\[
 \omega_0(\zeta_1^*),\ \lim_{n\to\infty}\omega_0(\zeta_3^{(n)})
 \in\{\alpha,\sqrt\B\},\qquad
 \lim_{n\to\infty}\omega_0(\zeta_3^{(n)})=2\omega_0(\zeta_1^*).
\]
This forces $\alpha=2\sqrt\B$ or $\sqrt\B=2\alpha$, contrary to the assumptions. Thus Case 2 is excluded.

We have proved that every sequence satisfying \eqref{eq:approximate-joint-resonance} has a subsequence converging to a point of $\mathfrak R$, with all three limiting frequencies nonzero. Applying this to sequences in $\mathfrak R$ shows that $\mathfrak R$ is compact. On $\mathfrak R$, the common signed derivative cannot vanish. Otherwise, the formula for $g$ implies that each vector lies on a coordinate axis. Since the three vectors are nonzero and sum to zero, they must lie on the same axis. Their frequencies are then equal, giving
\[
 0=\left|\sum_j\sigma_j\omega_0(\zeta_j)\right|
 =\left|\sum_j\sigma_j\right|\omega_0(\zeta_1)
 \ge\min\{\alpha,\sqrt\B\}>0,
\]
a contradiction. If $\mathfrak R\ne\varnothing$, the minima
\begin{equation}\label{eq:weight-positive-minimum}
 \begin{aligned}
 r_*(\alpha,B_\beta)&:=\min_{(\sigma,\zeta_1,\zeta_2,\zeta_3)\in\mathfrak R}
             \min_{1\le j\le3}|\zeta_j|>0,\\
 d_*(\alpha,B_\beta)&:=\min_{(\sigma,\zeta_1,\zeta_2,\zeta_3)\in\mathfrak R}
             \min_{1\le j\le3}|\zeta_j|\,|g(\zeta_j)|>0
 \end{aligned}
\end{equation}
are therefore well defined. Suppressing their parameter dependence, set
\[
 \mathcal U_1:=\left\{(\sigma,\zeta_1,\zeta_2,\zeta_3)\in\mathfrak S:
 \min_j|\zeta_j|>r_*/2,\quad
 \min_j|\zeta_j|\,|g(\zeta_j)|>d_*/2\right\}.
\]
This is a relatively open neighborhood of $\mathfrak R$. Its conic extension $\mathcal U$ satisfies $|\xi_j|\asymp R$ and $|\xi_j|\,|g(\xi_j)|>d_*/2$ for $j=1,2,3$. If $\mathfrak R=\varnothing$, set $\mathcal U_1=\mathcal U=\varnothing$.

We now prove that $c_{\mathrm{sep}}(\alpha,B_\beta)>0$. Suppose that the infimum in \eqref{eq:phase-separation-constant} is zero. There is a sequence in $\mathfrak S\setminus\mathcal U_1$ for which
\[
 |\Psi_n|+\sum_{i<j}|\sigma_i^{(n)}g(\zeta_i^{(n)})-\sigma_j^{(n)}g(\zeta_j^{(n)})|\le n^{-1}.
\]
By the preceding argument, a subsequence converges to a point of $\mathfrak R$. If $\mathfrak R=\varnothing$, this is impossible. Otherwise, since $\mathcal U_1$ is a relatively open neighborhood of $\mathfrak R$, this subsequence lies in $\mathcal U_1$ for all sufficiently large $n$, again a contradiction. This proves the positivity of \eqref{eq:phase-separation-constant}. For $\zeta_j=\xi_j/R$ outside $\mathcal U_1$, homogeneity gives
\begin{align*}
 R\sum_{i<j}|\sigma_i g(\xi_i)-\sigma_jg(\xi_j)|
 &=\sum_{i<j}|\sigma_i g(\zeta_i)-\sigma_jg(\zeta_j)|\\
 &\ge2c_{\mathrm{sep}}(\alpha,B_\beta)-|\Psi|
 \ge c_{\mathrm{sep}}(\alpha,B_\beta),
 \qquad |\Psi|\le\rho_{\mathrm{sep}}(\alpha,B_\beta),
\end{align*}
which proves \eqref{eq:phase-gap-outside-weight}.

It remains to construct the weight. When $\mathfrak R\ne\varnothing$, choose a smooth positive function $\gamma$ on $[0,\infty)$, constant on $[0,d_*/4]$ and satisfying $\gamma(v)=v$ for $v\ge d_*/2$. Set
\[
 L_s(\xi)=\left(1+\frac{|\xi|^2}{\gamma(|\xi|\,|g(\xi)|)^2}\right)^s,
 \qquad L_s(0)=1.
\]
For a triple in $\mathcal U$, $|\xi_j|\,|g(\xi_j)|>d_*/2$, so this formula gives \eqref{eq:common-weight}. With $u=\eta^2/|\xi|^2$, we have
\begin{align*}
 (|\xi|\,|g(\xi)|)^2
 &=\frac{(\alpha^2-\B)^2u(1-u)^2}{\B+(\alpha^2-\B)u},\\
 F(u)&=\gamma\left(\sqrt{\frac{(\alpha^2-\B)^2u(1-u)^2}{\B+(\alpha^2-\B)u}}\right)^{-2},\\
 L_s(\xi)&=(1+|\xi|^2F(u))^s\asymp\langle\xi\rangle^{2s},\\
 |\na L_s(\xi)|
 &\lesssim\langle\xi\rangle^{2s-2}
       \big(|\xi|\,|F(u)|+|\xi|^2|F'(u)|\,|\na u|\big)
 \lesssim\langle\xi\rangle^{2s-1}.
\end{align*}
Since $\gamma$ is positive and constant near zero, $F$ is smooth and positive on $[0,1]$. Moreover, with $q=\eta^2$ and $u=q/(q+l^2)$,
\[
 \partial_q L_s(\sqrt q,l)
 =s\big(1+(q+l^2)F(u)\big)^{s-1}\big[F(u)+(1-u)F'(u)\big],
\]
which proves the last bound in \eqref{eq:weight-properties}. The value $L_s(0)=1$ gives a $C^1$ extension at the origin. All choices depend only on $\alpha,B_\beta$. If $\mathfrak R=\varnothing$, or if either critical relation holds, we take $F\equiv1$. This proves \eqref{eq:weight-properties} and finishes the proof.
\end{proof}

\begin{re}
If $\sqrt\B>\alpha$, the nonzero-frequency three-wave space--time resonant set is empty. To see this, write the phase relation as $\omega_3=\omega_1+\omega_2$. A nonzero common signed derivative gives
\[
 \eta_1\eta_2>0,\qquad
 |\eta_3|=|\eta_1|+|\eta_2|,\qquad
 \left|\frac{l_3}{\eta_3}\right|>
 \max\left\{\left|\frac{l_1}{\eta_1}\right|,
            \left|\frac{l_2}{\eta_2}\right|\right\},
\]
where the last inequality follows from the strict increase of $\omega_0$ in $|l/\eta|$. Hence $|l_3|>|l_1|+|l_2|$, contrary to $\sum_j l_j=0$. The zero common signed derivative is excluded by the equal-frequency argument in the preceding proof. The low-frequency limit at $\sqrt\B=2\alpha$ is estimated in Case IV of lemma~\ref{lem:phase-blocks}.
\end{re}


\addcontentsline{toc}{section}{References}
\begin{thebibliography}{99}
\setlength{\itemsep}{0pt}
\setlength{\parsep}{0pt}
	\bibitem{Arbon2025}
	R. Arbon, Enhanced dissipation, Taylor dispersion, and inviscid damping of Couette flow in the Boussinesq system on the plane. 	arXiv:2501.00690 (2025).

	\bibitem{AB2025}	
	R. Arbon, J. Bedrossian, Quantitative hydrodynamic stability for Couette flow on unbounded domains with Navier boundary conditions. Comm. Math. Phys. 406(6): Paper No. 129, 57 pp (2025).

	\bibitem{BMN1999}
	A. Babin, A. Mahalov, B. Nicolaenko, On the regularity of three-dimensional rotating Euler-Boussinesq equations. Math. Models Methods Appl. Sci. 9(7): 1089–1121 (1999).

	\bibitem{BBCD2023}
	J. Bedrossian, R. Bianchini, M. Coti Zelati, M. Dolce, Nonlinear inviscid damping and shear-buoyancy instability in the two-dimensional Boussinesq equations. Comm. Pure Appl. Math. 76(12): 3685–3768 (2023).

	\bibitem{MR3612004}
	J. Bedrossian, P. Germain, N. Masmoudi, On the stability threshold for the 3D Couette flow in
	Sobolev regularity. Ann. of Math.  185(2): 541--608 (2017).

	\bibitem{MR3974608}
	J. Bedrossian, P. Germain, N. Masmoudi, Stability of the Couette flow at high Reynolds numbers in two dimensions and three dimensions. Bull. Amer. Math. Soc.  56(3):  373--414 (2019).

	\bibitem{MR4126259}
	J. Bedrossian, P. Germain, N. Masmoudi, Dynamics near the subcritical transition of the 3D Couette flow I: Below threshold case. Mem. Amer. Math. Soc. 266(1294), v+158pp (2020).

	\bibitem{BGM2022}
	J. Bedrossian, P. Germain, N. Masmoudi, Dynamics near the subcritical transition of the 3D Couette flow II: Above threshold case. Mem. Amer. Math. Soc. 279(1377) v+135pp (2022).

	\bibitem{BHIW1-2024}
	J. Bedrossian, S. He, S. Iyer, F. Wang, Uniform inviscid damping and inviscid limit of the 2D Navier-Stokes equation with Navier boundary conditions. arXiv:2405.19249  (2024).

	\bibitem{BHIW2-2024}
	J. Bedrossian, S. He, S. Iyer, F. Wang, Stability threshold of nearly-Couette shear flows with Navier boundary conditions in 2D.  Comm. Math. Phys. 406(2): Paper No. 28, 42 pp (2025).

	\bibitem{MR3448924}
	J. Bedrossian, N. Masmoudi, V. Vicol, Enhanced dissipation and inviscid damping in the inviscid limit of the Navier-Stokes equations near the two dimensional Couette flow. Arch. Ration. Mech. Anal. 219(3): 1087--1159  (2016).

	\bibitem{MR3867637}
	J. Bedrossian, V. Vicol, F. Wang, The Sobolev stability threshold for 2D shear flows near Couette. J. Nonlinear Sci. 28(6): 2051--2075 (2018).

	\bibitem{Bradshaw1969}
	P. Bradshaw, The analogy between streamline curvature and buoyancy in turbulent shear flow. J. Fluid Mech. 36(1): 177–191 (1969).

	
	
	\bibitem{MR4121130}
	Q. Chen, T. Li, D. Wei, Z. Zhang, Transition threshold for the 2-D Couette flow in a finite channel. Arch. Ration. Mech. Anal. 238(1): 125--183 (2020).

	\bibitem{CWZ2024}
	Q. Chen, D. Wei, Z. Zhang, Transition threshold for the 3D Couette flow in a finite channel. Mem. Amer. Math. Soc. 296(1478), v+178pp (2024).

	\bibitem{Zelati2023}
	M. Coti Zelati, A. Del Zotto, Suppression of lift-up effect in the 3D Boussinesq equations around a stably stratified Couette flow. Quarterly of Applied Mathematics, 83(2): 389–401 (2024).

	\bibitem{ZZW2024}
	M. Coti Zelati, A. Del Zotto, K. Widmayer, Stability of viscous three-dimensional stratified Couette flow via dispersion and mixing. Comm. Pure Appl. Math. 0:1–56,   (2026)

	
	\bibitem{CZDZW2025}
	M. Coti Zelati, A. Del Zotto, K. Widmayer, On the stability of viscous three-dimensional rotating Couette flow. arXiv:2501.17735 (2025).	

	\bibitem{CWW2025}
	S. Cui, L. Wang, W. Wang, Stability threshold of Couette flow for 3D Boussinesq system in Sobolev spaces. 	arXiv:2504.16401 (2025).

	\bibitem{DR}
	P. G. Drazin and W. H. Reid, Hydrodynamic stability, Cambridge University Press, Cambridge, 1981.

	\bibitem{DWZ2021}
	W. Deng, J. Wu, P. Zhang, Stability of Couette flow for 2D Boussinesq system with vertical dissipation. J. Funct. Anal. 281(12): Paper No. 109255, 40 pp (2021).

	\bibitem{HZ2023}
	C. Hao, W. Zhang, A priori estimates for free boundary problem of 3D incompressible inviscid rotating Boussinesq equations. Z. Angew. Math. Phys. 74(2): Paper No. 80, 21 pp (2023).

\bibitem{HSX2024}
	W. Huang, Y. Sun, X. Xu, On the Sobolev stability threshold for 3D Navier-Stokes equations with rotation near the Couette flow. Calc. Var. Partial Differential Equations 65 (9): Paper No. 251, 52 pp (2026).

	\bibitem{HSX2024-2}
	W. Huang, Y. Sun, X. Xu, Stability of the Couette flow for 3D Navier-Stokes equations with rotation. arXiv:2412.11005 (2024).

	\bibitem{Huang2018}
	Y. Huang, On the classical Bradshaw–Richardson number: Its generalized form, properties, and application in turbulence. Phys. Fluids. 30(12): 125110 (2018).

	\bibitem{IMT2017}
	T. Iwabuchi, A. Mahalov, R. Takada, Global solutions for the incompressible rotating stably stratified fluids. Math. Nachr. 290(4): 613–631 (2017).

	\bibitem{K1887}
	L. Kelvin, Stability of fluid motion-rectilinear motion of viscous fluid between two parallel plates, Phil.
	Mag., 24: 188–196  (1887).
	\bibitem{Kno2505}
	N. Knobel, Suppression of fluid echoes and Sobolev stability threshold for 2D dissipative fluid equations around Couette flow, arXiv:2505.23391.
	\bibitem{KMY2012}
	H. Koba, A. Mahalov, T. Yoneda, Global well-posedness of the rotating Navier-Stokes-Boussinesq equations with stratification effects, Adv. Math. Sci. Appl. 22: 61–90 (2012).

	\bibitem{LLZ2025}
	H. Li, N. Liu, W. Zhao, Stability threshold of the two-dimensional Couette flow in the whole plane. 	arXiv:2501.10818 (2025).

	\bibitem{LMZ2025}
	H. Li, N. Masmoudi, W. Zhao, Asymptotic stability of two-dimensional Couette flow in a viscous fluid. Arch. Ration. Mech. Anal. 249(5): Paper No. 56 (2025).

	\bibitem{LSWZ2025}
	M. Li, C. Sun, C. Wang, D. Wei, Z. Zhang, Transition threshold for the Navier–Stokes–Coriolis system at high Reynolds numbers. Proc. Lond. Math. Soc. 133 (3): Paper No. e70194 (2026).

	\bibitem{Majda2003}
	A. Majda, Introduction to PDEs and waves for the atmosphere and ocean. Courant Lecture Notes in Mathematics, 9. New York University, Courant Institute of Mathematical Sciences, New York; American Mathematical Society, Providence, RI, 2003. x+234 pp.

	\bibitem{MSZ2022}
	N. Masmoudi, B. Said-Houari, W. Zhao, Stability of the Couette flow for a 2D Boussinesq system without thermal diffusivity. Arch. Ration. Mech. Anal. 245(2): 645–752 (2022). 

	\bibitem{MR4176913}
	N. Masmoudi, W. Zhao, Enhanced dissipation for the 2D Couette flow in critical space. Comm. Partial Differential Equations. 45(12): 1682--1701(2020).

	\bibitem{MS2022}
	P. Mu, S. Schochet, Dispersive estimates for the inviscid rotating stratified Boussinesq equations in the stratification-dominant three-scale limit. J. Math. Pures Appl.  158(9): 90–119 (2022).

	\bibitem{NZ2024}
	B. Niu, W. Zhao, Improved stability threshold of the two-dimensional Couette flow for Navier-Stokes-Boussinesq systems via quasi-linearization.  Commun. Pure Appl. Anal, (2024).

	\bibitem{OK2024}
	W. Oxley, R. Kerswell, Linear stability of stratified, rotating, viscous plane Couette-Poiseuille flow. J. Fluid Mech. 991: Paper No. A16, 33 pp, (2024).

	\bibitem{Pedlosky1987}
	J. Pedlosky, Geophysical fluid dynamics, 2nd edn. Springer-Verlag, New York (1987).

	\bibitem{R1879}
	L. Rayleigh, On the stability, or instability, of certain fluid motions, Proc. Lond. Math. Soc., 11:  57–70 (1879/80).

	\bibitem{R1973}
	V. A. Romanov, Stability of plane-parallel Couette flow, Funkcional. Anal. i Prilozen, 7 (1973), pp.62–73.


\bibitem{SLF2025}
	J. Sun, N. Li, S. Fu, Global solutions of 3D Boussinesq equations for rotating stratified fluids. Z. Angew. Math. Phys. 76(3): Paper No. 112, 18 pp (2025).

	\bibitem{T1993}
	L. N Trefethen, A. E Trefethen, S. C Reddy and T. A Driscoll, Hydrodynamic stability without
	eigenvalues, Science, 261 (1993), 578–584.

	\bibitem{MR4373161}
	D. Wei, Z. Zhang, Transition threshold for the 3D Couette flow in Sobolev space. Comm. Pure Appl. Math. 74(11): 2398--2479  (2021).

	\bibitem{WZ2023}
	D. Wei, Z. Zhang, Nonlinear enhanced dissipation and inviscid damping for the 2D Couette flow. Tunis. J. Math. 5(3): 573–592 (2023).

	\bibitem{WZ2026}
	D. Wei and Z. Zhang, Asymptotic stability threshold of the 2D Couette flow in a finite channel, Ann. Inst. H. Poincar$\mathrm{\acute{e}}$  C Anal. Non Lin$\mathrm{\acute{e}}$aire 43 , no. 3, pp. 501-540 (2026).

	\bibitem{ZhaiZ2023}
	C. Zhai, W. Zhao, Stability threshold of the Couette flow for Navier-Stokes Boussinesq system with large Richardson Number $\gamma^2 > \frac{1}{4}$. SIAM J. Math. Anal. 55(2): 1284–1318 (2023).

	\bibitem{ZZ2023}
	Z. Zhang, R. Zi, Stability threshold of Couette flow for 2D Boussinesq equations in Sobolev spaces. J. Math. Pures Appl. 179(9): 123–182 (2023).
\end{thebibliography}
\end{document}